\documentclass[11pt]{amsart}
\usepackage[margin=0.9in]{geometry}
\usepackage{amsmath}
\usepackage{amsfonts}
\usepackage{amssymb}
\usepackage{graphicx}
\usepackage{hyperref}
\usepackage{mathrsfs}
\usepackage{graphicx,color}
\usepackage[mathscr]{eucal}
\usepackage{caption}
\usepackage{subcaption}
\usepackage{color}
\usepackage[all]{xy}
\usepackage{comment}
\usepackage{tikz-cd}
\usepackage{graphicx,color}
\usepackage{mathrsfs}
\usepackage{marvosym}
\usepackage{mathtools,tikz-cd}
\usepackage{xcolor}
\usepackage{enumitem}
\usepackage{graphicx}

\usepackage{tikz}
\usetikzlibrary{shapes.geometric, arrows}
\tikzstyle{startstop} = [rectangle, rounded corners, minimum width=2.5cm, minimum height=0.5cm,text centered, text width=2cm, draw=black, fill=white!30]
\tikzstyle{startstop2} = [rectangle, rounded corners, minimum width=2.5cm, minimum height=0.5cm,text centered, text width=2cm, draw=black, fill=white!30]
\tikzstyle{startstop3} = [rectangle, rounded corners, minimum width=2.5cm, minimum height=0.5cm,text centered, text width=2.5cm, draw=black, fill=white!30]
\tikzstyle{arrow} = [thick,->,>=stealth]

\newtheorem{theorem}{Theorem}[section]
\newtheorem{lemma}[theorem]{Lemma}

\newtheorem{corollary}[theorem]{Corollary}
\newtheorem{proposition}[theorem]{Proposition}
\newtheorem{definition}[theorem]{Definition}
\newtheorem{example}[theorem]{Example}
\newtheorem{remark}[theorem]{Remark}

\numberwithin{equation}{section}
\numberwithin{figure}{section}

\newcommand{\ZZ}{\mathbb{Z}}

\newcommand{\RR}{\mathbb{R}}

\newcommand{\PP}{\mathbb{P}}

\newcommand{\kg}{\mathfrak{g}}

\newcommand{\FF}{\mathbb{F}}

\newcommand{\CC}{\mathbb{C}}
\newcommand{\CP}{\mathbb{CP}}
\newcommand{\cS}{\mathcal{S}}
\newcommand{\cZ}{\mathcal{Z}}

\newcommand{\cO}{\mathcal{O}}

\newcommand{\cY}{\mathcal{Y}}

\newcommand{\cR}{\mathcal{R}}

\newcommand{\cC}{\mathcal{C}}

\newcommand{\cX}{\mathcal{X}}

\newcommand{\fkm}{\mathfrak{m}}

\newcommand{\Tot}{\mathit{Tot}}
\newcommand{\Hess}{\mathit{Hess}}

\newcommand{\diag}{\operatorname{diag}}

\newcommand{\AKh}{\operatorname{AKh}}
\newcommand{\Kh}{\operatorname{Kh}}

\def\eL{\EuScript{L}}

\newcommand{\fko}{\mathfrak{o}}
\newcommand{\frh}{\mathfrak{h}}
\newcommand{\frp}{\mathfrak{p}}
\newcommand{\frg}{\mathfrak{g}}
\newcommand{\frsl}{\mathfrak{sl}}
\newcommand{\frgl}{\mathfrak{gl}}

\newcommand{\HF}{\mathit{HF}}

\newcommand{\fix}{\mathrm{fix}}

\newcommand{\co}{\nobreak\mskip2mu\mathpunct{}\nonscript
  \mkern-\thinmuskip{:}\penalty300\mskip6muplus1mu\relax}

\usepackage{placeins}

\newcommand{\keh}{\textcolor{violet}}
\begin{document}

\title
{Symmetries of symplectic Khovanov--Rozansky homology}

\author{Akram Alishahi}
\address{Department of Mathematics\\ University of Georgia\\ Athens, GA 30602}
\email{akram.alishahi@uga.edu}
\thanks{AA was partially supported by NSF grant DMS-2238103}

\author{Kristen Hendricks}
\address{Department of Mathematics\\Rutgers University\\  New Brunswick, NJ 08854}
\email{kristen.hendricks@rutgers.edu}
\thanks{KH was partially supported by NSF grant DMS-2505573 and a Simons Travel Fellowship.}

\author{Cheuk Yu Mak}
\address{School of Mathematical and Physical Sciences\\ University of Sheffield\\ Sheffield S10 2TN, UK\\}
\email{c.mak@sheffield.ac.uk}
\thanks{CYM was partially supported by the Royal Society University Research Fellowship}

\begin{abstract}
We study the behavior of Manolescu's construction of symplectic Khovanov--Rozansky homology under symmetries. Specifically, for $k<n$ we construct a localization spectral sequence from the symplectic $\frsl(n)$ homology of a knot to the direct sum of the symplectic $\frsl(k)$ and $\frsl(n-k)$  homologies of the same knot, giving a geometric analog of rank-reducing spectral sequences arising from deformations on the non-symplectic theory. We also prove the existence of a localization spectral sequence relating the symplectic Khovanov--Rozansky homology of a two-periodic link to the symplectic Khovanov--Rozansky homology of the quotient link, in analogy to a theorem of Seidel and Smith for symplectic Khovanov homology. The proofs rely on a generalization of the localization theorem of Seidel--Smith and Large for equivariant Lagrangian Floer cohomology with respect to order two actions; in particular, we extend the theory to the case that the fixed set of the action carries a stable normal real polarization.
\end{abstract}

\maketitle

\tableofcontents

\section{Introduction} Khovanov--Rozansky homology is an invariant of links which associates to a link $L$ and a field $\mathbb F$ a bigraded homology theory $H_n(L; \mathbb F)$\footnote{The theory was originally defined over $\mathbb Q$, but has since been extended to arbitrary coefficients \cite{Kh:sl3, QR:combinatorial}.}. For every rank $n\geq 2$, the bigraded Euler characteristic of the theory is the $\frsl(n)$ polynomial of the link \cite{KR:Part1, KR:Part2}. The theory generalizes the $n=2$ case, which is Khovanov's link homology theory, and has Euler characteristic the Jones polynomial \cite{Khovanov:Jones}.

Khovanov--Rozansky homology admits algebraically-defined rank-reducing spectral sequences which generalize the well-known Lee spectral sequence for Khovanov homology \cite{Leedeformation, Rasmussen:Milnor}, as follows. For simplicity let $K$ be a knot. Beginning by focusing on complex coefficients, the underlying Frobenius algebra of Khovanov--Rozansky homology is $\mathbb C[x]/(x^n)$. If one replaces $x^n$ with a monic degree $n$ polynomial $p(x)$ having distinct roots $(r_1,\dots, r_s)$ of multiplicities $(m_1, \dots, m_s)$, the deformed theory splits into a direct sum $\bigoplus_{i=1}^s H_{m_i}(K;\mathbb C)$; see for example \cite[Section 4]{GW:Kauffman} or \cite[Section 1.1]{RW:foams}. One obtains a deformation spectral sequence
\begin{equation} \label{eq:krssbig} H_n(K;\mathbb C) \rightrightarrows \bigoplus_{i=1}^s H_{m_i}(K;\mathbb C). \end{equation}
Some special cases of this spectral sequence merit particular remark. In the case $s=n$, such that $p(x)$ has $n$ distinct roots, in analogy to the Lee spectral sequence Gornik \cite{Gornik} and Wu \cite{Wu:generic} show that this spectral sequence converges to $\mathbb C^{\oplus n}$. These deformations have been used to give new slice genus bounds by Lobb and Lewark \cite{Lobb:slice, LL:quantum}. On the other hand, in the case $s=2$, one has spectral sequences associated to each $k$ such that $n > k >0$ of the form
\begin{equation} \label{eq:krsssplit} H_n(K; \mathbb C) \rightrightarrows H_{k}(K; \mathbb C) \oplus H_{n-k}(K; \mathbb C) \end{equation}
governed by the deformation $p(x)=x^k(x-1)^{n-k}$. There is furthermore a similar story over fields $\mathbb F$ other than $\mathbb C$ for polynomials which split over $\mathbb F$ into a product of linear terms. For example, one obtains an analog of \eqref{eq:krsssplit} over $\mathbb F_2$, which in the case that $n=2$ and $k=1$ is the Bar-Natan spectral sequence \cite{BN:tangles}. However, over $\mathbb F_2$ one cannot obtain any analog of \eqref{eq:krssbig} with $s\geq 3$, since the maximum number of distinct roots of a polynomial over $\mathbb F_2$ is two. Rank-reducing spectral sequences have also been considered more broadly in colored versions of $\frsl(n)$ homology; see for example \cite{RW:foams, Wedrich:exponential}. 

Khovanov homology has a Lagrangian Floer theoretic analog defined by Seidel and Smith \cite{SS06,SS10} and often studied using a reformulation of Manolescu \cite{Manolescu:nilpotent}. This theory, called symplectic Khovanov homology, is known to be equivalent to the original theory in characteristic zero after a collapse of the bigrading \cite{AS16, AS19}. The construction of symplectic Khovanov homology is significantly more demanding than that of combinatorial Khovanov homology. However, because of the more intrinsically geometric nature of the symplectic theory, it has proved to be fertile ground for studying group actions on Khovanov homology and realizing algebraic constructions geometrically \cite{SS10,HLS:flexible, HMR}. Inspired by symplectic Khovanov homology, Manolescu \cite{Manolescu:sln} has also constructed a symplectic $\frsl(n)$ homology theory $H_{n,symp}(L)$, which is a singly-graded theory conjecturally equivalent to combinatorial $\frsl(n)$ homology after a collapse of the bigrading. This theory associates to a link the Lagrangian Floer cohomology of a triple consisting of an exact symplectic manifold $\mathcal Y_{m,n,\tau}$ together with two exact Lagrangians $\eL_0$ and $\eL_1$, each of which is diffeomorphic to a product of $m$ copies of $\mathbb{CP}^{n-1}$. (Here $m$ is the index of a choice of braid representative of the link.) The theory is definable over arbitrary coefficients, but we will only study it over $\mathbb F_2$.

In this paper we use equivariant Lagrangian Floer theory to produce a geometrically-defined analog of the spectral sequence \eqref{eq:krsssplit} for Manolescu's symplectic $\frsl(n)$ homology. Our main result is the following. Let $\mathbb F_2[\theta, \theta^{-1}] = \theta^{-1}H^*(\mathbb{RP}^{\infty}; \mathbb F_2)$ denote the cohomology of the classifying space $B\ZZ/ 2\ZZ \simeq \mathbb{RP}^{\infty}$ localized at the degree one generator.

\begin{theorem} \label{thm:mainss} Given a knot $K \subset S^3$ and $n,k$ such that $0< k < n$, there is a spectral sequence with $E^1$ page
\[ H_{n, symp}(K) \otimes \FF_2[\theta, \theta^{-1}]\]
and $E^{\infty}$ page isomorphic to
\[ (H_{k,symp}(K) \oplus H_{n-k, symp}(K)) \otimes \FF_2[\theta, \theta^{-1}].\]
\end{theorem}

Note that in the above $H_{1,symp}(K) \simeq \FF_2$. Since our present tools require working over the field with two elements, this is the strongest statement for a knot that one could currently hope to obtain; see Remark~\ref{rem:orderp} for further discussion. Our techniques also obtain a spectral sequence for a link $L$, which appears to parallel existing results for links on the combinatorial side \cite{GW:Kauffman, RW:foams, Wedrich:exponential}. However,  understanding the final page of the spectral sequence in the link case is more complicated than in the knot case, as we discuss in Section~\ref{sec:link-discussion}.

We now discuss how Theorem \ref{thm:mainss} arises. The key technical input of the theorem is a mild generalization of a localization theorem of P. Seidel and I. Smith \cite{SS10} and T. Large \cite{Large}. The basic setup of these theorems and our extension is as follows. Let $M$ be an exact symplectic manifold, convex at infinity, containing two exact compact Lagrangians $\eL_0$ and $\eL_1$, and let $\iota$ be a symplectic involution on $M$ preserving each $\eL_i$ setwise for $i=0,1$. Let $M^{\fix}$ be the fixed set of $M$ under $\iota$. Suppose that $M^{\fix}$ has finitely many connected components, and let $\{M_{j}^{\fix}: j\in A\}$ be the set of connected components of $M^{\fix}$. Each connected component $M_{j}^{\fix}$ is an exact symplectic manifold with two (possibly empty) exact Lagrangians $\eL_{i,j}^{\fix}$ consisting of the intersections of the fixed sets $\eL_i^{\fix}$ of the Lagrangians with $M_{j}^{\fix}$. To each component of the fixed set we may associate the triple of normal bundles $(N M_{j}^{\fix}, N\eL_{0,j}^{\fix}, N\eL_{1,j}^{\fix})$, which we refer to as the normal polarization $\mathfrak{p}_{N,j}$. (We discuss polarizations at greater length below.) Under suitable hypotheses on the behavior of these normal polarizations, one obtains a localization spectral sequence
\begin{equation} \label{eq:general-loc} \HF(M, \eL_0, \eL_1) \otimes \mathbb F_2[\theta, \theta^{-1}] \rightrightarrows \bigoplus_{j \in A} \HF(M_{j}^{\fix}, \eL_{0,j}^{\fix}, \eL_{1,j}^{\fix}) \otimes \mathbb F_2[\theta, \theta^{-1}].\end{equation}
In the case of symplectic $\frsl(n)$ homology, given $k$ such that $0<k < n$, we consider an involution $\iota_k$ on $\mathcal Y_{m,n,\tau}$ which preserves the Lagrangians. While we defer a full definition of $\iota_k$ to Section~\ref{sec:intrinsic}, it is helpful to observe its behavior on the Lagrangians $\eL_i$ for $i=0,1$. Each of these Lagrangians is up to isotopy a product $(\mathbb{CP}^{n-1})^{\times m}$, and the action of $\iota_k$ on each factor of $\mathbb{CP}^{n-1}$ in the product is
\[ [z_0 : \dots : z_{n-1}] \mapsto [-z_0: \dots : -z_{k-1} : z_{k} : \dots : z_{n-1}].\]

Notice that negating the first $k$ coordinates is the same as negating the final $n-k$ coordinates. The fixed set of each factor $\mathbb{CP}^{n-1}$ is the union of a copy of $\mathbb{CP}^{k-1}$ and a copy of $\mathbb{CP}^{n-k-1}$. The fixed set of the product Lagrangian therefore has $ m \choose t$ components of real dimension $2( t(k-1) + (m-t)(n-k-1))$ for each $0\leq t \leq m$. However, in Section~\ref{sec:lagrangian-involution} we will show that only two components of the fixed set $\cY_{m,n, \tau}^{\fix}$ contain a component of the fixed set of both Lagrangians. If $k \neq n-k$, these are the unique components of highest and lowest dimension. In one of these two special components, the fixed set is precisely an embedded copy of $\cY_{m,k,\tau}$, with the Lagrangians embedded so as to compute $H_{k, symp}(K)$. Similarly in the other special component, the fixed set is precisely an embedded copy of $\mathcal Y_{m,n-k,\tau}$ with fixed Lagrangians which compute $H_{n-k, symp}(K)$. 

In order to obtain a localization theorem of the form~\eqref{eq:general-loc}, such as Theorem~\ref{thm:mainss}, one must (1) compute the normal polarization class of the fixed set of the involution $\iota_k$, and (2) show that the result satisfies some hypothesis under which localization holds. Most of the work of this paper is devoted to (1). In particular, we prove the following result (Theorem \ref{thm:stable_real}), which is of independent interest.

Recall that if $(M, \eL_0, \eL_1)$ is a triple consisting of a symplectic manifold and two Lagrangians, then a polarization $\mathfrak{p} = (E, F_0, F_1)$ is a symplectic vector bundle $E$ over $M$ together with $F_i$ a Lagrangian subbundle of $E|_{\eL_i}$ for $i=0,1$. We say that $E$ has a real structure if it admits a global Lagrangian subbundle $E_{\RR}$ and that $\mathfrak{p}$ admits a real polarization if the restriction of $E_{\RR}$ to $\eL_i$ is isotopic to $F_i$ as Lagrangian subbundles for $i=0,1$. We say that $E$ admits a stable real structure if $E\oplus \underline{\CC}^N$ admits a real structure for some $N\geq 0$ and that $\mathfrak{p}$ admits a stable real polarization if the stabilized polarization $(E\oplus\underline{\CC}^N, F_0 \oplus\underline{\RR}^N, F_1 \oplus i\underline{\RR}^N)$ admits a real polarization for some $N\geq 0$. We review these definitions at greater length in Section~\ref{sec:localization}.

\begin{theorem}\label{thm:stable_real}
The tangent bundle $T\mathcal Y_{m,n,\tau}$ of the exact symplectic manifold $\mathcal Y_{m,n,\tau}$ admits a stable real structure $T_{\RR}\mathcal Y_{m,n,\tau}$. Moreover, this real structure gives a stable real polarization of the tangent polarization $\mathfrak{p}_T = (T\mathcal Y_{m,n,\tau}, T\eL_0, T\eL_1)$. If $M^{\fix}_{j}$ is a connected component of the fixed set of the involution $\iota_k$ containing connected components $\eL_{0,j}^{\fix}$ and $\eL_{1,j}^{\fix}$ of the fixed Lagrangians, there is likewise a stable real polarization of the normal polarization $\mathfrak{p}_{N, j} = (NM_{j}^{\fix}, N\eL_{0,j}^{\fix}, N\eL_{1, j}^\fix)$.
\end{theorem}

Seidel and Smith's localization theorem originally required that the normal polarization to the fixed set of the involution be stably trivial in an appropriate sense; Large extended the theorem to the case that there existed a stable isomorphism of polarization classes between the normal polarization in each component and its tangent polarization $(T M_{j}^{\fix}, T\eL_{0,j}^{\fix}, T\eL_{1,j}^{\fix})$. Unfortunately, an elementary characteristic class computation shows that neither of these assumptions holds in our case.
Therefore, in order to prove Theorem~\ref{thm:mainss}, we supply the following extension.

\begin{theorem}\label{thm:localization-intro}
Suppose that $(M, \eL_0, \eL_1)$ and $\iota$ satisfy the symplectic hypotheses above. Further suppose that $M^{\fix}$ has finitely many connected components $M^{\fix}_{j}$ for $j \in A$, and suppose that for each component $(M^{\fix}_{j}, \eL_{0,j}^{\fix}, \eL_{1,j}^{\fix})$ in which both subsets $\eL_{i,j}^{\fix} = M^{\fix}_{j} \cap \eL_i^{\fix}$ are nonempty, the normal polarization $\mathfrak{p}_{N,j}=(NM^{\fix}_{j}, N\eL_{0,j}^{\fix}, N\eL_{1,j}^{\fix})$ admits a stable real polarization. Then there is a spectral sequence with $E^1$ page
\[ \HF(M, \eL_0, \eL_1) \otimes \FF_2[\theta, \theta^{-1}]\]
and $E^{\infty}$ page isomorphic to 
\[\bigoplus_{j\in A} (\HF(M_{j}^{\fix}, \eL_{j,0}^{\fix}, \eL_{j,1}^{\fix})) \otimes \FF_2[\theta, \theta^{-1}].\]
\end{theorem}

We refer to the hypothesis above as $M^{\fix}$ being equipped with a set of stable normal real polarizations. We note that this extension of the localization theorem was expected to hold by the authors of the earlier versions of the Theorem \cite[Section 3.3]{Large}. The proof, which appears in Section~\ref{subsec:localization-realpolar} as the proof of Theorem~\ref{t:genSS} and Corollary~\ref{cor:components}, closely follows  the strategy of Seidel and Smith's original work~\cite[Section 3]{SS10}.

We also study extrinsic symmetries of symplectic Khovanov--Rozansky homology, which is to say symmetries arising from topological actions on the original link. In particular, recall that a link $L \subset S^3$ is said to be two-periodic if it is preserved setwise by an orientation-preserving involution on $S^3$ which has fixed set an unknot disjoint from the link. Taking the quotient by the action produces the quotient link $\overline{L}$ in $S^3$. Seidel and Smith proved there is a spectral sequence relating the symplectic Khovanov homology of the periodic link and that of its quotient \cite[Section 5]{SS10}, to wit
\begin{equation} \label{eq:ssperiodic} \Kh_{symp}(L; \FF_2) \otimes \FF_2[\theta, \theta^{-1}] \rightrightarrows \Kh_{symp}(\overline{L}; \FF_2)\otimes \FF_2[\theta, \theta^{-1}]. 
\end{equation}

We lift the action on a two-periodic link $L$ to the space $\mathcal Y_{m,n,\tau}$, and obtain the following analog.
\begin{theorem} \label{thm:periodicss} Let $L$ be a two-periodic link, and $\overline{L}$ its quotient link. There is a spectral sequence with $E_1$ page 
\[ H_{n, symp}(L) \otimes \FF_2[\theta, \theta^{-1}]\]
and $E_{\infty}$ page isomorphic to 
\[ H_{n, symp}(\overline{L}) \otimes \FF_2[\theta, \theta^{-1}].\]
\end{theorem}

Although a proof of Theorem~\ref{thm:periodicss} requires a careful analysis of the tangent and normal polarizations of the fixed set of our lift of the involution, it does not require the full strength of Theorem~\ref{thm:localization-intro}; in particular, it fulfills the hypotheses of Large's localization theorem \cite{Large}.

\begin{remark} The expert reader may be aware that the spectral sequence \eqref{eq:ssperiodic} is the first in a family of localization results for symplectic Khovanov homology. Recent work of the second and third authors together with S. Raghunath shows that there is also a spectral sequence
\begin{equation} \label{eqn:hmrakh} \Kh_{symp}(L; \FF_2) \otimes \FF_2[\theta, \theta^{-1}] \rightrightarrows \AKh_{symp}(\overline{L}; \FF_2) \otimes \FF_2[\theta, \theta^{-1}] \end{equation}
where $\AKh_{symp}$ is a construction of annular symplectic Khovanov homology \cite[Theorem 1.5]{HMR}. This spectral sequence is the analog of a spectral sequence constructed for combinatorial Khovanov homology by M. Stoffregen and M. Zhang \cite{SZ:localization} and by M.\,Borodzik, W.\,Politarczyk, and M.\,Silvero \cite{BPS}. However, the proof of \eqref{eqn:hmrakh} relies heavily on the fact that symplectic Khovanov homology has a reformulation due to Manolescu \cite{Manolescu:nilpotent} in terms of the Hilbert scheme of a Milnor fibre; as far as the authors are aware, no such description of symplectic $\frsl(n)$ homology exists. It is therefore presently unclear to the authors how to obtain an analog of the spectral sequence \eqref{eqn:hmrakh}.
\end{remark}

\begin{remark} \label{rem:Khovanov-case} Let us briefly discuss the special case of Theorem~\ref{thm:mainss} for $n=2$ and $k=1$. In this case, the spectral sequence of Theorem~\ref{thm:mainss} specializes to
\[ \Kh_{symp}(K) \otimes \FF_2[\theta, \theta^{-1}] \rightrightarrows \FF_2^{\oplus 2} \otimes \FF_2[\theta, \theta^{-1}].\]
Indeed, the arguments of Section~\ref{sec:link-discussion} show that for a link $L$ having $\ell$ components, one has a spectral sequence
\begin{equation} \label{eq:symplecticbn} \Kh_{symp}(L) \otimes \FF_2[\theta, \theta^{-1}] \rightrightarrows \FF_2^{\oplus 2^{\ell}} \otimes \FF_2[\theta, \theta^{-1}].\end{equation}
This spectral sequence and the action inducing it were previously studied in work of the second author, Lipshitz, and Sarkar \cite[Theorem 7.16]{HLS:flexible}. The fixed set of the symplectic manifold under the relevant involution is in fact zero-dimensional, and the Lagrangians intersect in precisely $2^{\ell}$ points, which may be understood to correspond to the possible orientations of the link. The authors conjectured \cite[Conjecture 7.28]{HLS:flexible} that the spectral sequence \eqref{eq:symplecticbn} is the analog of the Bar-Natan spectral sequence associated to the deformation $p(x) = x(x-1)$ in combinatorial Khovanov homology \cite{BN:tangles}. The appearance of the same spectral sequence within a generalized set of rank-reducing spectral sequences for $\frsl(n)$ homology seems to support the conjecture.

We note that the same set of $2^{\ell}$ generators is also the fixed set of the symplectic manifold computing symplectic Khovanov homology under a natural intrinsic action of $S^1$ (from which the involution above may be recovered by restricting to an order two subgroup). It has occasionally been proposed that there may be a connection between the (localized) $S^1$-equivariant symplectic Khovanov homology, partially defined in \cite{HLS:simplicial}, and the Lee homology; we comment further on this in the next remark.
\end{remark}

\begin{remark} \label{rem:orderp} The extension of existing localization theory for Floer cohomology given in this paper is not the only new such extension: recently E.\,Bao, T.\,Lawson, and the third author extended Seidel and Smith's original work to the case of actions of order $p$, for $p$ an odd prime, on Lagrangian Floer cohomology in the presence of a suitably equivariant stable normal trivialization on the fixed set \cite[Theorem 1.1]{BLM}. 
The authors apply their theorem to extend Seidel and Smith's spectral sequence \eqref{eq:ssperiodic} for two-periodic links to $p$-periodic links \cite[Theorem 1.2 and Section 3]{BLM}.

The cannot immediately apply the $\mathbb Z/p\mathbb Z$ localization theorem of \cite{BLM} to extend either Theorem~\ref{thm:mainss} or Theorem~\ref{thm:periodicss}, since in neither case is there a stable normal trivialization on the fixed set; further technical development is therefore required in order to combine these ideas. However, we comment briefly on what one may expect. The construction of the intrinsic involution and analysis of its fixed sets in Section~\ref{sec:intrinsic} generalizes to the order $p$ case: for every choice partition of $n$ into the sum of $s$ integers $k_s$ with $2 \leq s \leq p$, one obtains an order $p$ action on $(\cY_{m,n,\tau}, \eL_0, \eL_1)$ with fixed sets whose Floer homology computes a direct sum $\oplus_{i=1}^s H_{k_s, symp}(K, \FF_p)$. We expect this action to give rise to the mod $p$ analog of the combinatorial sequence \eqref{eq:krssbig}. (One would also expect a straightforward analog of Theorem~\ref{thm:periodicss} for $p$-periodic links.)

More speculatively, one may also modify the action of Section~\ref{sec:intrinsic} to be an $S^1$-action with the same fixed set, and similarly for the order $p$ actions of the previous paragraph. We conjecture that $S^1$-localization with respect to these actions produces the analog of the characteristic zero spectral sequences \eqref{eq:krssbig}, in analogy with the conjecture for the $n=2$ case mentioned in Remark~\ref{eq:symplecticbn}.

Another recent extension of $\ZZ/2\ZZ$ localization may be found in recent work of S.\,Bai, E.\,Shelukhin, Y.\,Wang, and G.\,Xu on the Arnol'd-Givental conjecture, which contains new techniques specifically focused on the case of Hamiltonian diffeomorphisms \cite{BSWX:arnold}.
\end{remark}

\subsection*{Conventions} We record here some notational and expositional conventions for the benefit of the reader.

\begin{itemize}

\item Throughout, \emph{combinatorial} and \emph{symplectic} are antonyms; in particular combinatorial $\frsl(n)$ homology indicates Khovanov and Rozansky's original invariant, with suitable coefficients.

\item Throughout, $L$ is a link and $\eL$ is a Lagrangian.

\item The \emph{rank} of $\frsl(n)$ homology is $n$. We will mostly not need to refer to the dimension of the homology as a vector space.

\item We use the cohomological construction of Lagrangian Floer cohomology; in particular, given a triple $(M, \eL_0, \eL_1)$, a topological map $\phi \co \mathbb R \times [0,1]\to M$ such that $\phi(s,0) \in \eL_0$, $\phi(s,1) \in \eL_1$, and furthermore such that $\lim_{s \rightarrow \pm \infty} \phi(s,t) = x_{\pm}$ is understood to run from $x_+$ to $x_-$. This is the convention ordinarily used in symplectic Khovanov homology (and the opposite of the convention used in Heegaard Floer theory). 
\end{itemize}

\subsection*{Organization} This paper is organized as follows. In Section~\ref{sec:localization} we discuss equivariant Lagrangian Floer theory for symplectic involutions; in particular, in Section~\ref{subsec:localization-existing}  we review the localization theorems of Seidel--Smith and Large for Lagrangian Floer cohomology and in Section~\ref{subsec:localization-realpolar} we prove a new refinement in the presence of a stable normal real polarization. 

In Section~\ref{sec:construction} we review Manolescu's construction of symplectic $\frsl(n)$ homology. In particular, in Section~\ref{sec:lie} we review preliminary constructions from Lie theory, which we then use in Section~\ref{sec:ambient} to introduce the ambient symplectic manifold in which symplectic Khovanov--Rozansky homology is computed. In Section~\ref{sec:parallel} we review the notion of rescaled parallel transport on the (singular) fibration whose fibre is our ambient symplectic manifold $\mathcal Y_{m,n,\tau}$, and in Section~\ref{sec:lagrangian} we use this construction to associate a Lagrangian in $\mathcal Y_{m,n,\tau}$ to any crossingless matching of $2m$ points. Finally in Section~\ref{sec:defn} we use this construction to define the symplectic Khovanov--Rozansky homology. 

Section~\ref{sec:intrinsic} is devoted to the study of an intrinsic involution on this theory and the proof of Theorem~\ref{thm:mainss}. In Section~\ref{sec:ambient-intrinsic} we define the action on the ambient symplectic manifold and study the connected components of its fixed set. In Section~\ref{sec:lagrangian-involution} we confirm that our action preserves the two Lagrangians and locate the components of the fixed Lagrangians within the fixed components of the ambient symplectic manifold. Section~\ref{sec:real} is devoted to the proof of Theorem~\ref{thm:stable_real}. Finally, in Section~\ref{sec:assemble-proof-thmmain} we assemble the elements of the previous sections into a proof of Theorem~\ref{thm:mainss}. We conclude the section with some discussion of the link case and future research in Section~\ref{sec:link-discussion}.

In Section~\ref{sec:extrinsic} we study an extrinsic involution on the theory which lifts the topological action on a two-periodic link; in particular, after reviewing some basic facts about link periodicity in Section~\ref{sec:periodic-action} we define our lift of the periodic action, and in Section~\ref{sec:periodic-tniso} we study the fixed sets of the symplectic manifold and Lagrangians under this action and prove the existence of a stable tangent normal isomorphism in the sense of Large. This allows us to supply a proof of Theorem~\ref{thm:periodicss} in Section~\ref{sec:proof-periodic}.

\subsection*{Acknowledgments} We are grateful to Robert Lipshitz, Ivan Smith, Linh Truong, Jake Rasmussen, Sriram Raghunath, Paul Wedrich, and Melissa Zhang for helpful conversations. Portions of this work were carried out at a Simons Center Workshop on Gauge Theory and Floer Homology in Low Dimensional Topology in April 2025, the 2025 Georgia International Topology Conference in May 2025, and the Rutgers Symplectic Summer School in August 2025. We are grateful to all three conferences for their hospitality. Portions of this work were also carried out while the second author was in residence at the ``New frontiers in gauge theory, topology, and physics'' Thematic Visiting Program at the Okinawa Institute of Science and Technology in Fall 2026; she is grateful to the institute and the program organizers.

\subsection*{Declaration of AI use} 

The infinite loop space argument in the second half of the proof of Lemma \ref{lem:isomorphism-is-isotopy} was suggested by ChatGPT-5.6 Luna and written in its present form by the authors.  The proof that Equation~\eqref{eqn:gemini-helped} holds was assisted by consultation with Google Gemini 3.6 Thinking. Claude Opus 5.5 was used for final proofreading of the paper. AI was not involved in the design of other proofs or in the writing of the paper. 

\section{Localization in equivariant Floer theory} \label{sec:localization}

In this section we review some of the ideas surrounding localization spectral sequences in equivariant Floer theory, and prove a mild extension of previous results.

\subsection{Background and existing results} \label{subsec:localization-existing} We now review the localization spectral sequence of Seidel--Smith and Large. We begin by recalling the following two definitions, stated here as in \cite[Section 3.2]{Large}.

\begin{definition} Let $M$ be a symplectic manifold containing Lagrangians $\eL_0$ and $\eL_1$. A \emph{set of polarization data} or \emph{polarization} for $(M, \eL_0, \eL_1)$ is a triple $\mathfrak p = (E, F_0, F_1)$ such that
\begin{itemize}
\item $E$ is a symplectic vector bundle over $M$, and
\item $F_i$ is a Lagrangian subbundle of $E|_{\eL_i}$ for $i=0,1$.
\end{itemize}
\end{definition}

As an example, one always has the tangent polarization $\mathfrak{p}_T=(TM, T\eL_0, T\eL_1)$. Given $\mathfrak p = (E, F_0, F_1)$ a set of polarization data for $(M, \eL_0, \eL_1)$, one may stabilize by a trivial bundle to obtain $\mathfrak p \oplus \underline{\CC}^k = (E \oplus  \underline{\CC}^k, F_0 \oplus \underline{\RR}^k, F_1 \oplus i \underline{\RR}^k)$.

\begin{definition} Let $\mathfrak p = (E, F_0, F_1)$ and $\mathfrak p' = (E', F_0', F_1')$ be two sets of polarization data for $(M, \eL_0, \eL_1)$. An \emph{isomorphism of polarization data} between $\mathfrak p$ and $\mathfrak p'$ is an isomorphism of symplectic vector bundles 
\[\psi \co E \rightarrow E'\]
such that there are homotopies through Lagrangian subbundles of $(E')|_{\eL_i}$ between $\psi(\eL_i)$ and $\psi(F_i)$ and $F_i'$ for $i=0,1$. A \emph{stable isomorphism of polarization data} between $\mathfrak p$ and $\mathfrak p'$ is an isomorphism of polarization data between $\mathfrak p \oplus \underline{\CC}^{k_1}$ and $\mathfrak p' \oplus \underline{\CC}^{k_2}$ for some $N_1$ and $N_2$.\end{definition}

We now turn our attention to the situation of involutions. Let $M$ be an exact symplectic manifold which is convex at infinity, containing two compact exact Lagrangians $\eL_0$ and $\eL_1$, which in the usual way may always be made to intersect transversely via some Hamiltonian isotopy. Suppose that $\iota$ is a symplectic involution on $M$, that is, that $\iota^2 = \mathrm{Id}$ and $\iota^*\omega = \omega$. Further let us suppose that $\iota$ preserves each of $\eL_0$ and $\eL_1$ setwise. Suppose that the fixed sets $(M^{\fix}, \eL_0^{\fix}, \eL_1^{\fix})$ under $\iota$ decompose into connected components $\bigcup_{j\in A} (M_j^{\fix}, \eL_{0,j}^{\fix}, \eL_{1,j}^{\fix})$, such that for each value of $j$ we have that $(M_j^{\fix}, \eL_{0,j}^{\fix}, \eL_{1,j}^{\fix})$ is a new triple consisting of an exact symplectic manifold with two exact compact Lagrangians. We refer to this set of assumptions as the symplectic hypotheses of the theory.

\begin{definition} \label{def:stable-normal} Let $(M, \eL_0, \eL_1)$ and $\iota$ satisfy the symplectic hypotheses above, and suppose that the fixed set decomposes into connected components $(M^{\fix}, \eL_0^{\fix}, \eL_1^{\fix}) = \bigcup_{j\in A} (M_j^{\fix}, \eL_{0,j}^{\fix}, \eL_{1,j}^{\fix})$. 
\begin{itemize}
\item A \emph{stable normal trivialization} for the fixed set consists of a stable isomorphism of polarization data between the normal polarization $\mathfrak p_{N_j} =(NM_j^{\fix}, N\eL_{0,j}^{\fix}, N\eL_{1,j}^{\fix})$ and the trivial polarization $\underline{\mathbb C}$ for each  $j \in A$. 
\item A \emph{stable tangent-normal isomorphism} for the fixed set consists of a stable isomorphism of polarization data between the normal polarization $\mathfrak p_{N_j} =(NM_j^{\fix}, N\eL_{0,j}^{\fix}, N\eL_{1,j}^{\fix})$ and tangent polarization $\mathfrak p_{T_j} = (TM_j^{\fix}, T\eL_{0,j}^{\fix}, T\eL_{1,j}^{\fix})$ for each $j \in A$. \end{itemize} \end{definition}

Before moving on, we recall to the reader's attention that since the symplectic group deformation retracts onto the unitary group, the theory of symplectic vector bundles is isomorphic to the theory of complex vector bundles. Therefore, the definitions above may be equivalently phrased in terms of complex bundles with totally real subbundles.

The original localization theorem of Seidel and Smith for order two actions on Lagrangian Floer cohomology is as follows.

\begin{theorem} \cite[Theorem 20]{SS10} \label{thm:ss-localization} Let $(M, \eL_0, \eL_1)$ and $\iota$ satisfy the symplectic hypotheses, and suppose that $(M^{\fix}, \eL_0^{\fix}, \eL_1^{\fix})$ admits a stable normal trivialization. Then there is a spectral sequence whose $E^1$ page is
\[ \HF(M, \eL_0, \eL_1) \otimes \FF_2[\theta, \theta^{-1}]\]
and whose $E^{\infty}$ page is isomorphic to 
\[\bigoplus_{j\in A} \HF(M_j^{\fix}, \eL_{0,j}^{\fix}, \eL_{1,j}^{\fix}) \otimes \FF_2[\theta, \theta^{-1}].\]
\end{theorem}

\noindent Large extends the hypotheses of this theorem as follows.

\begin{theorem} \cite[Sections 1 and 10]{Large} \label{thm:large-localization} Let $(M, \eL_0, \eL_1)$ and $\iota$ satisfy the symplectic hypotheses, and suppose that $(M^{\fix}, \eL_0^{\fix}, \eL_1^{\fix})$ admits a stable tangent-normal isomorphism. Then there is a spectral sequence whose $E^1$ page is
\[ \HF(M, \eL_0, \eL_1) \otimes \FF_2[\theta, \theta^{-1}]\]
and whose $E^{\infty}$ page is isomorphic to 
\[\bigoplus_{j\in A}\HF(M_j^{\fix}, \eL_{0,j}^{\fix}, \eL_{1,j}^{\fix}) \otimes \FF_2[\theta, \theta^{-1}].\]
\end{theorem}

\begin{remark} Previous applications of the localization theorems of Seidel--Smith and Large have most commonly focused on the case that the fixed set is connected. However, in neither case does the argument assume this fact. Both proofs involve arguments local to the fixed set and go through without incident for the case that the fixed set is disconnected and has finitely many components of possibly differing dimension. \end{remark}

\subsection{Localization in the presence of a stable normal real polarization} \label{subsec:localization-realpolar}

In this section, we generalize the localization theorem of Seidel--Smith and Large to the case that the normal bundle $NM^\fix$ admits a stable real polarization. In particular, we prove Theorem~\ref{thm:localization-intro}, restated and slightly expanded in the form of Theorem~\ref{t:genSS} and 
Corollary~\ref{cor:components} below.

We continue the symplectic hypotheses and notation of the previous section: $M$ is an exact symplectic manifold which is convex at infinity, containing two transversely-intersecting compact exact Lagrangians $\eL_0$ and $\eL_1$. The action $\iota$ is a symplectic involution preserving each of $\eL_i$ for $i=0,1$ setwise. The fixed sets $(M^{\fix}, \eL_0^{\fix}, \eL_1^{\fix})$ under $\iota$ decompose into connected components $\bigcup_{j\in A} (M_j^{\fix}, \eL_{0,j}^{\fix}, \eL_{1,j}^{\fix})$. 

We begin by reviewing the definitions of a real structure and a stable real polarization from the introduction.

\begin{definition} Let $\mathfrak{p} = (E,F_0,F_1)$ be a set of polarization data for $(M, \eL_0, \eL_1)$. The bundle $E$ has a real structure if it admits a global Lagrangian subbundle $E_{\RR}$. The set of polarization data $\mathfrak{p}$ admits a real polarization if $E$ admits a real structure $E_{\RR}$ and the restriction of $E_{\RR}$ to $\eL_i$ is isotopic to $F_i$ as Lagrangian subbundles for $i=0,1$, or equivalently if the polarization $(E, E_{\RR}|_{\eL_0}, E_{\RR}|_{\eL_1})$ is isomorphic to $\mathfrak{p}$.  Finally, $E$ admits a stable real structure if $E \oplus \underline{\CC}^k$ admits a real structure for some $k\geq 0$, and $\mathfrak{p}$ admits a stable real polarization if $\mathfrak{p}\oplus \underline{\CC}^k$ admits a real polarization for some $k\geq 0$.
\end{definition}

\begin{definition}\label{def:sR-polarization} Let $(M, \eL_0, \eL_1)$ be a triple together with an involution $\iota$ satisfying the symplectic hypotheses above. 
\begin{itemize} \item If $M^{\fix}$ has a single connected component, then $M^{\fix}$ is said to be equipped with a stable normal real polarization if $\mathfrak{p}_N = (NM^{\fix}, N\eL_0^{\fix}, N\eL_1^{\fix})$ admits a stable real polarization.

\item If $M^{\fix}$ has multiple connected components $M_j^{\fix}$ for $j \in A$, then $M^{\fix}$ is said to be equipped with a set of stable normal real polarizations if every component $M_j^{\fix}$ containing a nonempty component of both $\eL_0^{\fix}$ and $\eL_1^{\fix}$ is equipped with a stable normal real polarization.
\end{itemize}
\end{definition}

We now state the most straightforward form of our new localization theorem.

\begin{theorem}\label{t:genSS}
Suppose that $(M, \eL_0, \eL_1)$ satisfies the symplectic hypotheses above, and $M^{\fix}$ has a single connected component. Suppose that $M^{\fix}$ is equipped with a stable normal real polarization. Then there is a spectral sequence whose $E^1$ page is
\[ \HF(M, \eL_0, \eL_1) \otimes \FF_2[\theta, \theta^{-1}]\]
and whose $E^{\infty}$ page is isomorphic to 
\[\HF(M^{\fix}, \eL_0^{\fix}, \eL_1^{\fix}) \otimes \FF_2[\theta, \theta^{-1}].\]
\end{theorem}

Because our argument, like that of Seidel--Smith and Large, is local to the fixed set, the case for multiple components will be an immediate corollary.

\begin{corollary}\label{cor:components}
Suppose that $(M, \eL_0, \eL_1)$ satisfies the symplectic hypotheses above, and $M^{\fix}$ has finitely many connected components and is equipped with a set of stable normal real polarizations. Then there is a spectral sequence whose $E^1$ page is
\[ \HF(M, \eL_0, \eL_1) \otimes \FF_2[\theta, \theta^{-1}]\]
and whose $E^{\infty}$ page is isomorphic to 
\[\bigoplus_{j=1}^\ell \HF(M_j^{\fix}, \eL_{0,j}^{\fix}, \eL_{1,j}^{\fix}) \otimes \FF_2[\theta, \theta^{-1}].\]
\end{corollary}

The localization theorem of the introduction, Theorem~\ref{thm:localization-intro}, is the special case of Corollary~\ref{cor:components} in the case that only a subset of the components of $M^{\fix}$ contain a nonempty component of both fixed Lagrangians.

In order to make progress toward the proof of Theorem~\ref{t:genSS}, suppose that $M^{\fix}$ has a stable normal real polarization and let $N_{\mathbb{R}} \rightarrow M^{\fix}$ be the underlying real structure, so that $N_{\mathbb{R}}$ is a Lagrangian subbundle of $NM^{\fix}\oplus \CC^N$ for some $N$ with the property that $N_{\mathbb{R}} \otimes_{\RR} \CC \simeq NM^{\fix}\oplus \CC^N$. Let $h$ be a fibrewise metric for the real bundle $N_{\mathbb{R}} \to M^{\fix}$. We may extend $h$ to a fibrewise Hermitian metric on $N_{\mathbb{R}} \otimes_{\mathbb{R}} \mathbb{C}$ and give it the standard associated symplectic structure.
Explicitly, $N_{\mathbb{R}} \otimes_{\mathbb{R}} \mathbb{C}$ is a symplectic vector bundle whose fibrewise symplectic structure is given by 
\[
\omega(a_1+\sqrt{-1}b_1,a_2+\sqrt{-1}b_2):=h(b_1,a_2)-h(a_1,b_2).
\]
Let $\phi$ be the isomorphism of symplectic bundles
\[
\phi: NM^{\fix} \oplus \mathbb{C}^k \to N_{\mathbb{R}} \oplus N_{\mathbb{R}} ^* \simeq N_{\mathbb{R}} \otimes_{\mathbb{R}} \mathbb{C} 
\]
which sends $v \in N_{\mathbb{R}}$ to $v \oplus 0 \simeq v \otimes 1$. By assumption, for $i=0,1$, there exist smooth families of Lagrangian subbundles $F_{i,t}$ of the bundle $NM^{\fix} \oplus \mathbb{C}^k|_{\eL_i^{\fix}}$, with $t \in [0,1]$, such that
\begin{align}\label{eq:Fs}
F_{0,0}=N\eL_0^{\fix} \oplus \mathbb{R}^k, \quad F_{0,1}=N_{\mathbb{R}}, \quad F_{1,0}=N\eL_1^{\fix} \oplus \sqrt{-1} \mathbb{R}^k, \quad  F_{1,1}=\phi^{-1}(N_{\mathbb{R}} \otimes \sqrt{-1}).
\end{align}

There is a standard procedure to equip a neighborhood of the zero section of $N_{\mathbb{R}} \otimes_{\mathbb{R}} \mathbb{C}$ with a symplectic form, which we now recall (see for example \cite[Section 36]{GSbook}).
Let $r$ be the rank of $N_{\mathbb{R}}$ and $\pi_P \co P \to M^{\fix}$ be the principal $O(r)$-bundle underlying $N_{\mathbb{R}}$.
For $g \in O(r)$ and $p \in P$, our convention for the action is
\[
R_g \co P \to P \quad R_g(p):=pg^{-1}.
\]
Let  $\alpha \in \Omega^1(P,\mathfrak{o}(r))^{O(r)}$ be a connection $1$-form for $P \to M^{\fix}$.
Furthermore let $W:=T^*\mathbb{R}^r =\mathbb{R}^r \oplus (\mathbb{R}^r)^*$ be equipped with the standard $O(r)$ action 
\[
g \cdot (q,p)=(gq,g^{-T}p)=(gq,gp) \text{ for } (q,p) \in \mathbb{R}^r \oplus (\mathbb{R}^r)^*,
\]
where the last equality uses that $g \in O(r)$.
This action is Hamiltonian with respect to the canonical symplectic form $\omega_W$ on $W$, and the moment map 
$\mu \co W \to \fko(r)^*$
 is given by
\[
\langle \mu(q,p), A \rangle=p^TAq \text{ for all } A \in \fko(r).
\]
We consider the closed 2-form on $\widetilde{E}:=P \times W$ given by
\[
\omega_{\widetilde{E}}:=\pi_P^*\omega_{M^{\fix}} + \omega_{W} + d\langle \mu,\alpha \rangle.
\]
For $v \in \fko(r)$, the induced vector field on $\widetilde{E}$ is the sum of the induced vector fields $v_P \in TP$ on $P$ and $v_W \in TW$ on $W$.
One can check that
\begin{align*}
\iota_{v_P+v_W} \omega_{\tilde{E}}&=\iota_{v_W} \omega_{W}+\iota_{v_P+v_W}d\langle  \mu, \alpha \rangle \\
&=d\langle \mu, v \rangle-d\iota_{v_P+v_W}\langle \mu, \alpha \rangle\\
&=d\langle \mu, v \rangle-d\langle \mu, v \rangle \\&=0.
\end{align*}
Here the second equality uses the definition of the moment map and the Cartan formula, and the last equality uses the fact that
$\iota_{v_P+v_W} \alpha=\iota_{v_P} \alpha=v$ for any connection $1$-form.
Since $\omega_{\widetilde{E}}$ is $O(r)$-equivariant, it then follows that this form descends to a closed $2$-form $\omega_E$ on
\[
E:=P \times_{O(r)} W \simeq \Tot(N_{\mathbb{R}} \oplus N_{\mathbb{R}} ^*)
\]
where the latter notation indicates the total space of the bundle $N_{\mathbb{R}} \oplus N_{\mathbb{R}} ^*$.
It is well-known that $\omega_E$ is non-degenerate near the zero section.

\begin{lemma}\label{l:isotopy_realization}
The total space of $N_{\mathbb{R}}|_{\eL^{\fix}_i}$ is a Lagrangian submanifold of $(E,\omega_E)$ near the zero section.
\end{lemma}

\begin{proof}
Let $(q,p) \in \mathbb{R}^r \oplus (\mathbb{R}^r)^* =W$ and regard $q$ and $p$ as column vectors.
By definition, we have
\begin{align*}
\omega_{\widetilde{E}}&=\pi_P^*\omega_{M^{\fix}} + dp^T \wedge dq + d(p^T \alpha q) \\
&=\pi_P^*\omega_{M^{\fix}} + dp^T \wedge dq +dp^T \wedge \alpha q +p^T d\alpha q+p^T \alpha  \wedge dq.
\end{align*}
By construction, the subbundle $\{p=0\}=P \times \mathbb{R}^r \times \{0\}$ of $\widetilde{E}$ is preserved under the $O(r)$ action. 
When we restrict $\omega_{\tilde{E}}$ to $\{p=0\}$, it reduces to
\[
\pi_P^*\omega_{M^{\fix}}.
\]
Therefore the total space of $N_{\mathbb{R}}|_{\eL^{\fix}_i}$ is Lagrangian with respect to $\omega_{\widetilde{E}}$ because $\eL^{\fix}_i$ is itself Lagrangian with respect to $\omega_{M^{\fix}}$.
\end{proof}

\begin{remark}
In general, the total space of a Lagrangian subbundle over a Lagrangian submanifold is not a Lagrangian submanifold of $(E,\omega_E)$.
\end{remark}

We now consider the involution on $F$ given by $(q,p) \mapsto (-q,-p)$, which is symplectic and furthermore commutes with the $O(r)$ action. The induced involution on $\widetilde{E}$ acts trivially on the $P$ factor, and
descends to a symplectic involution on $(E,\omega_E)$ whose fixed locus is exactly the zero section.

Fix a relatively compact open subset $U \subset M^{\fix}$ containing $\eL_i^{\fix}$.
By a Moser argument, we can find an invariant open subset $V \subset M \times \mathbb{C}^k$ with $U \subseteq V \cap (M^{\fix}\times\{0\})$ and an equivariant symplectic embedding
\[
\Phi \co V \to E
\]
such that $\Phi|_{V \cap M^{\fix}}$ is the inclusion to the zero section, and $d\Phi|_{NM^{\fix} \oplus \mathbb{C}^k}=\phi$.

The following lemma should be compared with~\cite[Equation 64]{SS10}.

\begin{lemma}
Under the assumptions of Theorem \ref{t:genSS}, for $i=0,1$, there exists an equivariant exact Lagrangian isotopy from $\eL_i \times \mathbb{R}^k$ to $\widetilde \eL_i$ such that $\eL_i^{\fix} \times \{0\}$ is fixed by the isotopy and $\Phi(\widetilde\eL_0)=\Tot(N_{\mathbb{R}}|_{\eL_0^{\fix}})$ and $\Phi(\widetilde\eL_1)=\Tot(N_{\mathbb{R}}^*|_{\eL_1^{\fix}})$.
\end{lemma}

\begin{proof}
Let $\Phi$ and $V$ be as above, and recall the bundles $F_{i,t}$ appearing in Equation~\ref{eq:Fs} for $i=0,1$. By slight abuse of notation, during this proof we will denote the image $\Phi((\eL_i \times \mathbb{R}^k)\cap V)$ by $\eL_i$ and $\phi(F_{i,t})$ by $F_{i,t}$. 
In other words, we begin the argument with two Lagrangian submanifolds $\eL_i$ in $E$ with fixed loci $\eL_i^{\fix} \subset M^{\fix}$, together with an isotopy of Lagrangian subbundles $F_{i,t}$ from $F_{i,0}=N\eL_i^{\fix}$ to $F_{0,1}=N_{\mathbb{R}}$ and $F_{1,1}=N_{\mathbb{R}}^*$ as appropriate in $N_{\mathbb{R}} \oplus N_{\mathbb{R}}^*$.

Notice that for any invariant Lagrangian submanifold $\eL$ whose fixed locus is given by $\eL^{\fix}$ and any point $x \in \eL$, the tangent space $T_{x}\eL^{\fix}$ is the eigenspace of $T_x\eL$ with respect to the induced involution on the tangent space with eigenvalue $1$, and $N_x\eL^{\fix}$ is canonically isomorphic to the $-1$-eigenspace of $T\eL$. It follows that $N\eL_i^{\fix}$ is naturally a Lagrangian subbundle of $N_{\mathbb{R}} \oplus N_{\mathbb{R}}^*$.
Our goal is to show that we have sufficient invariant Hamiltonian functions that we may realize any isotopy of Lagrangian subbundles of $N_{\mathbb{R}} \oplus N_{\mathbb{R}}^*$ as an equivariant isotopy of Lagrangian submanifolds of $\eL_i$ relative to $\eL_i^{\fix}$.

Let $U_0 \subseteq M^{\fix}$ be an open subset and $P|_{U_0} \simeq U_0 \times O(r)$ be a trivialization.
Note that $U_0 \times \{e \} \times W$ is a slice of the $O(r)$ action, so we have a natural diffeomorphism 
\[U_0 \times \{e \} \times W \simeq U_0 \times O(r) \times_{O(r)} W=E|_{U_0},\] and we identify these sets throughout.
Moreover under this identification, we have \[\omega_{\tilde{E}}|_{U_0 \times \{e \} \times W}=\omega_E|_{U_0 \times O(r) \times_{O(r)} W}.\]

Let $\mathcal Q$ be the space of quadratic forms $Q$ on the vector space $W$.
For any $g \in O(r)$, the pullback of $Q$ by $g$ is still a quadratic form, so $O(r)$ acts on $\mathcal Q$.
Moreover, the pullback of $Q$ by the involution on $F$, which is to say by $g=-\mathrm{Id}$, preserves $Q$.

Let $\mathcal H$ be the space of smooth Hamiltonian functions $H \co E \to \mathbb{R}$ such that $H$ is fibrewise a quadratic form; that is to say, such that when we restrict $H$ to a trivialization
$ E|_{U_0} \simeq U_0 \times \{e \} \times W$, it takes the form 
\[
H(u, -) \in \mathcal Q \text{ for all }u \in U_0.
\]
Note that, since $O(r)$ acts on $\mathcal Q$, $H$ being fibrewise a quadratic form is a well-defined notion that is independent of the choice of trivialization.

For $H \in \mathcal H$, the function $H$ is invariant under the involution, and hence its time-$s$ Hamiltonian flow $\psi^s$ is equivariant.
Let $u_0 \in \eL^{\fix}_i$ and $v_0 \in N_{(u_0,e,0)}\eL^{\fix}_i$.
Since $v_0$ lies in the $-1$-eigenspace, we know that  $v_0=(0,v) \in T_{(u,e,0)} U_0 \times \{e\} \times W$  for some $v \in W$.
We are interested in the infinitesimal deformation of $\eL_i$ by $\psi^s$ near $\eL_i^{\fix}$. We therefore compute $\frac{d}{ds}|_{s=0}\psi^s(v)$. 
Let $\gamma(t)$ be a curve such that $\gamma(0)=(u,e,0)$ and $\gamma'(0)=v$. Then we have
\[
\frac{d}{ds}|_{s=0}\psi^s(v)=\frac{d}{ds}|_{s=0} \frac{d}{dt}|_{t=0} \psi^s(\gamma(t))= \frac{d}{dt}|_{t=0}\frac{d}{ds}|_{s=0} \psi^s(\gamma(t))=\frac{d}{dt}|_{t=0}X_H(\gamma(t))
\]
where $X_H$ is the Hamiltonian vector field of $H$.

To understand the vector field $X_H$, we compute
\[
dH=D_u H du+ D_pH dp+D_qH dq
\]
where $D_u, D_p, D_q$ are the total partial derivatives in the $u, p, q$ coordinates respectively.
Note that $D_pH$ and $D_qH$ are valued in linear forms as functions of $u$ because $H(u, -) \in \mathcal Q$.
We are only interested in $X_H$ in a small neighborhood of the zero section, so we may assume $p$ and $q$ are close to $0$.
Therefore the term 
\[
\epsilon:=d(p^T \alpha q)=dp^T \wedge \alpha q +p^T d\alpha q+p^T \alpha  \wedge dq,
\]
which is zero when $p=q=0$, is small in comparison to the forms
\[
\omega_0:=\pi_P^*\omega_{M^{\fix}} + dp^T \wedge dq
\]
and
\[
\omega_{\widetilde{E}}|_{U_0 \times \{e\} \times W}=\omega_0+\epsilon.
\]
Let $X_0$ be the $\omega_0$-dual of $dH$. Explicitly, 
if we define $X_{H|_{U_0}}$ to be the $\omega_{M^{\fix}}$-dual of $H|_{U_0 \times \{0\} \times \{0\}}$
and $X_{H(u,-)}$ to be the $dp^T \wedge dq$-dual of the quadratic form $H(u,-)$, then we have
\[
X_0=X_{H|_{U_0}}+X_{H(u,-)}.
\]
As a result, for the curve $\gamma(t)$ defined as above,
\[
\frac{d}{dt}|_{t=0}X_0(\gamma(t))= \frac{d}{dt}|_{t=0}X_{H(u,-)}(\gamma(t))=\Hess(H(u,-)) \cdot v
\]
where $\Hess(H(u,-))$ is the Hessian of the quadratic form $H(u,-)$ with respect to the $(p,q)$ coordinates.

We claim that for $X_H$, we also have
\begin{align}\label{eq:quad_var}
\frac{d}{dt}|_{t=0}X_H(\gamma(t))=\frac{d}{dt}|_{t=0}X_0(\gamma(t))=\Hess(H(u,-)) \cdot v.
\end{align}
To see this, we use that for square matrices $B=B_0+C$ such that $B$ is invertible and $C$ is small compared to $B$, we have
\[
B^{-1}=(B_0+C)^{-1}=B_0^{-1}-B_0^{-1}C B_0^{-1}+B_0^{-1}CB_0^{-1}CB_0^{-1}-\dots
\]
If we regard of $\omega, \omega_0,$ and $\epsilon$ as linear maps from tangent bundle to cotangent bundle giving the matrices $B$, $B_0$, $C$, respectively, we see that
\[
X_H=B^{-1} \cdot dH=X_0-B_0^{-1}CX_0+B_0^{-1}CB_0^{-1}CX_0- \dots
\]
from which it follows that
\begin{align*}
\frac{d}{dt}|_{t=0}(X_H-X_0)(\gamma(t))&=-\frac{d}{dt}|_{t=0} \left(B_0^{-1}(\gamma(t))C(\gamma(t))X_0(\gamma(t))\right)=0
\end{align*}
where the last equality follows from the product rule and the fact that $C(\gamma(0))=X_0(\gamma(0))=0$.

Let $\mathcal P$ be the space of Lagrangian subbundles of $E|_{\eL^{\fix}_i}$.
In particular, $N\eL^{\fix}_i$ defines a point in  $\mathcal P$.
By Equation \eqref{eq:quad_var} and our freedom on the choice of quadratic forms $H(u,-)$, we have enough Hamiltonians in $\mathcal H$ to achieve any deformation direction of an element in $\mathcal P$, in other words, for any Lagrangian subbundle $P_0 \in \mathcal P$, the map $T_0\mathcal H \to T_{P_0} \mathcal P$ is surjective.

Recall that $F_{i,t}$ is our assumed choice of Lagrangian isotopy from $\eL_i \oplus \mathbb{R}^k$ to $\widetilde \eL_i$.
Without loss of generality, we can assume $F_{i,t}$ is smooth and hence defines a smooth path in $\mathcal P$. 
We can use a time-dependent Hamiltonian function $H=(H_t)$ such that $H_t \in \mathcal H$ for all $t$ to get an exact Lagrangian isotopy whose induced path on the Lagrangian subbundle is exactly $F_{i,t}$.
\end{proof}

The connection $1$-form $\alpha$ induces a connection $1$-form for the bundle $\pi_E \co E \to M^{\fix}$.
For any point $x \in E$, we have the induced splitting
\begin{align}\label{eq:splitT}
T_xE=T_x^{hor}E \oplus T_x^{vert}E
\end{align}
where $T_x^{vert}E$ is the kernel of $d\pi_E|_x$ and $T_x^{hor}E$ is the kernel of the connection $1$-form.
We have the canonical identifications $T_x^{vert}E \simeq W$ and $d\pi_E:T_x^{hor}E \simeq T_{\pi_E(x)} M^{\fix}$.
Let $J_{std}$  be the standard complex structure on $F$ given by $\partial_p \mapsto \partial_q$ and $\partial_q \mapsto -\partial_p$.
We say that an almost complex structure $J_E$ of $E$ has \emph{diagonal type} if there is an almost complex structure $J_{M^{\fix}}$ of $M^{\fix}$ such that $J_E = J_{M^{\fix}} \oplus J_{std}$ with respect to \eqref{eq:splitT} and the canonical identifications.

\begin{lemma}\label{l:regular_u}
Let $\eL_0^{\fix}, \eL_1^{\fix}$ be two transversely intersecting compact exact Lagrangian submanifolds of $M^{\fix}$.
Let $\eL_0:=\Tot(N_{\mathbb{R}}|_{\eL_0^{\fix}})$ and $\eL_1:=\Tot(N_{\mathbb{R}}^*|_{\eL_1^{\fix}})$.
For a diagonal type almost complex structure on $E$, a finite energy Floer solution $u \co \mathbb{R} \times [0,1] \to M^{\fix}$ with boundary on $(\eL_0^{\fix}, \eL_1^{\fix})$ is regular if and only if it is regular as a Floer solution to $E$ with boundary on $(\eL_0,\eL_1)$.
\end{lemma}

\begin{proof}
The linearized operator $D_u$ for $u$ as a map to $E$ splits into the direct sum of the $\iota$-invariant part and  the $\iota$-anti-invariant part.
The invariant part is the same as $D_u$ for  $u$ as a map to $M^{\fix}$, which is assumed to be surjective.
Since $N_{\mathbb{R}}$ and $N_{\mathbb{R}}^*$ are subbundles of $E$ which are preserved under parallel transport with respect to the connection, 
the anti-invariant part $D_u^{\perp}: W^{1,2}(\mathbb{R} \times [0,1], u^*E, u^*N_{\mathbb{R}}, u^*N_{\mathbb{R}}^*) \to L^2(\mathbb{R} \times [0,1], u^*E,u^*N_{\mathbb{R}}, u^*N_{\mathbb{R}}^*)$ can be identified with the linearized operator for the constant map to $0 \in T^*\mathbb{R}^r$ with constant Lagrangian boundary condition $\mathbb{R}^r$ and $T^*_0 \mathbb{R}^r$ along the two boundary components.
The latter is well-known to be surjective.
\end{proof}

\begin{proof}[Proof of Theorem \ref{t:genSS}] We closely follow the proof of Theorem \ref{thm:ss-localization} appearing in \cite[Section 3]{SS10}.

Given $\iota$-equivariant Lagrangian submanifolds $\eL_0, \eL_1 \subset M$, the products $\eL_0':=\eL_0 \times \mathbb{R}^k$ and $\eL_1':=\eL_1 \times \sqrt{-1} \mathbb{R}^k$ are equivariant Lagrangians in $M \times \mathbb{C}^k$.
By Lemma \ref{l:isotopy_realization}, $\eL_0'$ is equivariantly Hamiltonian isotopic via a compactly supported Hamiltonian to a new Lagrangian 
$\widetilde{\eL}_0$ such that near $M^{\fix}$, $\widetilde{\eL}_0$ can be identified with $\Tot(N_{\mathbb{R}}|_{\eL_0^{\fix}})$.
Similarly, $\eL_1'$ is equivariantly Hamiltonian isotopic to $\widetilde{\eL}_1$ such that near $M^{\fix}$, $\widetilde{\eL}_0$ can be identified with $\Tot(N_{\mathbb{R}}^*|_{\eL_1^{\fix}})$.

As in \cite[Equation 65]{SS10}, we have quasi-isomorphisms
\begin{align*}
CF(M, \eL_0,\eL_1) &\simeq CF(M\times \mathbb{C}^k,\widetilde{\eL}_0,\widetilde{\eL}_1) \\
CF_{borel}(M, \eL_0,\eL_1) &\simeq CF_{borel}(M \times \mathbb{C}^k, \widetilde{\eL}_0,\widetilde{\eL}_1) 
\end{align*}
If we use a family of compatible almost complex structures on $M \times \mathbb{C}^k$ that are of diagonal type near $M^{\fix}$, we can then apply Lemma \ref{l:regular_u} to achieve regularity of Floer solutions that lie completely in $M^{\fix}$.
This allows us to run the construction of \cite[Sections 2b, 2c and 3e]{SS10} to get a chain complex $CF_{equiv}(\widetilde{\eL}_0,\widetilde{\eL}_1) $ together with a canonical quasi-isomorphism
\[
 CF_{borel}(M\times \mathbb{C}^k,\widetilde{\eL}_0,\widetilde{\eL}_1)  \to CF_{equiv}(M \times \mathbb{C}^k,\widetilde{\eL}_0,\widetilde{\eL}_1)
\]
and a localization map
\[
 CF_{equiv}(M\times \mathbb{C}^k, \widetilde{\eL}_0,\widetilde{\eL}_1) \to CF(M^{\fix}, \eL_0^{\fix},\eL_1^{\fix})[[\theta]].
\]
The result now follows from a homological argument identical to that appearing in \cite[Section 2c and 3f]{SS10}. Note that the original uses the variable $q$ rather than $\theta$.
\end{proof}

\section{Manolescu's construction of symplectic Khovanov--Rozansky homology} \label{sec:construction}

In this section, we give a quick introduction to the symplectic Khovanov--Rozansky homology defined by Manolescu \cite{Manolescu:sln}, generalizing Seidel--Smith's construction of symplectic Khovanov homology \cite{SS06}. Before beginning, we note that in symplectic link homologies, there is generally a choice of whether to describe the ambient symplectic manifold as a subset of the Lie algebra $\frg =\frsl(mn)$ or $\frg = \frgl(mn)$; Manolescu \cite{Manolescu:sln} uses $\frsl(mn)$, but here we use $\mathfrak{gl}(mn)$ to simplify the discussion in Sections~\ref{sec:intrinsic} and \ref{sec:extrinsic}. The two constructions are carried out nearly identically.  A similar choice is made in \cite[Section 5]{SS10}, where $\frgl(2m)$ is used instead of  $\frsl(2m)$ to simplify the study of an action on the space.

\subsection{Lie theory preliminaries} \label{sec:lie} Let $\frg=\frgl(mn)$. Let $\pi = (\pi_1,\dots, \pi_s)$ be a partition of $mn$, arranged so that $\pi_1 \ge \dots \ge \pi_s>0$. We sometimes also write $\pi$ in the form $((a_1)^{k_1}\dots(a_{\ell})^{k_{\ell}})$, where $a_1 > \dots > a_{\ell}$, and $k_i$ is the number of indices $j$ with $\pi_j=a_i$. We consider the partial flag variety associated to $\pi$, that is
\[
\mathcal{F}^{\pi}=\{0=F_0 \subset F_1  \subset \dots \subset F_s =\mathbb{C}^{mn}| \dim(F_j)-\dim(F_{j-1})=\pi_j\}.
\]
Given a flag $F \in \mathcal{F}^{\pi}$, let $\frp(F) \subset \frg$ be the parabolic subalgebra corresponding to $F$, consisting of the matrices $x$ preserving the flag, that is, such that $x(F_j) \subset F_j$ for all $j$.
Let
\begin{align*}
\hat{\frg}^{\pi}&=\{ (x,F) | F \in \mathcal{F}^{\pi}, x \in \frp(F)\}.
\end{align*}
Then for any $(x,F) \in \hat{\frg}^{\pi}$, there is an induced action $x^F_j$ on $F_j/F_{j-1}$. We consider the elements whose action is central, and set
\[
\tilde{\frg}^{\pi}=\{ (x,F) \in \hat{\frg}^{\pi}| x^F_j \in Z(\frgl(F_j/F_{j-1}))\}.
\]
The image of the natural projection $\tilde{\frg}^{\pi} \to \frg$ is denoted by $\frg^{\pi}$. 
Let $\frh \simeq \mathbb{C}^{mn}$ be the space of diagonal matrices, $W_{\pi}:=S_{\pi_1} \times \dots \times S_{\pi_s}$ be the product of permutation groups, and $\frh^{\pi} \subset \frh/W_{\pi}$ be the equivalence class of diagonal matrices of the type $\alpha=diag(\alpha_1,\dots, \alpha_1,\alpha_2, \dots, \alpha_2,\dots, \alpha_s,\dots,\alpha_s)$ where $\alpha_j$ appears exactly $\pi_j$ times. Let $m_{\pi}(k)$ be the number of times that $k$ appears among $\pi_j$ and $W^{\pi}=S_{m_{\pi}(1)} \times \dots \times S_{m_{\pi}(mn)}$. Then, we get a well-defined map $\chi^{\pi}:\frg^{\pi}\to \mathfrak{h}^{\pi}/W^{\pi}$, sending each matrix in $\frg^{\pi}$ to its eigenvalues (with multiplicity).
The \emph{partial resolution} of $\chi^{\pi}$ is the commutative diagram
\[
        \begin{tikzcd}
    \tilde{\frg}^{\pi} \arrow[r] \arrow[d, "\tilde{\chi}^{\pi}"] & \frg^{\pi} \arrow[d, "\chi^{\pi}"] \\
     \frh^{\pi} \arrow[r] & \frh^{\pi}/W^{\pi}
        \end{tikzcd}
\]
A general element of $\frh^{\pi}/W^{\pi}$ will be denoted $\tau$. The map $\tilde{\chi}^{\pi}$ is a $C^{\infty}$ fibre bundle which admits an explicit (albeit non-canonical) trivialization.
On the other hand, the map $\chi^{\pi}$ is not a fibre bundle and $\frg^{\pi}$ is a singular variety.
The subspace of regular values of $\chi^{\pi}$ is $\frh^{\pi, reg}/W^{\pi}$, where $\frh^{\pi, reg}\subset \frh^{\pi}$ consists of matrices such that the values of $\alpha_j$ are pairwise distinct. We denote the preimage of this subspace in $\frg^{\pi}$ by 
$\frg^{\pi, reg}$, so that $\chi^{\pi}\co\frg^{\pi,reg}\to \frh^{\pi, reg}/W^{\pi}$ is a differentiable fibre bundle \cite[Proposition 2.4]{Manolescu:sln}.

We are primarily interested in the standard partition $\pi=((n-1)^m 1^m)$. In this case, $\frh^{\pi,reg}/W^{\pi}$ is identified with the following bipartite configuration space 
\[\mathrm{BConf}_m=\{(\{\mu_1,\mu_2,\cdots,\mu_m\}, \ \{\lambda_1,\lambda_2,\cdots,\lambda_{m}\}) \ |\ \mu_i,\lambda_j\in\mathbb{C}\ \text{distinct} \}.\]

Here each $\tau=(\vec{\mu},\vec{\lambda})\in\mathrm{BConf}_m$ is identified with the equivalence class of the diagonal matrices with eigenvalues $\vec{\mu}=\{\mu_1,\mu_2,\cdots,\mu_{m}\}$ of multiplicity $n-1$ and $\vec{\lambda}=\{\lambda_1,\lambda_2,\cdots,\lambda_{m}\}$ of multiplicity $1$. Each $\mu_i$ is called a thick eigenvalue, whereas each $\lambda_i$ is called a thin eigenvalue.

\subsection{The ambient symplectic manifold} \label{sec:ambient}
Let $\cS_{m,n}\subset \frgl(mn)$ be the affine subspace
\[
\cS_{m,n}=\left\{  \begin{pmatrix} Y_{11} & I_n & &\\
Y_{21} & 0 & I_n &&\\
\dots & & &\dots   &\\
Y_{(m-1)1} & & & &I_n\\
Y_{m1}& 0 & 0 & 0&0
\end{pmatrix} : Y_{j1} \in \frgl(n) \text{ for }j \ge 1   \right\}. 
\] 
This is a transverse slice to the adjoint orbit of the nilpotent element with $m$ Jordan blocks of size $n$ \cite[Lemma 3.1]{Manolescu:sln}. 

By \cite[Proposition 2.17]{Manolescu:sln},
\[
\chi^{\pi}|_{\cS_{m,n}} : \frg^{\pi,reg} \cap \cS_{m,n} \to \frh^{\pi,reg}/W^{\pi}
\]
is a differentiable fibre bundle.
For $\tau \in \frh^{\pi, reg}/W^{\pi}$, we consider the following smooth subspace in $\frg^{\pi}$
\[
\cY_{m,n,\tau}^\pi:=\cS_{m,n} \cap (\chi^{\pi})^{-1}(\tau).  
\]
For readability, if we are using the standard partition or we expect that the partition is otherwise clear, we write $\cY_{m,n,\tau}^\pi$ as $\cY_{m,n,\tau}$. If we wish to emphasize that we are using the standard partition for the values of $m$ and $n$ specified (for example when multiple possible values of $m$ and $n$ are under consideration) we sometimes write $\pi_{std}$.
We will also want to understand certain cases in which $\tau$ is not a regular value, and the space $\cY_{m,n,\tau}$ is singular. In particular, if the value of a thick eigenvalue $\mu_j$ and a thin eigenvalue $\lambda_i$ coincide, then $\tau\in \frh^{\pi}/W^{\pi}$ is not a regular value of $\chi^{\pi}$. For instance, suppose $\mu_1=\lambda_1=\lambda$. Then the singular set of $\cY_{m,n,\tau}$, denoted by $\cC_{m,n,\tau}$, consists of matrices such that the eigenspace of $\lambda$ is $n$-dimensional. The union of all such singular sets $\cC_{m,n,\tau}$ as $\tau$ varies is equal to $\cC_{m,n}=\cS_{m,n}\cap \frg^{\pi',reg}$ for the partition  $\pi'=(n(n-1)^{m-1}1^{m-1})$.

The following lemma of Manolescu will feature significantly in our arguments, so we recall its statement here.

\begin{lemma} \cite[Lemma 3.2]{Manolescu:sln} \label{lem:injectivity} Let $Y \in \cS_{m,n}$, and for any $\alpha$ an eigenvalue of $Y$ consider the eigenspace $E_{Y}(\alpha) = \mathrm{Ker}(\alpha I - Y) \subset \CC^{mn}$. The map $E_{Y}(\alpha) \rightarrow \CC^n$ given by projection onto the first $n$ coordinates of $\CC^{mn}$ is an injection. \end{lemma}

\begin{proof}
The original \cite[Lemma 3.2]{Manolescu:sln} is stated with respect to $\frsl(mn)$, but the same proof applies for  $\frgl(mn)$.
\end{proof}

We particularly observe the following immediate corollary.

\begin{corollary}\label{cor:too-big-is-trivial}
If any entry in the partition $\pi$ is greater than $n$, then $\cY_{m,n,\tau}^\pi$ is necessarily empty for any $\tau$.
\end{corollary}

We now choose a symplectic form on all of $\cS_{m,n} \cap \frg^{\pi,reg}$. First let $\beta>m$ be a real number and $\xi_k(z)=|z|^{\frac{2\beta}{k}}$. 
We can construct a function $\eta_k:\mathbb{C} \to \mathbb{R}_{\ge 0}$ such that 
\[
\psi_k:=\xi_k+\eta_k
\]
is real analytic near the origin, smooth and strictly plurisubharmonic everywhere, as in \cite[Section 7.1]{Manolescu:sln}. We additionally may require that $\eta_k(z)=\eta_k(-z)$, so that in total $\psi_k(z)=\psi_k(-z)$.
Define a function $\hat{\psi}:\cS_{m,n} \to \mathbb{R}$ by summing $\psi_k(a)$ for every entry $a \in Y_{k1}$. Restricting this function to $\frg^{\pi,reg}$, we get a function 
\begin{align}\label{eq:potential}
\psi:\cS_{m,n} \cap \frg^{\pi,reg} \to \mathbb{R}
\end{align}
We then equip $\cS_{m,n} \cap \frg^{\pi,reg}$  with the symplectic form $\Omega=-dd^c \psi$.

Let $\bar{\tau}=(\{\mu_2,\dots,\mu_m\}, \{\lambda_2,\dots,\lambda_m\} )$ and $\tau_0=(\{ 0,\mu_2,\dots,\mu_m\}, \{0,\lambda_2,\dots,\lambda_m\})$.
By \cite[Lemma 3.3]{Manolescu:sln} $\cY_{m-1,n,\bar{\tau}}$ is identified with the singular set of $\cY_{m,n,\tau_0}$. Recall that this set is denoted by $\mathcal{C}_{m,n,\tau_0}$ and is the subspace consisting of matrices $Y\in\cY_{m,n,\tau_0}$ such that $\ker(Y)$ is $n$-dimensional.  More precisely, $Y\in \mathcal{C}_{m,n,\tau_0}$ if and only if $Y_{m1}=0$. So, one can define a bijection $f$ from $\cY_{m-1,n,\bar{\tau}}$ to $\cC_{m,n,\tau_0}$ by setting
\begin{align}\label{eq:embed-sln}
f\left(\begin{matrix}Y_{11}&I_n&&&\\Y_{21}&0&I_n&&\\\cdots&&&\cdots&\\Y_{(m-1)1}&0&0&0&0\end{matrix}\right)=\left(\begin{matrix}Y_{11}&I_n&&&\\Y_{21}&0&I_n&&\\\cdots&&&\cdots&\\Y_{(m-1)1}&&&&I_n\\0&0&0&0&0\end{matrix}\right).
\end{align}
Considering the definition of $\psi$ from \eqref{eq:potential}, it is straightforward that $f$ is a symplectomorphism.

Before moving on, we recall the special case $m=1$ from \cite[Section 3.3]{Manolescu:sln}. In this case the slice $\cS_{1,n}$ is all of $\frgl(n)$.  
Consider the subvariety $\cZ_{\frgl}=\frgl(n)^{((n-1)1)}$ consisting of matrices with an eigenspace of dimension at least $n-1$. Then the spaces $\cY_{1,n,\tau}$ are the fibres of the map
\begin{equation} 
 \chi^{((n-1)1)}_{\frgl} \co  \frgl(n)^{((n-1)1)} \rightarrow \CC \times \CC \end{equation}
which takes a matrix into its ordered eigenvalues. More precisely, the first factor is the eigenvalue of multiplicity $\ge n-1$  and the second factor is the remaining eigenvalue. If there is in fact one eigenvalue of multiplicity $n$ then the values coincide.
We will also have cause to consider the subvariety $\cZ=\frsl(n)^{((n-1)1)} \subset \cZ_{\frgl}$ consisting of traceless matrices. In total these spaces fit into a commutative diagram as follows.
\begin{equation}\label{eqn:m=1}
\begin{tikzcd}
		\CC \times \cZ\arrow[r,"\simeq"] \arrow[d, "id_\CC \times \chi^{((n-1)1)}"]&  \cZ_{\frgl} \arrow[d,"\chi^{((n-1)1)}_\frgl"]\\
		\CC \times \mathbb{C}  \arrow[r, "\simeq"] & \mathbb{C}  \times\mathbb{C} 
	\end{tikzcd}
\end{equation}
The maps are as follows. The map $\chi^{((n-1)1)}$ takes a matrix into the eigenvalue of multiplicity at least $n-1$. The top horizontal arrow sends $(a,A)$ to $A+aI_n$,
and the bottom arrow sends $(a,\mu)$ to $(\mu+a, -(n-1)\mu+a)$.
This diagram will play an important role in constructing local models for neighborhoods of $\cS_{m,n}$ for general $m$ in Section~\ref{sec:lagrangian}.

\subsection{Rescaled parallel transport} \label{sec:parallel} From the previous section, we have a fibration \[\cS_{m,n} \cap \frg^{\pi,reg} \to \frh^{\pi,reg}/W^{\pi}\] which is symplectic. We now proceed to discuss a notion of rescaled parallel transport for this fibration following \cite[Sections 4.1 and 7.1]{Manolescu:sln}. Toward this end, given $\tau \in \frh^{\pi,reg}/W^{\pi}$, we let $Z_{\tau}$ be the induced Liouville vector field on $\cY_{m,n,\tau}$; that is, let $Z_{\tau}$ be the vector field such that \[\iota_{Z_{\tau}} \Omega|_{\cY_{m,n,\tau}}=-d^c \psi|_{\cY_{m,n,\tau}}\] where $\psi$ is the potential function defined in \eqref{eq:potential}.

For a path $\gamma\co[0,1] \to \frh^{\pi,reg}/W^{\pi}$, there is some constant $c>0$ such that the critical points of $\psi|_{\cY_{m,n,\gamma(t)}}$ are contained in $\psi^{-1}([0,c))$.
Let $H_{\gamma}$ be the symplectic parallel transport vector field associated to $\gamma$. In some cases, integrating $H_{\gamma}$ directly supplies a symplectomorphism $h_{\gamma} \co Y_{\gamma(0)} \rightarrow Y_{\gamma(1)}$; this map is called naive parallel transport. However, in general the lack of compactness means that integral flow lines may not exist everywhere, and in particular naive parallel transport is only necessarily defined on compact subsets $P$ for a short time $\epsilon$ which depends on $P$. Instead one must rely on rescaled parallel transport, as follows.

For $s >0$ large enough, the integral lines of
\[
\widetilde{H}_{\gamma,s}:=H_{\gamma}- s Z_{\gamma(t)}
\]
remain inside $\psi^{-1}([0,c))$. The rescaled parallel transport map
\[
h^{resc}_{\gamma}: \cY_{m,n,\gamma(0)} \cap \psi^{-1}([0,c)) \to \cY_{m,n,\gamma(1)} 
\]
 is defined to be the composition of the flow of $\widetilde{H}_{\gamma,s}$ with a time $s$ Liouville flow on $\cY_{m,n,\gamma(1)}$.  The map $h^{resc}_{\gamma}$ is a symplectic embedding which is independent of the parameter $s$ up to exact isotopy, as in \cite[Section 4(A)]{SS06} and \cite[Section 4.1]{Manolescu:sln}.

We will also want to consider rescaled parallel transport on the singular set. Toward this end, consider the partition $\pi'=(n(n-1)^{m-1}1^{m-1})$.
As explained in the last paragraph of \cite[Section 7.1]{Manolescu:sln}, the rescaled parallel transport map can be similarly defined for the fibration from the corresponding singular sets 
\[
\mathcal C_{m,n} \to \frh^{\pi',reg}/W^{\pi'} \subset \frh^{\pi}/W^{\pi}.
\]
Moreover, the rescaled parallel transport agrees with the parallel transport map on $\cY_{m-1,n, \bar{\tau}}$ under the identification given in \eqref{eq:embed-sln}.

\subsection{Crossingless matching Lagrangians} \label{sec:lagrangian}

In this section we review the construction of the Lagrangians used to define symplectic $\frsl(n)$ homology via vanishing projective cycles.

\subsubsection{The vanishing projective cycle construction and the case of $m=1$} \label{sec:vanishing-projective} We begin by reviewing the basic construction of vanishing projective cycles, following \cite[Section 4.3]{Manolescu:sln}. Recall the space $\cZ = \frsl(n)^{((n-1)1)}$ introduced at the close of Section~\ref{sec:ambient}, which consists of matrices $A \in \frsl(n)$ which have an eigenvalue $\alpha$ with an eigenspace of dimension at least $n-1$. As previously, we have the map
\[ \chi^{((n-1)1)} \co \frsl(n)^{((n-1)1)} \rightarrow \mathbb C\]
which goes by projection to the eigenvalue $\alpha$. The following alternate description of this space is a minor rephrasing of \cite[Lemmas 4.4 and 4.8]{Manolescu:sln}.

\begin{lemma} \label{lem:cotangent} The space $\cZ$ is a K{\"a}hler quotient of $T^*\CC^n$ by the action of the circle on $T^*\CC^n$ having weights $1$ and $-1$, each with multiplicity $n$. Furthermore, \keh{if $\alpha \neq 0$} then there is a symplectomorphism 
\begin{align}\label{eq:phi_symp}
 \phi \co \cZ_\alpha = \cY_{1,n, \tau} \rightarrow (T^* \CC \PP^{n-1}, \omega_{can})
\end{align}
where $\omega_{can}$ is the canonical symplectic form of a cotangent bundle and $\tau = (\alpha, (1-n)\alpha)$ is a configuration with multiplicities corresponding to the standard partitition $\pi=((n-1)1)$ of $n$.

\end{lemma}

\begin{proof} Let $(z,w) \in \mathbb{C}^n \times \mathbb{C}^n \simeq T^*\mathbb{C}^n$.
	We equip $T^*\mathbb{C}^n$ with the holomorphic symplectic form $\Omega:=\sum dz_i \wedge dw_i$ and the real symplectic form 
	$\omega:=\frac{i}{2}\sum_i (dz_i \wedge d\bar{z}_i+dw_i \wedge d\bar{w}_i)$.
	Consider the $S^1$ action on $(z,w)$ given by $t \cdot (z,w):=(e^{\sqrt{-1}t}z, e^{-\sqrt{-1}t}w)$; this is the action described in the statement of the lemma.
	The real and complex moment maps are given by
	\[
	\mu_{\mathbb{R}}(z,w)=\sum_{i=1}^n (|z_i|^2-|w_i|^2) \text{ and } \mu_{\mathbb{C}}(z,w)=\sum_{i=1}^n z_iw_i.
	\]
	
	As observed in \cite[Lemmas 4.4 and 4.8]{Manolescu:sln}, the map
	\[
	(z,w) \mapsto z^Tw-\frac{1}{n}(z \cdot w)I_n
	\]
	
	from $T^*\CC^n$ to $\frsl(n)^{((n-1)1)}$ is invariant under the action and descends to an isomorphism
	\[
	f:T^*\mathbb{C}^n//S^1:=\mu^{-1}_{\mathbb{R}}(0)/S^1 \to \frsl(n)^{((n-1)1)}
	\]
	Furthermore, as observed in \cite[Lemma 4.8]{Manolescu:sln}, the composition of $f$ and $\chi^{((n-1)1)}$ is 
	\[
	(z,w) \mapsto -\frac{1}{n}(z \cdot w)
	\]
	which is exactly $-\frac{1}{n}\mu_{\mathbb{C}}$.
	Therefore, $\cY_{1,n, \tau}$ is precisely the hyperK{\"a}hler quotient of $T^*\mathbb{C}^n$ by $S^1$, which is well-known to be $T^*\mathbb{CP}^{n-1}$.

We also give an explicit proof of the last statement without using hyperK\"ahler rotation.
Let $\psi(z,w)=|z|^2+|w|^2$ be the K\"ahler potential of the standard K\"ahler form on $T^*\mathbb{C}^n$.
The map $\psi$ descends to the quotient $T^*\mathbb{C}^n//S^1=\mu^{-1}_{\mathbb{R}}(0)/S^1$.

One may check that the $U(n)$ action $U \cdot (z,w)=(zU^{T},wU^{-1})$ preserves the locus $\mu^{-1}_{\mathbb{R}}(0) \cap \mu_{\mathbb{C}}^{-1}(1)$ and $\psi$ is $U(n)$-invariant.
For  $(z,w) \in \mu^{-1}_{\mathbb{R}}(0) \cap \mu_{\mathbb{C}}^{-1}(1)$, let $U \in U(n)$ be such that only the first coordinate of $zU^{T}$ is non-zero.
Let $(z',w')=U \cdot (z,w)$.
Then $\psi(z,w)=\psi(z',w')=|z_1'|^2+|w_1'|^2+\dots +|w_n'|^2$ where $z'=(z_1',0\dots,0)$ and $w'=(w_1',\dots,w_n')$.
Since $(z',w') \in \mu^{-1}_{\mathbb{R}}(0)$, we can further simplify the formula as $\psi(z',w')=2|z_1'|^2$.
From this equation, one may check that over the locus $\mu^{-1}_{\mathbb{R}}(0) \cap \mu_{\mathbb{C}}^{-1}(1)$, the only critical points of $\psi$ are global minima which occur when $w=\bar{z}$. Therefore, the set of global minima is diffeomorphic to $S^{2n-1}$.
After taking the quotient by $S^1$, we see that $\psi$ descends to an exhausting Morse-Bott K\"ahler potential whose only critical manifold is $\mathbb{CP}^{n-1}$. Therefore, the backward Liouville vector field (which is to say the gradient vector field) contracts the hyperK\"ahler quotient to $\mathbb{CP}^{n-1}$. It follows that the hyperK\"ahler quotient is symplectomorphic to $T^*\mathbb{CP}^{n-1}$.
\end{proof}

We now furthermore note that the map $\chi^{((n-1)1)} \co \cZ \rightarrow \CC$ is regular away from the origin, so given a K{\"a}hler metric on $\cZ$ using the construction of Section~\ref{sec:parallel} we can consider parallel transport with respect to $\chi^{((n-1)1)}$ on $\cZ \setminus \{0\}$. For this construction we only in fact need naive parallel transport. Given $\alpha \in \CC^*$, let $\gamma \co [0,1) \rightarrow \CC^*$ be given by $\gamma = (1-t)\alpha$. Let $h_{\gamma}$ be the naive parallel transport map. Then the vanishing projective cycle associated to $\gamma$ is
\[ \eL_{\alpha} = \{y \in \cZ_{\alpha} : h_{\gamma}|_{[0,t]} \text{ is defined near } y \text { for all } t<1, \ h_{\gamma}|_{[0,t]}(y)\rightarrow 0 \text{ as } t \rightarrow 1\}.\]

This is a Lagrangian submanifold. If the K{\"a}hler form $\Omega$ is the pullback of the standard form on $\mathbb C^{n^2}$ under the usual inclusion $\cZ \hookrightarrow \mathfrak{gl}(n,\CC)\simeq \CC^{n^2}$, it admits a particularly simple description, as follows. For $\alpha \in \CC^*$, consider the diagonal matrix $D_{\alpha} = \mathrm{diag}(\alpha,\cdots, \alpha, (1-n)\alpha) \in \cZ_{\alpha}$. Manolescu shows \cite[Lemmas 4.5 and 4.8] {Manolescu:sln} that we have 
\[ \eL_{\alpha} = \{UD_{\alpha}U^{-1} : U \in U(n)\} \simeq U(n)/(U(n-1)\times U(1)) \simeq \CC \PP^{n-1}.\]

Note that on $\cY_{1,n, \tau}$, the pullback of the standard form on $\mathbb C^{n^2}$ is different from the K\"ahler form in Lemma \ref{lem:cotangent}.
Indeed, for $z,w \in \CC^n$ and the matrix $A:=z^Tw-\frac{1}{n}(z \cdot w)I_n$,  the K\"ahler potential from  $\mathbb C^{n^2}$ is $\psi_{std}:=\mathit{Tr}(A^*A)$  while the 
K\"ahler potential $\psi$ from Lemma \ref{lem:cotangent} is $|z|^2+|w|^2$.
However, we have the following.

\begin{lemma}
Let $\alpha \in \RR \setminus \{0\}$ and  $\tau=(\alpha, (1-n)\alpha)$.  
 There is a $U(n)$-equivariant symplectomorphism $(\cY_{1,n, \tau}, -dd^c \psi_{std}) \to (\cY_{1,n, \tau}, -dd^c \psi)$ preserving $\eL_{\alpha}$.
Moreover, the symplectomorphism \eqref{eq:phi_symp} sends  $\eL_{\alpha}$ to the zero section of the cotangent bundle.
\end{lemma}

\begin{proof} 

The isomorphism 
\[
f:T^*\mathbb{C}^n//S^1 \to \frsl(n)^{((n-1)1)}, \quad f(z,w)=z^Tw-\frac{1}{n}(z \cdot w)I_n
\]
is $U(n)$-equivariant with respect to the actions $U \cdot (z,w)=(zU^T, wU^{-1})$ on the domain and $U \cdot A=UAU^{-1}$ on the codomain for $U \in U(n)$.
Note also that $\psi(z,w)$ and $\psi_{std}(A)$ are $U(n)$-invariant.
Therefore, $\psi_t:=tf_*\psi+(1-t)\psi_{std}$ is a family of $U(n)$-invariant K\"ahler potentials inducing a family of $U(n)$-equivariant K\"ahler forms. 

Similarly to the computation in Lemma~\ref{lem:cotangent}, for any $A=z^Tw-\frac{1}{n}(z \cdot w)I_n$ we can find $U \in U(n)$ such that  for $(z',w'):=U \cdot (z,w)$, all entries of $z'$ other than the first one are $0$, and hence all entries of $A'=UAU^{-1}=(z')^Tw'-\frac{1}{n}(z' \cdot w')I_n$ other than the diagonal entries and the entries on the first row are $0$.
Using the notation $z'=(z_1',0\dots,0)$ and $w'=(w_1',\dots,w_n')$, we see that over the locus $\mu^{-1}_{\mathbb{R}}(0) \cap \mu_{\mathbb{C}}^{-1}(1)$,
\begin{align*}
\psi_{std}(A)=\psi_{std}(A')&=\left(\left(\frac{n-1}{n}\right)^2+\frac{n-1}{n^2}\right)|z_1'w_1'|^2+\sum_{i=2}^n |z_1'w_i'|^2 \\
&=|z_1'|^2\left(\frac{n-1}{n} |w_1'|^2+ \sum_{i=2}^n |w_i'|^2\right)\\
&=|z_1'|^2\left(|z_1'|^2-\frac{1}{n} |w_1'|^2 \right)\\
&=|z_1'|^2\left(|z_1'|^2-\frac{1}{n|z_1'|^2} \right)
\end{align*}
where the second last equality uses $\mu_{\mathbb{R}}=0$ and the last equality uses $\mu_{\mathbb{C}}=1$.

On the other hand, by Lemma~\ref{lem:cotangent}, we have $f_*\psi(A)=2|z_1'|$.
One may check that the family of $U(n)$-invariant K\"ahler potentials $\psi_t:=tf_*\psi+(1-t)\psi_{std}$ only has critical points along the $U(n)$-orbit of  $z = w = (1, 0, \dots, 0)$ consisting of all pairs $(z,w)$  with $\bar{w}=z$, where $\psi_t$ achieves its global minimum.
Therefore, we can apply Moser's argument to obtain the $U(n)$-equivariant symplectomorphism from $-dd^c \psi_0$ to $-dd^c\psi_1$.

Note that the matrix $D_{\alpha}$,  for $\alpha = -\frac{1}{n}$, is the image of the point $(z,w)$ such that $z = w = (0, \dots, 0, 1)$ under $f$.
We therefore see that the symplectomorphism \eqref{eq:phi_symp} sends  $\eL_{\alpha}$, or more precisely, $f^{-1}(\eL_{\alpha})$, to the zero section of the cotangent bundle.

\end{proof}

More generally, for any K{\"a}hler form $\Omega$ on $\cZ-\{0\}$ which is the pullback of a form on $\CC^{n^2}$ which is real analytic on some neighborhood of zero, for sufficiently small values of $t$, $\eL_\alpha$ is diffeomorphic to $\CC\PP^{n-1}$ \cite[Lemma 4.8]{Manolescu:sln}. Such a K{\"a}hler form is called \emph{effective}.

This construction additionally has a relative version, as follows. 

\begin{proposition} \cite[Section 4.4]{Manolescu:sln}
Let $X$ be a smooth complex manifold. Suppose that $W \simeq X \times \frsl(n)^{((n-1)1)}$, equipped with $\Omega$ a K{\"a}hler form which is the pullback of a K{\"a}hler form on $X \times \frsl(n)$ which is real analytic near $X \times \{0\}$. Consider the map $W \rightarrow \CC$ which maps $(x,A)$ to $\chi^{((n-1)1)}(A)$. Given a compact Lagrangian submanifold $\Lambda \subset X$ and $\alpha \in \CC^*$ with $|\alpha|$ sufficiently small, there is a relative vanishing projective cycle in $W_{\alpha}$ which is diffeomorphic to $\Lambda \times \CC\PP^{n-1}$.
\end{proposition}

\subsubsection{Lagrangians associated to crossingless matchings} \label{sec:construct-lagrangians} We now iterate the procedure described above in order to associate a Lagrangian to a crossingless matching. 

Let $\tau = (\vec{\mu}, \vec{\lambda})$ be an element of $\mathrm{BConf}_m$. We say that $\mathfrak{m}=(\delta_1,\delta_2,\cdots,\delta_m)$ is a crossingless matching with endpoints on $\tau$ if $\mathfrak{m}$ is a collection of pairwise disjoint embedded arcs in $\mathbb{C}$ such that $\delta_{i}(0)=\lambda_i$ and $\delta_i(1)=\mu_{\rho(i)}$ for some permutation $\rho$ of $\{1,2,\cdots,m\}$. Two examples appear on the right-hand side of Figure~\ref{fig:matching_example}. We will associate a Lagrangian $\eL_{\mathfrak{m}}\subset\cY_{m,n,\tau}$ to $\mathfrak{m}$ via an inductive procedure on $m$, as we now explain.

For simplicity of notation, assume $\rho$ is the identity permutation. Let $\bar{\mathfrak{m}}=(\delta_2,\cdots,\delta_m)$ and $\bar{\tau}=(\{\mu_2,\cdots,\mu_m\}, \{\lambda_2,\cdots,\lambda_m\})$.  
We will use the map $\chi^{((n-1)1)}$ in Subsection~\ref{sec:vanishing-projective} iteratively.
Assume inductively that we have already constructed a Lagrangian $\eL_{\bar{\mathfrak{m}}}\subset \cY_{m-1,n,\bar{\tau}}$ corresponding to $\bar{\mathfrak{m}}$; when $m=1$ so that $m-1=0$ we may take both the Lagrangian $\eL_{\bar{\mathfrak{m}}}$ and the symplectic manifold $\cY_{m-1,n,\bar{\tau}}$ to be points.  Denote the midpoint of $\delta_1$ by $\lambda$; that is, let $\lambda=\delta_1(1/2)$. Let
\[\tau_0:=(\{0,\mu_2,\cdots,\mu_m\}, \{0,\lambda_2,\cdots,\lambda_m\}), \quad\text{and}\quad \tau_{\lambda}:=(\{\lambda,\mu_2,\cdots,\mu_m\},\{\lambda,\lambda_2,\cdots,\lambda_m\}).\]
Identifying $\cY_{m-1,n,\bar{\tau}}$ with $\cC_{m,n,\tau_0}$ and using parallel transport in the singular set fibration $\mathcal{C}_{m,n}$ from $\cC_{m,n,\tau_0}$ to $\cC_{m,n,\tau_\lambda}$, the Lagrangian $\eL_{\bar{\mathfrak{m}}}$ gets mapped to a Lagrangian $\eL_{\bar{\mathfrak{m}},\lambda}\subset\mathcal{C}_{m,n,\tau_{\lambda}}$. (When $m=1$ all these identifications are identifications of points).

Consider a disk $D\subset \frh^{\pi}/W^{\pi}$ around $\tau_{\lambda}$ 
consisting of the points $(\{\lambda+z,\mu_2,\cdots,\mu_m\},\{\lambda-(n-1)z,\lambda_2,\cdots,\lambda_m\})$ for $z\in \mathbb{C}$ small. Using the path $\delta_1$, we may construct a new path $\gamma \co [0,1]\to \mathfrak{h}^{\pi}/W^{\pi}$ connecting $\tau$ to $\tau_{\lambda}$ 
by setting  
\[\gamma(t)=(\{\delta_1(1-t/2),\mu_2,\cdots,\mu_m\},\{\delta_1(t/2),\lambda_2,\cdots,\lambda_{m}\}).\]
Assume that $\delta_1$ is linear near its midpoint such that for $t$ small, $\gamma(1-t)$ lies on the disk $D$. By \cite[Lemma 3.5]{Manolescu:sln}, we have a commutative diagram 
\begin{equation}\label{eq:local-isom-1}
	\begin{tikzcd}[ column sep=large]
		(\chi^{\pi})^{-1}(D)\cap \cS_{m,n} \arrow[rr, "\simeq \text{ near } \mathcal{C}_{m,n,\tau_{\lambda}}"] \arrow[d, "\chi^{\pi}"]& &
		\mathcal{C}_{m,n,\tau_{\lambda}}\times \mathfrak{sl}(n)^{((n-1)1)} \arrow[d,"\chi^{((n-1)1)}"]\\
		D \arrow[rr,"z"] & &
		\mathbb{C} 
	\end{tikzcd}
\end{equation}
where the top horizontal map is a local isomorphism near $\cC_{m,n,\tau_{\lambda}}$, denoted by $\Upsilon$. We will briefly review the definition of $\Upsilon$, and refer the reader to the proof of \cite[Lemma 3.5]{Manolescu:sln} for more details. Note that $\cC_{m,n,\tau_{\lambda}}$ consists of the matrices $Y\in \cY_{m,n,\tau_{\lambda}}$ where the eigenvalue $\lambda$ has an $n$-dimensional eigenspace. In Manolescu's notation, we have $\cC_{m,n,\tau_{\lambda}} = \cS_{m,n} \cap \cO^{sub}$, where $\cO^{sub}$ is the subregular adjoint orbit in the fibre $(\chi^{\pi})^{-1}(\tau_{\lambda})$ consisting of matrices with $n$ independent $\lambda$ eigenvectors.

Consider the holomorphic subbundle $T\cC_{m,n,\tau_{\lambda}}$ of $T\cS_{m,n}$ over $\cC_{m,n,\tau_{\lambda}}$. Since, $\cC_{m,n,\tau_{\lambda}}$ is affine, there is a holomorphic subbundle $\cR$  of $T\cS_{m,n}$ such that $\cR_Y$ is complementary to $T_Y\cC_{m,n,\tau_{\lambda}}$ for any $Y\in \cC_{m,n,\tau_{\lambda}}$. Then $Y+\cR_Y$ is a transverse slice to the adjoint orbit $\cO(Y)$, and the intersection $(Y+\cR_Y)\cap \frg^{\pi}$ forms a tubular neighborhood of $\cC_{m,n,\tau_{\lambda}}$ in $\cS_{m,n}\cap\frg^{\pi}$. The local isomorphism $\Upsilon$ is the composition of two maps. The first map is the local canonical isomorphism between $(Y+\cR_Y)\cap \frg^{\pi}$ and $\cS^{ss}_Y\cap \frg^{\pi}$, where $\cS^{ss}_Y:=Y + \frg_Y=\{Y+A:[A,Y]=0\}$ is the standard transverse slice to the adjoint orbit $\cO(Y)$ of the semi-simple element $Y$. This isomorphism is given by sending any element of the form $Y+R\in (Y+\cR_Y)\cap \frg^{\pi}$ to a sum of the form $Y+A_R\in \cS^{ss}_Y\cap \frg^{\pi}$ with the property that $Y+A_{R}=g(Y+R)g^{-1}$ for the unique $g\in \exp(V) \subset \mathrm{GL}(mn)$, where $V$ is a small neighborhood of $0$ inside $\operatorname{ad}(\frg) Y$. By \cite[Equation (8)]{Manolescu:sln} on a small neighborhood of $\cC_{m,n,\tau_{\lambda}}$ this map is well-defined and gives a canonical local isomorphism as claimed. This first isomorphism depends on the choice of the subbundle $\cR$ but the Lagrangian constructed later does not (up to Hamiltonian isotopy).

To define the second map, recall from Lemma~\ref{lem:injectivity} that the map $p \co \CC^{mn} \to \CC^n$ which projects onto the first $n$ coordinates is an injection on the eigenspace $E_{Y}(\lambda)$. Let $e_i$ be the $i$th standard basis vector in $\CC^n$, and consider a basis $B=\{v_1,v_2,\cdots,v_n\}$ for the $E_{Y}(\lambda)$ such that the $p(v_i)=e_i$. Now consider the element $A_R$, and recall that since $A_R\in\cS^{ss}_Y$ we have $[Y,A_R]=0$. So $A_{R}$ acts on $E_{Y}(\lambda)$. Let the matrix of this action in the basis $B$ by $\widetilde{A}_R$. We let the second map send $Y+A_R$ to the pair $(Y,\widetilde{A}_R)$.
Note that $\widetilde{A}_R$ is traceless because for any $(\vec{\mu}, \vec{\lambda}) \in D$, the sum $\sum (n-1)\mu_i+ \sum \lambda_i$ is the same.

Having established the isomorphism $\Upsilon$, we define the Lagrangian $\eL_{\mathfrak{m}}$ as follows.  Choose a small $t_0$ such that $\gamma(1-t)\in D$ for all $0\le t\le t_0$ and consider the Lagrangian \[\eL_{z(\gamma(1-t_0))}\subset (\chi^{((n-1)1)})^{-1}(z(\gamma(1-t_0)))\] obtained from the vanishing projective space construction. Then, the isomorphism from \eqref{eq:local-isom-1} maps $\eL_{\bar{\mathfrak{m}},\lambda}\times \eL_{z(\gamma(1-t_0))}$ to a Lagrangian in $\cY_{m,n,\gamma(1-t_0)}$. Finally we may use parallel transport along $\gamma$ to get the Lagrangian $\eL_{\mathfrak{m}}$ in $\cY_{m,n,\tau}$.

Alternatively, we can use a simultaneous collapse instead of induction to define $\eL_{\mathfrak{m}}$.  By \cite[Lemma 3.7]{Manolescu:sln}, we have a commutative diagram 
\begin{equation}\label{eq:local isom}
	\begin{tikzcd}[ column sep=large]
		(\chi^{\pi})^{-1}(P)\cap \cS_{m,n} \arrow[rr, "\simeq \text{ near } \mathcal{O}^{min} \cap  \cS_{m,n}"] \arrow[d,"\chi^{\pi}"]& &
		(\mathcal{O}^{min} \cap  \cS_{m,n})\times \left(\mathfrak{sl}(n)^{((n-1)1)}\right)^\ell \arrow[d,"\left(\chi^{((n-1)1)}\right)^\ell"]\\
		P \arrow[rr,"z"] & &
		\mathbb{C}^\ell 
	\end{tikzcd}
\end{equation}
generalizing \eqref{eq:local-isom-1}. The notation is as follows: firstly, $P\subset \frh^{\pi}/W^{\pi}$ is a polydisk consisting of the set of points
 \[(\{\lambda_1+z_1,\cdots,\lambda_\ell+z_\ell,\mu_{\ell+1},\cdots, \mu_m\}, \{\lambda_1-(n-1)z_1,\dots,\lambda_\ell-(n-1)z_\ell,\lambda_{\ell+1},\dots, \lambda_m\})\] 
 for sufficiently small $(z_1,\cdots,z_\ell)\in \CC^\ell$. Moreover, letting $P_0\in P$ be the point with $z_1=z_2=\cdots=z_\ell=0$, then  $\cO^{min}\subset (\chi^{\pi})^{-1}(P_0)$ is the singular stratum consisting of the matrices $Y$ with $\dim E_{Y}(\lambda_i)=n$ for all $i=1,\cdots, \ell$. If $\ell=1$ this map is precisely the same as \eqref{eq:local-isom-1}, so for simplicity we denote the new isomorphism by $\Upsilon$ as well. 

In particular, when $\ell=m$, the intersection $\mathcal{O}^{min} \cap  \cS_{m,n}$ is a point and the righthand column of \eqref{eq:local isom} is the $m$-fold product of 
$\chi^{((n-1)1)}\co \mathfrak{sl}(n)^{((n-1)1)} \to \mathbb{C}$. To see this, notice that if $Y\in \mathcal{O}^{min}\cap\mathcal{S}_{m,n}$ and the component $n \times n$ matrices which form the first column of $Y$ are denoted by $Y_{11}, Y_{21},\dots, Y_{m1}$ then 
\[\lambda^{m}I_n-\lambda^{m-1}Y_{11}-\lambda^{m-2}Y_{21}-\cdots-Y_{m1}=0\]
for $\lambda=\lambda_1,\lambda_2,\dots,\lambda_m$. This is a linear system of equations with an invertible coefficients matrix, and thus has a unique solution. 

Given this isomorphism, one may construct a Lagrangian by carrying out the vanishing projective cycles construction in each factor on the righthand side. Indeed, by the following lemma the results of the iterated construction and the simultaneous construction agree up to Hamiltonian isotopy.

\begin{lemma}\cite[Lemmas 4.12 and 7.4]{Manolescu:sln}\label{l:productLag}
	Let $\fkm=(\delta_1,\cdots,\delta_m)$ be a perfect matching with endpoints on $\tau\in P$ such that for any $t\in [0,1]$, $\lambda_i=\delta_i(t/2)$ and $\mu_i=\delta_i(1-t/2)$, then if $\ell=m$ in the local isomorphism \eqref{eq:local isom}, the configurations $\tau_t=(\vec{\mu}_t, \vec{\lambda}_t)$  lie in the polydisk $P$. It follows that the Lagrangian $\eL_{\fkm}$ is Hamiltonian isotopic to the Lagrangian identified with the product $(\CC \PP^{n-1})^m$ of vanishing projective spaces of $\mathfrak{sl}(n)^{((n-1)1)} \to \mathbb{C}$ under the isomorphism \eqref{eq:local isom}.
\end{lemma}

\subsection{Braid closures and the definition of the invariant} \label{sec:defn}

\begin{figure}
\begingroup%
  \makeatletter%
  \providecommand\color[2][]{%
    \errmessage{(Inkscape) Color is used for the text in Inkscape, but the package 'color.sty' is not loaded}%
    \renewcommand\color[2][]{}%
  }%
  \providecommand\transparent[1]{%
    \errmessage{(Inkscape) Transparency is used (non-zero) for the text in Inkscape, but the package 'transparent.sty' is not loaded}%
    \renewcommand\transparent[1]{}%
  }%
  \providecommand\rotatebox[2]{#2}%
  \newcommand*\fsize{\dimexpr\f@size pt\relax}%
  \newcommand*\lineheight[1]{\fontsize{\fsize}{#1\fsize}\selectfont}%
  \ifx\svgwidth\undefined%
    \setlength{\unitlength}{343.46395111bp}%
    \ifx\svgscale\undefined%
      \relax%
    \else%
      \setlength{\unitlength}{\unitlength * \real{\svgscale}}%
    \fi%
  \else%
    \setlength{\unitlength}{\svgwidth}%
  \fi%
  \global\let\svgwidth\undefined%
  \global\let\svgscale\undefined%
  \makeatother%
  \begin{picture}(1,0.46113075)%
    \lineheight{1}%
    \setlength\tabcolsep{0pt}%
    \put(0,0){\includegraphics[width=\unitlength,page=1]{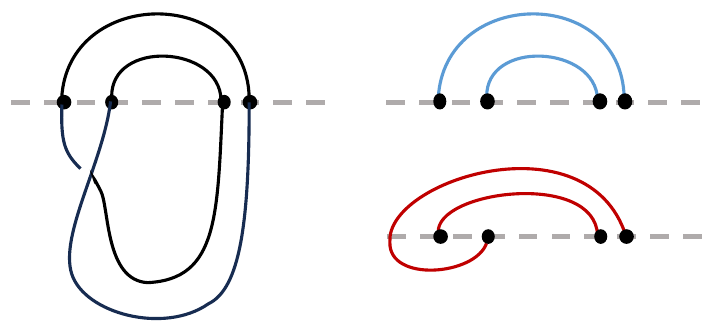}}%
    \put(0.82650902,0.21141675){\color[rgb]{0,0,0}\makebox(0,0)[lt]{\lineheight{1.25}\smash{\begin{tabular}[t]{l}$\mathfrak{m}_\beta$\end{tabular}}}}%
    \put(0.35512368,0.25441697){\color[rgb]{0,0,0}\makebox(0,0)[lt]{\lineheight{1.25}\smash{\begin{tabular}[t]{l}$L$\end{tabular}}}}%
    \put(0.85424035,0.41077739){\color[rgb]{0,0,0}\makebox(0,0)[lt]{\lineheight{1.25}\smash{\begin{tabular}[t]{l}$\mathfrak{m}_0$\end{tabular}}}}%
  \end{picture}%
\endgroup%

\caption{On the left, a link $L$ is represented as the closure of a braid. On the right, we decompose the link into two matchings, one of which is the standard ``rainbow matching'' $\mathfrak{m}_0$ and the other of which determines the link. The endpoints of the matchings lie on $\tau_{std}=\left(
\vec{\mu}=
\left\{\frac{1}{n-1}, \frac{2}{n-1}\right\}, \vec{\lambda}=\{-1,-2\}
\right)$.}
\label{fig:matching_example}
\end{figure}

We are finally ready to state the definition of Manolescu's link invariant. Let $L$ be a link and express $L$ as the closure of an $m$-strand braid $\beta$. Let $\tau_{std}= \left (\vec{\mu}=\left\{\frac{1}{n-1}, \cdots, \frac{m}{n-1}\right\}, \vec{\lambda} =\left\{-1, \dots, -m+1, -m\right\} \right )$ be a standard element of $\mathrm{BConf}_m$.  We use the braid representative $\beta$ to write $L$ as a union of two crossingless matchings $\mathfrak{m}_0$ and $\mathfrak{m}_\beta$ of $\tau_{std}$. To do this, we take $\mathfrak{m}_0$ to be the standard rainbow crossingless matching which consists of $m$ paths beginning on $-i$ and ending on $\frac{i}{n-1}$, as in Figure~\ref{fig:matching_example}. We set $\mathfrak{m}_\beta$ to be the crossingless matching determined by the bottom half of the link, such that allowing the paths in $\mathfrak{m}_\beta$ to undercross those in $\mathfrak{m}_0$ produces a flattening of the closure of $\beta$.  We let $\eL_0=\eL_{\mathfrak{m}_0}$ be the \emph{standard Lagrangian}, and $\eL_\beta= \eL_{\mathfrak{m}_\beta}$ be the \emph{braid Lagrangian}.

We are now prepared to state the definition of Manolescu's symplectic $\frsl(n)$ homology.

\begin{definition} Let $L$ be a link. The symplectic $\frsl(n)$ homology is 
\[ H_{n,symp}^*(L; \mathbb F_2) = HF^{*+(n-1)(m+w)}(\cY_{m,n,\tau_{std}}, \eL_0, \eL_\beta)\]
where $w$ is the writhe of $\beta$,  $\eL_0$ is the standard Lagrangian and $\eL_\beta$ is the braid Lagrangian associated to an $m$-strand braid representative of $L$.
\end{definition}

Manolescu shows that this theory is a link invariant \cite[Section 8]{Manolescu:sln}.

\section{Intrinsic involutions on symplectic Khovanov--Rozansky homology} \label{sec:intrinsic}

In this section, we study an involution on $\cY_{m,n,\tau}$  which, when combined with Theorem~\ref{thm:localization-intro}, will lead to Theorem \ref{thm:mainss}.

\subsection{The ambient symplectic manifold and intrinsic symmetries} \label{sec:ambient-intrinsic}

Let $t_k$ be the $n \times n$ diagonal matrix 
\[
t_k:=diag(-1, \dots, -1, 1,\dots,1)
\]
where $-1$ appears $k$ times.
Let $T_k$ be the $mn \times mn$ diagonal matrix $diag(t_k, \dots, t_k)$, where the matrix $t_k$ appears $m$ times.
We consider the involution
\begin{equation} \label{eq:intrinsic}
\iota_k \co \frg \to \frg, \quad x \mapsto T_kxT_k^{-1}=T_k x T_k
\end{equation}
It is immediate to check the following.

\begin{lemma}\label{lem:fixedptslice}
The involution preserves $\cS_{m,n}$, and therefore descends to an involution on $\cS_{m,n}$.
Moreover, the fixed point set is 
\[
\cS_{m,n}^{\fix}=\left\{  \begin{pmatrix} Y_{11}^a & 0 & I_k & & && &\\
0&Y_{11}^b & 0& I_{n-k} &  & &&\\
Y_{21}^a &0& 0& 0 & I_k & &&\\
0&Y_{21}^b &0 &0 &0 &I_{n-k}&&\\
\dots & &&  & &\dots   &\\
Y_{m1}^a & 0&& & && 0 &\\
0& Y_{m1}^b& 0 & 0 & 0&0&0 &0
\end{pmatrix} : \begin{split} Y_{j1}^a\in \frgl(k), Y_{j1}^b \in \frgl(n-k) \text{ for }1\le j\le m\end{split}   \right\} \]
\end{lemma}


The involution can be lifted to an involution on $\tilde{\frg}^{\pi}$ given by
\[
\tilde{\iota}_k \co \tilde{\frg}^{\pi} \to \tilde{\frg}^{\pi} \quad  (x,F) \mapsto (T_kxT_k^{-1}, T_kF)
\]
where $T_kF$ is the partial flag such that $(T_kF)_j:=T_kF_j$.
Since $\tilde{\chi}^{\pi}$ is conjugation invariant, it is $\tilde{\iota}_k$-invariant.
As a consequence, we get the following.

\begin{lemma}
For any partition $\pi$ of $mn$, the involution $\iota_k$ descends to an involution on $\frg^{\pi}$, with respect to which the map $\chi^{\pi}$ is $\iota_k$-invariant. 
\end{lemma}

We now describe the fixed point set $\cY_{m,n,\tau}^{\fix}$ when $\tau \in \frh^{\pi,reg}/W^{\pi}$. First, we introduce some definitions.
Let $a \le c$, and furthermore let $\pi'=(\pi_1',\dots,\pi'_{s'})$ be a partition of $a$ and $\pi=(\pi_1,\dots,\pi_{s})$ be a partition of $c$. An embedding of $\pi'$ to $\pi$ is an injective function $f:\{1,\dots,s'\} \to \{1, \dots,s\}$  such that $\pi'_i \le \pi_{f(i)}$ for all $i$. We also require that if $i<j$ and $\pi_i'=\pi_j'$, then we have $f(i)<f(j)$. Moreover, corresponding to a partition $\pi'$ embedded in $\pi$, let $\pi\setminus \pi'$ be the partition of $c-a$ whose entries are $\pi_{f(i)}-\pi_i'$ for $i\in\{1,2,\cdots,s'\}$ and $\pi_j$ for $j\notin \mathrm{Im}(f)$, with zero entries removed. By definition, this partition comes with a natural embedding in $\pi$.

Let $\alpha_i$ be the pairwise distinct eigenvalues associated to $\tau$ so that $\alpha_i$ has multiplicity $\pi_i$.
There are $|W^{\pi}|$ many ways to label the subscripts of $\alpha_i$ so that $\alpha_i$ has multiplicity $\pi_i$.
We fix one such choice; this is equivalent to fixing a lift of $\tau$ to $\sigma \in \frh^{\pi,reg}$, thought of as an ordered list of eigenvalues $\alpha_i$ with multiplicities $\pi_i$.

Given an embedding $f$, we let $(f^*\alpha)_i:=\alpha_{f(i)}$. If we then regard $(f^*\alpha)_i$ as an eigenvalue of multiplicity $\pi_i'$ for all $i \in \{1,\dots,s'\}$, this specifies an element $\sigma' \in \frh^{\pi',reg}$ and correspondingly an element 
$\tau' \in \frh^{\pi',reg}/W^{\pi'}$.
We call $\sigma'$ the pullback of $\sigma$ under $f$, and $\tau'$ the pullback of $\tau$ under $f$.

Now let $\pi$ be a partition of $mn$ and $\cX_{k,\pi}$ be the set of all embeddings $(\pi',f)$
from a partition $\pi'$ of $m(n-k)$ into $\pi$. 

\begin{lemma}
For any $\tau \in \frh^{\pi,reg}/W^{\pi}$ the fixed point set $\cY_{m,n,\tau}^\fix$ under the involution \eqref{eq:intrinsic} is isomorphic to 
\[
\amalg_{(\pi',f)\in\cX_{k,\pi}}\cY_{m,k,\tau_{\pi\setminus\pi'}}\times\cY_{m,n-k,\tau_{\pi'}}
\]
where  $\tau_{\pi'} \in \frh^{\pi',reg}/W^{\pi'}$ and $\tau_{\pi \setminus \pi'} \in \frh^{\pi\setminus \pi',reg}/W^{\pi \setminus \pi'}$ denote the pullback of $\tau$ under $f$ and under the induced embedding of $\pi\setminus \pi'$ in $\pi$, respectively.
\end{lemma}

\begin{proof}
The fixed point set $\cY_{m,n,\tau}^{\fix}$ is the intersection $\cS_{m,n}^{\fix}\cap (\chi^{\pi})^{-1}(\tau)$ of the fixed point set of $\cS_{m,n}$ with the preimage of $\tau$ under $\chi^{\pi}$.  For any $Y\in\cS_{m,n}^{\fix}$ as in Lemma \ref{lem:fixedptslice} let
\[Y^a=\begin{pmatrix}Y_{11}^a&I_k&&&\\
Y_{21}^a&0&I_k&&\\ 
\cdots&&&\cdots&\\
 Y_{(m-1)1}^a&0&0&0&I_k\\
  Y_{m1}^a&0&0&0&0
  \end{pmatrix}  \quad\text{and}\quad Y^b=\begin{pmatrix}Y_{11}^b&I_{n-k}&&&\\
Y_{21}^b&0&I_{n-k}&&\\ 
\cdots&&&\cdots&\\
 Y_{(m-1)1}^b&0&0&0&I_{n-k}\\
  Y_{m1}^b&0&0&0&0
  \end{pmatrix}.\]
Note that $\chi(Y)=\chi(Y^a)\chi(Y^b)$, where $\chi$ denotes the characteristic polynomial,  and both $Y^a$ and $Y^b$ are diagonalizable. As above, we fix an ordering $\sigma$ of $\tau$. For each $\alpha_i$ an element of $\sigma$, let $\pi'_i$ be its multiplicity as an eigenvalue of $Y^b$, and let $\pi'$ be the partition of $m(n-k)$ formed by concatenating the nonzero values $\pi'_i$. Likewise let $\tau_{\pi'}$ be the unordered configuration consisting of the elements $\alpha_i$ for which $\pi'_i$ nonzero. There exists an ordered lift $\sigma_{\pi'}$ of $\tau_{\pi'}$ determined by preserving the relative order of the eigenvalues in $\sigma$, and a natural embedding $f \co \pi' \rightarrow \pi$ with the property that $f^{*}(\sigma)=\sigma_{\pi'}$. With this in mind, $f^*(\tau) = \tau_{\pi'} = \chi^{\pi'}(Y^b)$ and so $Y^b\in\cY_{m,n-k,\tau_{\pi'}}$. Similarly, $Y^a\in\cY_{m,k,\tau_{\pi\setminus\pi'}}$.  Furthermore, we may reverse this procedure: having fixed an ordered lift $\sigma$ of $\tau$, we see that for any given embedding $(\pi',f)$, if we construct the pullbacks $\sigma_{\pi'}$ and $\sigma_{\pi \setminus \pi'}$ as described above, this specifies a copy of $\cY_{m,k,\tau_{\pi\setminus\pi'}}\times\cY_{m,n-k,\tau_{\pi'}}$ which lies inside of $\cY_{m,n,\tau}^\fix$.
\end{proof}

As an aid to intuition, we now examine the decomposition of the fixed set into disjoint components in a few examples.

\begin{example} We let $n=4$ and $m=1$, so that we may take the ambient symplectic manifold to be $\cY_{m,n,\tau}^{\pi}$ for $\pi = (3,1)$ and $\tau = (\mu_1,\lambda_1)$, where  $\mu_1$ is an eigenvalue of multiplicity three and $\lambda_1$ is an eigenvalue of multiplicity one. We let $k=2$, so that we are interested in embeddings of partitions $\pi'$ of $m(n-k)= 2$ into $\pi$. The possible embeddings are as follows.
\begin{itemize}

\item The partition $\pi'=(2)$ embeds into $\pi = (3,1)$ in exactly one way, along the map $f \co \{1\} \rightarrow \{1,2\}$ given by $f(1)=1$. In this case we have $\pi \setminus \pi' = (1,1)$ and
\[\tau_{\pi \setminus \pi'} = (\mu_1, \lambda_1)  \qquad \qquad   \tau_{\pi'} = (\mu_1). \]

\item The partition $\pi'=(1,1)$ embeds into $\pi=(3,1)$ in exactly one way, along the map $f \co \{1,2\} \rightarrow \{1,2\}$ given by the identity. In this case we have $\pi \setminus \pi' = (2)$ and
\[
\tau_{\pi \setminus \pi'} = (\mu_1)  \qquad \qquad   \tau_{\pi'} = (\mu_1,\lambda_1).
\]
\end{itemize}
The fixed set therefore has two connected components. We will describe these components geometrically in Lemma \ref{l:m=1model}.
\end{example}

\begin{example} \label{example:n5m2} We let $n=5$ and $m=2$, so that we may take the ambient symplectic manifold to be $\cY_{m,n,\tau}$ for $\pi= \pi_{std} = (4,4,1,1)$ and $\tau = (\mu_1, \mu_2, \lambda_1, \lambda_2)$. We let $k=3$ and $n-k=2$, so that we are interested in embeddings of partitions $\pi'$ of $m(n-k) = 4$ into $\pi$. Before enumerating the possible embeddings, we recall that by Corollary~\ref{cor:too-big-is-trivial}, the fixed component corresponding to an embedding $(\pi',f)$ will necessarily be empty if any entry of $\pi'$ is greater than $n-k=2$ or any entry of $\pi \setminus \pi'$ is greater than $k=3$. We therefore list only those embeddings that result in nonempty components, which are as follows.
\begin{itemize}

\item The partition $\pi' = (2,2)$ embeds into $\pi = (4,4,1,1)$ in exactly one way, along the map $f_1 \co \{1,2\} \rightarrow \{1,2,3,4\}$ given by $f_1(1)=1$ and $f_1(2)=2$. In this case we have $\pi \setminus \pi' = (2,2,1,1)$ and 
\[
\tau_{\pi \setminus \pi'} = (\mu_1, \mu_2, \lambda_1, \lambda_2)  \qquad \qquad   \tau_{\pi'} = (\mu_1, \mu_2).
\]
This component of the fixed set is of the form $\cY^{\pi_{std}}_{2,3,\tau} \times \{\mathrm{pt}\}$.

\item The partition $\pi' = (1,1,1,1)$ embeds into $\pi = (4,4,1,1)$ in exactly one way, along the identity map $f_2 \co \{1,2,3,4\}\rightarrow \{1,2,3,4\}$. In this case we have $\pi \setminus \pi' = (3,3)$ and 
\[
\tau_{\pi \setminus \pi'} = (\mu_1, \mu_2)  \qquad \qquad   \tau_{\pi'} = (\mu_1, \mu_2, \lambda_1, \lambda_2).
\]
This component of the fixed set is of the form $\{\mathrm{pt}\} \times \cY^{\pi_{std}}_{2,2,\tau}$.

\item There are four embeddings of the partition $\pi' = (2,1,1)$ into $\pi = (4,4,1,1)$ which result in a nonempty fixed component. In all four cases the complementary partition is $\pi \setminus \pi' = (3,2,1)$. However, there are nevertheless four distinct connected components, which we now discuss. In the below, we list elements of $\tau_{\pi'}$ so that the first entry is an eigenvalue of multiplicity two in the lift of the configuration to an equivalence class of diagonal matrices in $\frh^{\pi',reg}$ and the remaining two entries have multiplicity one. Likewise for $\tau_{\pi \setminus \pi'}$ the first entry is the eigenvalue of multiplicity three, the second entry is an eigenvalue of multiplicity two, and the final entry is an eigenvalue of multiplicity one.

\noindent The first such embedding is along $f_3 \co \{1,2,3\} \rightarrow \{1,2,3,4\}$ where $f_3(1)=1$, $f_3(2)=2$, and $f_3(3)=3$, in which case we have
\[
\tau_{\pi'} = (\mu_1, \mu_2, \lambda_1)  \qquad \qquad   \tau_{\pi \setminus \pi'} = (\mu_2, \mu_1,\lambda_2).
\]
The second such embedding is along the map given by $f_4(1)=1$, $f_4(2) =2$, $f_4(3)=4$, in which case we have
\[
\tau_{\pi'} = (\mu_1, \mu_2, \lambda_2)  \qquad \qquad   \tau_{\pi \setminus \pi'} = (\mu_2, \mu_1,\lambda_1).
\]
The third such embedding is along the map given by $f_5(1)=2$, $f_5(2)=1$, $f_5(3)=3$. In this case we have
\[ \tau_{\pi'} = (\mu_2, \mu_1, \lambda_1)  \qquad \qquad   \tau_{\pi \setminus \pi'} = (\mu_1, \mu_2, \lambda_2).\]
Finally the fourth such embedding is along the map given by $f_6(1)=2$, $f_6(2)=1$, $f_6(3)=4$. In this case we have
\[ \tau_{\pi'} = (\mu_2, \mu_1, \lambda_2)  \qquad \qquad   \tau_{\pi \setminus \pi'} = (\mu_1, \mu_2, \lambda_1).\]

\end{itemize}
We see that in the case $n=5, m=2, k=3$ the fixed set has six nonempty connected components. The reader can easily check that whenever $\pi=\pi_{std}$ the total number of components in fact depends only on $m$, in particular it equals six whenever $m=2$. 
\end{example}

Finally, we observe that we can also describe the involution and the fixed set geometrically when $m=1$. Despite its limitations this perspective is often an important aid to intuition, as we will see in the next section. 

\begin{lemma}\label{l:m=1model} The symplectomorphism
	\[
	\phi:\cY_{1,n,\tau}\to (T^*\mathbb{CP}^{n-1},\omega_{can})
	\]
of Lemma~\ref{lem:cotangent} is $\ZZ/2\ZZ$-equivariant with respect to the action $\iota_k$ on $\cY_{1,n, \tau}$ and the action on $T^*\mathbb{CP}^{n-1}$ induced by the action on $\mathbb{CP}^{n-1}$ of the form
	\begin{align}\label{eq:actionCP}
		[z_0:\dots:z_{n-1}] \mapsto [-z_0:\dots:-z_{k-1}:z_k: \dots :z_{n-1}] 
	\end{align}
In particular, $\cY_{1,n, \tau}^{\fix}$ consists of two components.
\end{lemma}

\begin{proof} We continue with the notation of Lemma~\ref{lem:cotangent}. Recall that we have the isomorphism
\[
	f:T^*\mathbb{C}^n//S^1:=\mu^{-1}_{\mathbb{R}}(0)/S^1 \to \frsl(n)^{((n-1)1)}
	\]
which has the property that $\chi^{((n-1)1)} \circ f = -\frac{1}{n} z \cdot w$.

	To compare the $\mathbb{Z}/2\mathbb{Z}$ actions, we consider the action on $T^*\mathbb{C}^n$ given by
	\begin{align}\label{eq:actionC2n}
		(z_0,\dots,z_{n-1},w_0,\dots,w_{n-1}) \mapsto (-z_0,\dots,-z_{k-1},z_k, \dots ,z_{n-1}, -w_0,\dots,-w_{k-1},w_k, \dots ,w_{n-1})  
	\end{align}
	and the $\iota_k$ action on $\frsl(n)^{((n-1)1)}$.
	It is clear from the formula that $f$ is an equivariant isomorphism. The action \eqref{eq:actionC2n} descends to an involution on $T^*\mathbb{CP}^{n-1}=(\mu_{\mathbb{R}}^{-1}(0) \cap \mu_{\mathbb{C}}^{-1}(1))/S^1$ which is compatible with the K\"ahler potential $\psi$ of Lemma~\ref{lem:cotangent} and its gradient flow. Moreover, the involution on the zero section is given by  \eqref{eq:actionCP}, and therefore agrees with the action on $T^*\mathbb{CP}^{n-1}$ induced by  \eqref{eq:actionCP}.
\end{proof}

\subsection{Equivariant neighborhoods, parallel transport, and the fixed Lagrangians}\label{sec:lagrangian-involution}

In this section we study the fixed set of the Lagrangian $\eL_{\mathfrak m}$ associated to a crossingless matching $\mathfrak{m}$ under the involution \eqref{eq:intrinsic} defined in the previous subsection.

We begin with a technical lemma. Recall that the construction of the Lagrangian $\eL_{\mathfrak{m}}$ via simultaneous collapse used a local isomorphism, given in \eqref{eq:local isom}, to define coordinates on a polydisk where the collapse occurs.

\begin{lemma}\label{l:equi_local_isom}
	The local isomorphism \eqref{eq:local isom} may be chosen to be equivariant with respect to $\iota_k$ on the left-hand side and the diagonal action given by the $\iota_{k}$ action on $\mathcal{O}^{min}\cap\cS_{m,n}$ and each factor of $\frsl(n)^{((n-1)1)}$ on the right-hand side.
\end{lemma}

To prove Lemma \ref{l:equi_local_isom}, we need to use the following general result:

\begin{lemma}\label{lem:equiv-normal}
Let $X$ be a smooth affine variety with a holomorphic involution $\iota$, $E$ an $\iota$-equivariant holomorphic vector bundle and $F \subset E$ an $\iota$-equivariant holomorphic subbundle.
Then there is an $\iota$-equivariant holomorphic subbundle $R$ of $E$ complementing $F$.
\end{lemma}
 
\begin{proof}
Let $Q:=E/F$.
The existence of $R$ is equivalent to the existence of an $\iota$-equivariant holomorphic splitting of the vector bundles
\[
0 \to F \to E \to Q \to 0.
\]
We will see that obstruction of the equivariant  splitting lies in $H^1(X, \mathcal{H}om(Q,F))^{\iota}$, where $\iota$ acts on the sections $C^{\infty}(X,\mathcal{H}om(Q,E))$ by 
\[(\iota\cdot s)_x:=\iota\circ s_{\iota(x)}\circ\iota:Q_x\to Q_{\iota(x)}\to E_{\iota(x)}\to E_x.\]
Because $X$ is affine, hence Stein, $H^1(X,\mathcal{H}om(Q,F))=0$. So this obstruction is trivial.

To explain why the obstruction lies in $H^1(X, \mathcal{H}om(Q,F))^{\iota}$, first note that we can use a partition of unity to produce an equivariant metric on $E$ and hence a smooth splitting of $E$ into $F$ and the orthogonal complement of $F$.
We can use this splitting to define an $\iota$-invariant section  $s \in C^{\infty}(X,Hom(Q,E))^{\iota}$ with $\pi\circ s=\mathrm{Id}_Q$.
The failure of $s$ being holomorphic is $\beta:= \bar{\partial}_E \circ s- s \circ \bar{\partial}_Q$ which takes values in $\mathcal{H}om(Q,F)$.
We have $\bar{\partial} \beta=0$ so $\beta \in H^1(X, \mathcal{H}om(Q,F))^{\iota}$.
Since $H^1(X, \mathcal{H}om(Q,F))=0$, we can find a primitive $\gamma \in C^{\infty}(X,Hom(Q,F))$ such that $\bar{\partial} \gamma=\beta$. Replacing $\gamma$ with $\frac{1}{2}(\gamma+\iota\cdot \gamma)$ we may assume $\gamma$ is $\iota$-invariant. Then, the section $\tilde{s}:=s-\gamma$ defines an equivariant holomorphic splitting.
\end{proof}

\begin{proof}[Proof of Lemma \ref{l:equi_local_isom}]
	Following the proof of \cite[Lemma 3.7]{Manolescu:sln}, we see that to upgrade \eqref{eq:local isom} to an equivariant isomorphism, we must construct an equivariant neighborhood arising from an equivariant splitting of the tangent bundle along  $\mathcal{O}^{min} \cap  \cS_{m,n}$. By Lemma \ref{lem:equiv-normal}, such an equivariant splitting must exist. Fix an $\iota_k$-equivariant holomorphic subbundle $\cR$. Notice that $\iota_k(\cS^{ss}_Y\cap \frg^{\pi})=\cS^{ss}_{\iota_k(Y)}\cap \frg^{\pi}$ because $\iota_k(Y+A)=T_kYT_k^{-1}+T_kAT_k^{-1}$ and $[Y,A]=0$ if and only if $[T_kYT_k^{-1},T_kAT_k^{-1}]=0$.
Under the first local isomorphism from $(Y+\cR_Y)\cap\frg^{\pi}$ to $\cS^{ss}_Y\cap \frg^{\pi}$, if $Y+R$ is mapped to $Y+A_R$ by conjugation with $g$ then $\iota_k(Y+R)$ will be mapped to $\iota_k(Y+A_R)$ by conjugation with $T_kgT_k^{-1}$. Then, $\iota_k(Y+A_R)$ is mapped to $\iota_k(Y)\in \mathcal{O}^{min} \cap  \cS_{m,n}$ along with the actions of $\iota_k(A_R)$ on $E_{\iota_k(Y)}(\lambda_i)$ under the identification $\mathfrak{sl}(E_{\iota_k(Y)}(\lambda_i))^{((n-1)1)}$ and $\mathfrak{sl}(n)^{((n-1)1)}$ by projection on the first $n$ coordinates for $i=1,\cdots, \ell$. Note that if $B=\{v_1,v_2,\cdots,v_n\}$ is a basis for $E_{Y}(\lambda_i)$ that $p(v_i)=e_i$ then $B'=\{-T_kv_1,\cdots,-T_kv_{k},T_kv_{k+1},\cdots,T_{k}v_n\}$ is the corresponding basis for $E_{\iota_k(Y)}(\lambda_i)$ that projects to the standard basis for $\CC^n$. It is easy to check that the actions of $\iota_k(A_R)$ in the basis $B'$ and $A_R$ in the basis $B$ differ by $t_k$-conjugation.
\end{proof}

We immediately obtain the following corollary.

\begin{corollary}
When $m=\ell$, so that $\mathcal{O}^{min} \cap  \cS_{m,n}$ is a point, there is an equivariant neighborhood $V$ of $\mathcal{O}^{min} \cap  \cS_{m,n}$ with the property that for $\tau$ in the polydisk $P$, the intersection $V \cap \mathcal Y_{m,n,\tau}^{\fix}$ has exactly $2^m$ components corresponding to the components of $((T^*\mathbb{CP}^{n-1})^{\fix})^m$.
\end{corollary}

Notice in particular that given a matching $\mathfrak m$, not all of the connected components of $\cY_{m,n,\tau}^{\fix}$ appear in a neighborhood of the point of simultaneous collapse corresponding to $\mathfrak{m}$.

We now discuss the symplectic structure of the fixed set. Since the involution $\iota_k$ sends every matrix entry to itself up to sign, the plurisubharmonic potential function constructed in \eqref{eq:potential} is clearly $\iota_k$-invariant, and generates an $\iota_k$-invariant symplectic form. Since $\psi$ is $\iota_k$-invariant, it follows immediately that rescaled parallel transport also behaves well with respect to the involution.

\begin{corollary} \label{cor:invariant-parallel}
	The rescaled parallel transport map is $\iota_k$-equivariant.
\end{corollary}

We now turn our attention to the restriction of the involution on the Lagrangians. Happily, our work so far implies that this map has a particularly clean form, which we now describe.

\begin{lemma} \label{lem:lagrangian-involution}
Let $\eL_{\mathfrak{m}}\subset \cY_{m,n,\tau}$ be the Lagrangian corresponding to a crossingless matching $\mathfrak{m}$ with endpoints on $\tau$. Then $\iota_k$ preserves $\eL_{\mathfrak{m}}$ and there is a $\ZZ/2\ZZ$-equivariant diffeomorphism $\phi$ between $\eL_{\mathfrak{m}}$ and $\mathbb{CP}^{n-1}\times ...\times \mathbb{CP}^{n-1}$ equipped with the $\ZZ/2\ZZ$ action $\mathit{diag}(r_k,\cdots,r_k)$ where $r_k:\CP^{n-1}\to\CP^{n-1}$ is given by 
\[r_k([z_0:\cdots:z_{n-1}])=[-z_0:\cdots:-z_{k-1}:z_k:\cdots:z_{n-1}].\]
\end{lemma}

\begin{proof} The discussion of simultaneous collapse in Section~\ref{sec:construct-lagrangians} and Lemma \ref{l:productLag} implies that up to Hamiltonian isotopy $\eL_{\mathfrak{m}}$ may be taken to be the product manifold $(\CC \PP^{n-1})^m$ arising from the product of vanishing projective spaces on the righthand side of the diagram~\eqref{eq:local isom}. This Hamiltonian isotopy arises from (1) rescaled parallel transport, which by Corollary~\ref{cor:invariant-parallel} can be taken to be $\iota_k$-equivariant, and (2) by using the Moser trick to change the strictly plurisubharmonic potential function from the push-forward of the potential function on $\mathcal S_{m,n}$ to $\left(\mathfrak{sl}(n)^{((n-1)1)}\right)^m$ to the product potential function on $\left(\mathfrak{sl}(n)^{((n-1)1)}\right)^m$ as in \cite[Lemma 4.12]{Manolescu:nilpotent} and \cite[Lemma 30]{SS06}. These copies of $\CC \PP^{n-1}$ in particular embed as the zero sections of the product $(T^*\CC \PP^{n-1})^m$ in the local model of \eqref{eq:local isom}. Lemma \ref{l:equi_local_isom} shows that the action on $(T^*\CC \PP^{n-1})^m$ may be taken to be diagonal, where the action on each term of the product is the $m=1$ action specified in Lemma \ref{l:m=1model}. The result now follows.
\end{proof}

Lemma~\ref{lem:lagrangian-involution} implies that the fixed set $\eL_{\mathfrak m}^{\fix}$ has exactly $2^{m}$ connected components. We now pin down which components of $\cY_{m,n,\tau}^{\fix}$ contain a nonempty connected component of $\eL_{\mathfrak m}^{\fix}$.

\begin{proposition} \label{prop:fixed-components} Let $\fkm$ be a crossingless matching, so that $\mathfrak{m}=(\delta_1,\cdots,\delta_m)$ has endpoints $\delta_i(0)=\lambda_i$ and $\delta_i(1)=\mu_{\rho(i)}$. Let $\eL_{\mathfrak{m}}$ be the Lagrangian associated to $\mathfrak{m}$. The $2^m$ connected components of the fixed Lagrangian $\eL_{\mathfrak{m}}^{\fix}$ lie in the components
\[\cY_{m,k,\tau_{\pi\setminus\pi'}}\times\cY_{m,n-k,\tau_{\pi'}}\]
corresponding to embeddings $(\pi',f)$ with the following properties:

\begin{itemize}
\item The partitions $\pi'$ and $\pi \setminus \pi'$ are of the forms \[ \pi' = ((n-k)^s(n-k-1)^{m-s}1^{m-s}) \qquad \qquad \pi \setminus \pi' = ((k)^{m-s}(k-1)^s 1^s)\] for some $0 \leq s \leq m$.

\item Choose an ordering $\sigma$ of $\tau$ such that $\sigma = (\mu_1, \dots, \mu_m,\lambda_1, \dots, \lambda_m)$; that is, so that $\alpha_{i} = \mu_i$ if $i=1,\dots, m$ and $\alpha_{i} = \lambda_{i-m}$ if $i=m+1,\dots, 2m$. Then we require that 
\begin{enumerate}[label=(\alph*)]
	\item If $1 \leq i \leq m$ and $i+m \in \mathrm{Im}(f)$, then we have that $\rho(i)= f(j)$ for some $j$ with $\pi'_{j} = n-k-1$. 
	\item If $1 \leq i \leq m$ and $i+m \notin \mathrm{Im}(f)$, then we have that $\rho(i)=f(j)$ for some $j$ with $\pi'_{j} = n-k$.
	\end{enumerate}
\end{itemize}
\end{proposition}

\begin{proof} Let $\pi = ((n-1)^m1^m)$ be the standard partition of $mn$ and as usual let $\tau = (\vec{\mu},\vec{\lambda})$ be a configuration of $2m$ distinct values with  $\vec{\mu} =(\mu_1,\dots, \mu_m)$ eigenvalues of multiplicity $n-1$ and $\vec{\lambda} = (\lambda_1, \dots, \lambda_m)$ eigenvalues of multiplicity one. Assume that the values $\mu_j$ and $\lambda_i$ are labelled so that $\rho$ is the identity permutation, that is, so that $\lambda_i$ is matched with $\mu_i$. Moreover, notice that the product structure in the local isomorphism ~\eqref{eq:local isom} encodes the matching: the $i$th component of the product $\left( \frsl(n)^{((n-1)1)} \right)^m$ corresponds to the matched eigenvalues $\mu_i$ with multiplicity $n-1$ and $\lambda_i$ with multiplicity one. 

It therefore suffices to check the statement for $m=1$ and the partition $\pi = (n-1,1)$ with configuration $(\mu_1,\lambda_1)$. It is clear that in this case the two fixed components of $T^*\CC\PP^{n-1}$ under the action contain the two fixed components of the Lagrangian. Moreover, there are exactly two embeddings $(\pi',f)$ with the property that $\pi'$ is a partition of $m(n-k)=n-k$ into $\pi=(n-1,1)$ as follows:

\begin{itemize}
\item We may have $\pi' = (n-k-1,1)$ and $f \co \{1,2\} \rightarrow \{1,2\}$ is the identity. Then $\pi \setminus \pi' = (k)$. In the notation of Lemma \ref{lem:fixedptslice}, this is the component of the fixed set for which the eigenvalues $Y_{11}^b$ are $\mu_1$ with multiplicity $n-k-1$ and $\lambda_1$ with multiplicity one, and the eigenvalue of $Y_{11}^a$ is $\mu_1$ with multiplicity $k$.

\item We have have $\pi'=(n-k)$ and $f \co \{1\} \rightarrow \{1,2\}$ has $f(1)=1$. Then $\pi \setminus \pi' = (k-1,1)$. In the notation of Lemma \ref{lem:fixedptslice}, this is the component of the fixed set for which the eigenvalue $Y_{11}^b$ is $\mu_1$ with multiplicity $n-k$, and the eigenvalues of $Y_{11}^a$ are $\mu_1$ with multiplicity $k-1$ and $\lambda_1$ with multiplicity one.
\end{itemize}

This is exactly the claim in the case that $m=1$ and $\rho$ is the identity permutation. The general result follows.
\end{proof}

\begin{proposition}\label{p:fixedFloer} For $K$ a knot, exactly two components of $\cY_{m,n,\tau}^{\fix}$ contain a nonempty component of both fixed sets $\eL_0^{\fix}$ and $\eL_{\beta}^{\fix}$. These two components are the unique components isomorphic to $\cY^{\pi_{std}}_{m, k,\tau}$ and $\cY^{\pi_{std}}_{m, n-k,\tau}$. Taking Lagrangian Floer cohomology in these components computes the $\frsl(k)$ and $\frsl(n-k)$ homology of $K$ respectively.
\end{proposition}

\begin{proof} Suppose that $(\pi',f)$ is an embedding of a partition $\pi'$ into $\pi = ((n-1)^m1^m)$. As usual, let $\mathfrak{m}_0$ be the crossingless matching determined by the identity permutation. If $\mathfrak m_{\beta}$ is a crossingless matching with associated permutation $\rho$, then $K$ is a knot if and only if $\rho$ is a cyclic permutation; that is, it has no cycle of length less than $m$. We may assume up to relabelling that $\rho$ is in fact the standard cyclic permutation, so that $\rho(i)=i+1$ modulo $m$. Furthermore, as in Proposition~\ref{prop:fixed-components}, we fix an ordering $\sigma$ of $\tau$ so that for $i=1,\dots, m$, we have $\alpha_i = \mu_i$ and $\alpha_{m+i} = \lambda_i$. 

We now consider the components of $\cY_{m,n,\tau}^{\fix}$ which contain both a component of $\eL_0^{\fix}$ and $\eL_{\beta}^{\fix}$. Let $(\pi',f)$ be an embedding of $\pi' = ((n-k)^{s}(n-k-1)^{m-s}1^{m-s})$ into $\pi$. We observe that the conditions of Proposition~\ref{prop:fixed-components} must be satisfied for both the identity permutation and the standard permutation $\rho$. Suppose that $1 \leq i \leq m$ and $i+m \in \mathrm{Im}(f)$. Applying Proposition~\ref{prop:fixed-components} to the identity permutation then implies that there is some $j$ such that $f(j) = i$ with $\pi'_j=n-k-1$. But then applying the same proposition to the standard permutation $\rho$ implies that if $i>1$ then $i-1+m \in \mathrm{Im}(f)$, and if $i=1$ then $2m \in \mathrm{Im}(f)$. Iterating this procedure implies that either every element of $\{m+1,\dots, 2m\} \subset \mathrm{Im}(f)$ or $\{m+1,\dots, 2m\} \cap \mathrm{Im}(f)=\emptyset$. The two such embeddings satisfying this condition are 
\begin{itemize}

\item The embedding $(\pi',f)$ where $\pi' = ((n-k-1)^m 1^m)$ and $f \co \{1, \dots, 2m\} \rightarrow \{1, \dots, 2m\}$ is the identity map, so that $\pi \setminus \pi' = (k)^m$. This embedding corresponds to the connected component $\{\mathrm{pt}\} \times \cY^{\pi_{std}}_{m,n-k,\tau}$ of $\cY_{m,n,\tau}^{\fix}$.

\item The embedding $(\pi',f)$ where $\pi' = (n-k)^m$ and $f \co \{1,\dots,m\} \rightarrow \{1,\dots, 2m\}$ is $f(j)=j$, so that $\pi \setminus \pi' = ((k-1)^m1^m)$. This embedding corresponds to the connected component $\cY^{\pi_{std}}_{m,k,\tau} \times \{\mathrm{pt}\}$ of $\cY_{m,n,\tau}^{\fix}$.

\end{itemize}

We therefore see that exactly the two described components of $\cY_{m,n,\tau}^{\fix}$ contain components of both $\eL_0^{\fix}$ and $\eL_{\beta}^{\fix}$. Since we know from Lemma~\ref{l:equi_local_isom} and Corollary \ref{cor:invariant-parallel} that the involution is compatible with both rescaled parallel transport and the local isomorphism \eqref{eq:local isom} used in constructing the Lagrangians, it follows that the fixed sets in these two components are in fact isomorphic to the Lagrangians associated to $\mathfrak{m}_0$ and $\mathfrak{m}_{\beta}$ in the ranks $k$ and $n-k$. The conclusion of the proposition follows.
\end{proof}

\begin{example} \label{example:n5m2unknot} As in Example~\ref{example:n5m2}, let $n=5$ and $m=2$, and $k=3$, so that $n-k=2$. The six connected components of $\cY^{\fix}_{2,5,\tau}$ are enumerated in Example~\ref{example:n5m2}. Let $\tau =( \{\mu_1,\mu_2\},\{\lambda_1, \lambda_2\})$. We consider the decomposition of the unknot $K$ given in Figure~\ref{fig:matching_example}, so that $\mathfrak{m}_0$ is the identity matching and $\mathfrak{m}_{\beta}$ is the crossingless matching depicted in the figure, for which $\rho$ is the nontrivial permutation. 

Since $m=2$, we see that each of the fixed Lagrangians has four connected components. We now locate them.

\begin{enumerate}

\item The component of $\cY_{2,5,\tau}^{\fix}$ corresponding to $(\pi',f_1)$ where $\pi'=(2,2)$ is diffeormorphic to $\cY^{\pi_{std}}_{2,3,\tau}\times \{\mathrm{pt}\}$ and contains a component of $\eL_{0}^{\fix}$ and of $\eL_{\beta}^{\fix}$, each of which is diffeomorphic to $\CC \PP^2 \times \CC \PP^2$. Lagrangian Floer cohomology in this component of the fixed set computes a copy of the $\frsl(3)$ homology of $K$.

\item The component of $\cY_{2,5,\tau}^{\fix}$ corresponding to $(\pi', f_2)$ where $\pi'=(1,1,1,1)$ is diffeomorphic to $\{\mathrm{pt}\} \times \cY^{\pi_{std}}_{2,2,\tau}$ and contains a component of $\eL_0^{\fix}$ and of $\eL_{\beta}^{\fix}$, each of which is diffeomorphic to $\CC \PP^1 \times \CC \PP^1$. Lagrangian Floer cohomology in this component of the fixed set computes a copy of the $\frsl(2)$ homology of $K$.

\item The component of $\cY_{2,5,\tau}^{\fix}$ corresponding to $(\pi',f_3)$ with $\pi'=(2,1,1)$ contains a component of $\eL_\beta^{\fix}$. This component of $\cY_{2,5,\tau}^{\fix}$ is of the form $\cY_{2,3,\tau_{\pi\setminus\pi'}}\times\cY_{2,2,\tau_{\pi'}}$ and contain a component of $\eL_\beta^{\fix}$ which appears as a product of a copy of $\CC \PP^{2}$ in the first factor with a copy of $\CC \PP^1$ in the second factor.

\item The component of $\cY_{2,5,\tau}^{\fix}$ corresponding to $(\pi',f_4)$ with $\pi'=(2,1,1)$ contains a component of $\eL_0^{\fix}$, which similarly appears in $\cY_{2,3,\tau_{\pi\setminus\pi'}}\times\cY_{2,2,\tau_{\pi'}}$ as a product of a copy of $\CC \PP^{2}$ in the first factor with a copy of $\CC \PP^1$ in the second factor.

\item The component of $\cY_{2,5,\tau}^{\fix}$ corresponding to $(\pi',f_5)$ with $\pi'=(2,1,1)$ contains a component of $\eL_0^{\fix}$, which similarly appears in $\cY_{2,3,\tau_{\pi\setminus\pi'}}\times\cY_{2,2,\tau_{\pi'}}$ as a product of a copy of $\CC \PP^{2}$ in the first factor with a copy of $\CC \PP^1$ in the second factor.

\item The component of $\cY_{2,5,\tau}^{\fix}$ corresponding to $(\pi',f_6)$ with $\pi'=(2,1,1)$ contains a component of $\eL_\beta^{\fix}$, which likewise appears in $\cY_{2,3,\tau_{\pi\setminus\pi'}}\times\cY_{2,2,\tau_{\pi'}}$ as a product of a copy of $\CC \PP^{2}$ in the first factor with a copy of $\CC \PP^1$ in the second factor.

\end{enumerate}
\end{example}

\subsection{Real structures} \label{sec:real}
Let $\tau=(\vec{\mu},\vec{\lambda})\in \frh^{\pi,reg}/W^{\pi}$. In this section we prove Theorem~\ref{thm:stable_real}, that given a braid element $\beta$, the triple $(\cY_{m,n,\tau},\eL_{0},\eL_{\beta})$ admits a stable real polarization. Moreover we will show that for the involution $\iota_k$, the fixed set $\cY_{m,n,\tau}^{\fix}$ admits a set of stable normal real polarizations, allowing us to use Corollary~\ref{cor:components} to deduce Theorem~\ref{thm:mainss}.

We begin by giving local descriptions of the tangent and normal bundles of $\cY_{m,n,\tau}$ and the tangent and normal bundles of $\cY_{m,n,\tau}^{\fix}$. We will need the following preliminary lemma.

\begin{lemma}\label{l:TS} For any $Y \in \cS_{m,n}$ and $Z \in \frgl(mn)$, the vector $[Y,Z] \in T_Y \cS_{m,n}$ if and only if $Z$ takes the form
\[
Z=(Y^{m-1}X|\cdots|YX|X)
\]
for some $X\in M_{(mn)\times n}$.
\end{lemma}

\begin{proof}
First note that, $T_Y\cS_{m,n}$ is the set of $(mn)\times (mn)$-matrices $R=(R_{ij})$, where $R_{ij}$ is a $n \times n$ matrix for all pairs $1 \leq i,j \leq m$ such that $R_{ij}=0$ for $j>1$.

 Assume $Z=(Z_1\mid Z_2\mid\cdots\mid Z_m)$ where each $Z_i\in M_{(mn)\times n}$. Then, 
		\[YZ=(YZ_1\mid YZ_2\mid\cdots\mid YZ_m)\quad\text{and}\quad ZY=(*\mid Z_1\mid Z_2\mid \cdots\mid Z_{m-1}).\]
Therefore  $[Y,Z]\in T_Y\cS_{m,n}$ if and only if $Z_i = YZ_{i+1}$ for $1 \leq i \leq m-1$, or equivalently if $Z=(Y^{m-1}X|\cdots|YX|X)$ for some $X\in M_{(mn)\times n}.$
\end{proof}

This leads us to the following set of descriptions. In the below, recall that the fixed set of $\cY_{m,n,\tau}$ under $\iota_k$ has multiple components of varying dimension, and therefore throughout this section $T\cY_{m,n,\tau}^{\fix}$ and $N\cY_{m,n,\tau}^{\fix}$ vary in dimension between the connected components of the fixed set.

\begin{corollary} \label{cor:bundle-descriptions} For any $Y\in \cY_{m,n,\tau}$ we have 
\[T_Y\cY_{m,n,\tau}=\{[Y,Z_X]\ |\ Z_X:=(Y^{m-1}X|\cdots|YX|X),\ X\in M_{(mn)\times n}\}.\]
Moreover, if $\cY_{m,n,\tau}^{\fix}$ denotes the fixed set of $\cY_{m,n,\tau}$ under the involution $\iota_k$, for any $Y\in \cY_{m,n,\tau}^{\fix}$ we have 
\[\begin{split}
	&T_Y\cY^{\fix}_{m,n,\tau}=\{[Y,Z_X]\ |\ Z_X:=(Y^{m-1}X|\cdots|YX|X),\ X\in M_{(mn)\times n}\ \text{and}\ \ T_kXt_k=X\}\\
	&N_Y\cY^{\fix}_{m,n,\tau}=\{[Y,Z_X]\ |\ Z_X:=(Y^{m-1}X|\cdots|YX|X),\ X\in M_{(mn)\times n}\ \text{and}\ \ T_kXt_k=-X\}
\end{split}\]
	\end{corollary}
	
	\begin{proof}
		By definition $\cY_{m,n,\tau}=(\chi^{\pi})^{-1}(\tau)\cap \cS_{m,n}$. Note that $(\chi^{\pi})^{-1}(\tau)$ is the orbit $\cO(Y)$ of $Y$ under the adjoint action, and therefore \[T_{Y}((\chi^{\pi})^{-1}(\tau))=ad(\kg)Y=\{[Y,Z]\mid Z\in\kg\}.\]
		Thus, $T_Y\cY_{m,n,\tau}$ consists of all $[Y,Z]\in T_{Y}\cS_{m,n}$ so $T_Y\cY_{m,n,\tau}$ is the claimed subspace by Lemma \ref{l:TS}.
 Note that this part of the argument does not require $Y$ to be in the fixed set of $\cY_{m,n,\tau}$; this description of $T_Y\cY_{m,n,\tau}$ holds generally.
		
	We next observe that $T_Y\cY^{\fix}_{m,n,\tau}$ and $N_Y\cY^{\fix}_{m,n,\tau}$ are the $1$ and $-1$ eigenspaces of $(\iota_k)_*$ on $T_Y\cY_{m,n,\tau}$. Note further that 
	\[(\iota_k)_*[Y,Z_X]=[Y,T_kZ_XT_k]=[Y,Z_{T_kXt_k}],\]
	and so $T_Y\cY^{\fix}_{m,n,\tau}$ and $N_Y\cY^{\fix}_{m,n,\tau}$ are the claimed subspaces.
		 \end{proof}

Consider an $mn^2$-dimensional trivial complex vector bundle $V$ over $\cY_{m,n,\tau}$ such that $V_Y\cong M_{mn\times n}$ and let $f \co V\to T\cY_{m,n,\tau}$ be the bundle morphism given by $f(X)=[Y,Z_X]$ at $Y \in \cY_{m,n,\tau}$. 
By Corollary \ref{cor:bundle-descriptions}, $f$ is surjective.
Define $K$ to be the kernel of $f$, which is to say that $K_Y$ consists of matrices $X$ such that $[Y,Z_X]=0$.

 
 Our strategy is to define a real structure on $K$ and use it to produce a stable real structure on $T\cY_{m,n,\tau}$. Before we discuss the general case, we consider the simplest example, as follows.

\begin{example}[$K_{\mathbb{R}}$ when $m=1$]\label{example:realstructurem1}
Let $m=1$, so that $\cS_{m,n} \cap \frg^{\pi}=\frgl(n)^{((n-1)1)}$. For simplicity we denote $\cY_{\tau}:=\cY_{1,n,\tau}=(\chi^{((n-1)1)}_\frgl)^{-1}(\tau)$ where $\tau=(\{\mu\},\{\lambda\})$. Suppose $ \mu \neq \lambda$. Then, for any $Y\in \cY_{\tau}$ we have $V_{Y}\cong M_{n\times n}$ and furthermore
	\[T_{Y}\cY_{\tau}=\{[Y,X]\mid X\in M_{n\times n}\}\quad\text{and}\quad K_Y\cong \{X\in M_{n\times n}\mid [Y,X]=0\}.\]
To define the fiber $K_{Y,\RR}$ of $K_{\RR}$ over $Y$, we pick orthonormal bases $\{v_1,\cdots,v_{n-1}\}$ and $\{w\}$ for the eigenspaces $E_Y(\mu)$ and $E_Y(\lambda)$, respectively, and let $B:=(v_{1}\mid \cdots\mid v_{n-1}\mid w)$. 
In the basis $B$, we have
\[
B^{-1}YB=diag(\mu, \dots,\mu, \lambda).
\]
Write any $X\in V$ in the basis $B$ in the block form \[X':=B^{-1}XB=\begin{pmatrix}X'_{11}&X_{12}'\\ X_{21}'&X_{22}'\end{pmatrix}\] where $X_{11}'\in M_{(n-1)\times (n-1)}$ and $X_{22}'\in \CC$. Then
$B^{-1}[Y,X]B=\begin{pmatrix}0&(\mu-\lambda)X_{12}'\\ (\lambda-\mu) X_{21}'&0\end{pmatrix}$ so 
\[ K_Y\cong \{X\in M_{n\times n}\mid X_{12}'=0  \text{ and }X_{21}'=0\}.\]

We define 
\[
K_{Y,\RR}:=\{X \in K_Y \mid X'_{11} \text{ is Hermitian},\ X_{22}'\in\RR \}.\]
It is easy to check that the definition of $K_{Y, \mathbb{R}}$ does not depend on the choice of orthonormal bases of $E_Y(\mu)$ and $E_Y(\lambda)$, and furthermore that these fibres vary smoothly with $Y$ to form a real vector bundle $K_{\RR}$ such that $K_{\RR} \otimes_{\RR} \mathbb{C} \simeq K$. Since $V$ is trivial consider the trivial real subbundle $V_{\RR}\subset V$ such that $V_{\RR}\otimes \CC=V$ and let $T_{\RR}\cY_{\tau}:=[V_{\mathbb{R}}]-[K_{\mathbb{R}}]=-[K_{\mathbb{R}}]$. 

\end{example}

Now we consider the general case.

\begin{proposition}\label{p:K_real}
The complex vector bundle $K$ admits a real structure. 
	\end{proposition}
	
	\begin{proof} Let $E_Y(\mu_i)$ and $E_Y(\lambda_i)$ be the eigenspaces of $Y$. Recall that by Lemma~\ref{lem:injectivity}, projection on the first $n$-coordinates gives an injective map $p$ from  $E_Y(\mu_i)$ to $\CC^n$, and likewise for $E_Y(\lambda_i)$.
For any $i=1,\cdots,m$, consider a basis	$\{v_{i1},\cdots,v_{i(n-1)}\}$ for $E_{Y}(\mu_i)$ such that $\{p(v_{i1}),\cdots,p(v_{i(n-1)})\}$ is a orthonormal basis for $p(E_Y(\mu_i))$ with the standard unitary structure on $\CC^n$. Putting these vectors together, let 
	\[B=(v_{11}|\cdots|v_{1(n-1)}|v_{21}|\cdots|v_{m(n-1)}|w_1|\cdots|w_m)\]
	where for each $i$, the vector $w_i\in E_Y(\lambda_i) \setminus \{0\}$ is a unit-length eigenvector for $\lambda_i$.

We then define the real subspace $K_{Y,\RR}$ to be 
\[K_{Y,\RR}=\{X\in M_{(mn)\times n}\ |\ [Y,Z_X]=0, B^{-1}Z_XB=\overline{B^{-1}Z_XB}^t\}\]
Recall that since $[Y,Z_X]=0$, the matrix $Z_X$ acts on the eigenspaces of $Y$. Therefore the matrix $B^{-1}Z_XB$ which presents $Z_X$ in the basis $B$ is of the diagonal block form $diag(M_1,\cdots,M_{m},c_1,\cdots, c_m)$ where $M_i$ and $c_i$ are the actions of $Z_X$ on $E_Y(\mu_i)$ and $E_Y(\lambda_i)$ in the bases $\{v_{i1},\cdots,v_{i(n-1)}\}$ and $\{w_i\}$, respectively. We therefore see that $B^{-1}Z_XB$ is Hermitian if and only if $M_i$ is Hermitian and $c_i\in\RR$ for all $i$.


We now check that $K_{Y,\RR}$ is independent of the choices of bases for the eigenspaces. If we choose another basis $\{v'_{i1},\cdots,v'_{i(n-1)}\}$ for $E_{Y}(\mu_i)$ such that $\{p(v'_{i1}),\cdots,p(v'_{i(n-1)})\}$ is also an orthonormal basis for $p(E_{Y}(\mu_i))$, then for each $i$ there is an $(n-1) \times (n-1)$ unitary matrix $P$ such that
\[B_i'=B_iP,\ \ \ \text{where}\ \ \ \ B_i=(v_{i1}|\cdots|v_{i(n-1)})\ \text{and}\ B_i'=(v'_{i1}|\cdots|v'_{i(n-1)})\]
Moreover,the  action of $Z_X$ in the new basis is given by a matrix $M_i'$ where $Z_XB_i'=B_i'M_i'$. Thus, $M_i=PM_i'P^{-1}=PM_i'\overline{P}^{t}$ and so $M_i=\overline{M_i}^t$ implies $\overline{M_i'}^t=M_i'$. It follows that the definition of $K_{Y,\RR}$ does not depend on the choice of basis. 

Furthermore, we observe that $Z_{iX}=iZ_X$ so $\overline{Z}_{iX}^{t}=-i\overline{Z}_X^t$. With this in mind, following the chain of equalities
 \[\overline{B}^t\overline{Z_{iX}}^t(\overline{B}^{-1})^t=-i\overline{B}^t\overline{Z_X}^t(\overline{B}^{-1})^t=-B^{-1}Z_{iX}B\]
shows us that $B^{-1}Z_XB$ is Hermitian if and only if $B^{-1}Z_{iX}B$ is skew-Hermitian.
Consequently, $iK_{Y,\RR}$ is 
\[iK_{Y,\RR}=\{X\in M_{(mn)\times n}\ |\ [Y,Z_X]=0, B^{-1}Z_XB=-\overline{B^{-1}Z_XB}^t\}.\]
and in particular $iK_{Y,\RR}$ intersects $K_{Y,\RR}$ only at the origin.

Finally, we need to check that $K_{Y,\RR}\oplus iK_{Y,\RR}$ is the entirety of $K_Y$.

We first claim that it suffices to check $K_{Y,\RR}\oplus iK_{Y,\RR}=K_Y$ for one $Y\in \cY_{m,n,\tau}$. For if $Y'\in\cY_{m,n,\tau}$ is any other element, there exists $g\in G$ such that $gYg^{-1}=Y'$.
In fact, any $g \in G$ such that $g(E_{Y}(\mu_i))=E_{Y'}(\mu_i)$ and $g(E_{Y}(\lambda_i))=E_{Y'}(\lambda_i)$ for all $i$  satisfies $gYg^{-1}=Y'$ because $Y$ and $Y'$ are semi-simple.
In particular, $g$ can be taken to be such that $p\circ g|_{E_{Y}(\mu_i)}=q_i \circ p$ where $q_i \co p(E_{Y}(\mu_i))\to p(E_{Y'}(\mu_i))$ is unitary in the orthonormal bases of $p(E_{Y}(\mu_i))$ and $p(E_{Y'}(\mu_i))$ for all $i$. For instance, one may choose $g$ such that $q_i$ is equal to the identity on $p(E_{Y}(\mu_i))\cap p(E_{Y'}(\mu_i))$ while mapping the orthogonal complement of $p(E_{Y}(\mu_i))\cap p(E_{Y'}(\mu_i))$ in $p(E_{Y}(\mu_i))$ to its orthogonal complement in $p(E_{Y'}(\mu_i))$ preserving the inner product. Fix such $g$ and let $Q_i$ be the matrix presentation of $q_i$.
Then, if $X\in K_{Y}$ it follows that 
$Z'=gZ_Xg^{-1}$
commutes with $Y'$.
 In particular this means that  $[Y',Z']=0 \in T\cS_{m,n}$
so by Lemma \ref{l:TS}, $Z'$ is necessarily of the form
\[
Z'=Z_{X'}=((Y')^{m-1}X' |\cdots| Y'X| X')
\]
for a unique choice of $X' \in M_{(mn) \times n}$.
Moreover, the actions of $Z_X$ on $E_{Y}(\mu_i)$ and $Z'$ on $E_{Y'}(\mu_i)$ are related by $M_i'=Q_iM_iQ_i^{-1}$. But since $Q_i$ is unitary, $X\in K_{Y,\RR}$ if and only if $X'\in K_{Y',\RR}$.  Therefore as claimed it suffices to check $K_{Y,\RR}\oplus iK_{Y,\RR}=K_Y$ for a single $Y\in \cY_{m,n,\tau}$.

With this in mind, we may make a particularly simple choice of element of $\cY_{m,n,\tau}$, as follows. Let $z^{m}-a_1z^{m-1}-\cdots-a_{m}:=\prod_{i=1}^m(z-\lambda_i)$ and  $z^{m}-b_1z^{m-1}-\cdots-b_{m}:=\prod_{i=1}^m(z-\mu_i)$. Let $D\in\mathcal{S}_{m,n}$ be the matrix which has the square matrices $D_i:=D_{i1}$ of the form $D_i=diag(b_i,\cdots,b_i,a_i)$ in its first column. It is straightforward to check that $D$ lies in in $(\chi^{\pi})^{-1}(\tau)$.
 We are going to show that  $K_{D,\RR}\oplus iK_{D,\RR}=K_D$.
 
 Toward this end, we observe that permuting rows and columns of $D$ we get the block diagonal matrix $D'=\begin{pmatrix}D_b&0\\0&D_a\end{pmatrix}$ where  

\[D_b=\begin{pmatrix}b_1I_{n-1}&I_{n-1}&0&\cdots&0\\ b_2I_{n-1}&0&I_{n-1}&\cdots&0\\\cdots&&&&\\ b_{m}I_{n-1}&0&0&\cdots&0\end{pmatrix}\ \ \text{and}\ \ \ D_a=\begin{pmatrix}a_1&1&0&\cdots&0\\ a_2&0&1&\cdots&0\\\cdots&&&&\\ a_{m}&0&0&\cdots&0\end{pmatrix}.\]
The block matrix $Z'=\begin{pmatrix}Z'_{11}&Z'_{12}\\ Z'_{21}&Z'_{22}\end{pmatrix}$, obtained from $Z_X$ by applying the same permutations, commutes with $D'$ if and only if 
\[[Z'_{11},D_b]=0,\ \ [Z'_{22},D_a]=0,\ \ D_bZ'_{12}=Z'_{12}D_a,\ \ \text{and}\ \ D_aZ'_{21}=Z'_{21}D_b.\]
The matrices $D_a$ and $D_b$ have no common eigenvalue, because the eigenvalues of $D_a$ and $D_b$ are $\{\lambda_1,\cdots,\lambda_m\}$ and $\{\mu_1,\cdots,\mu_m\}$, respectively, so the Sylvester equations $D_bW=WD_a$ and $D_aW=WD_b$ each have the unique solution $W=0$. In other words, $Z'$ commutes with $D$ if and only if  $Z'_{12}=0$ and $Z_{21}'=0$. Reversing the permutations of rows and columns used to turn $Z_X$ into $Z'$, we conclude that the matrix $X$ used to generate $Z_X=(D^{m-1}X|\cdots|DX|X)$ has a block form 
\[X=\begin{pmatrix}X_1\\ X_2\\...\\ X_m\end{pmatrix}\ \ \text{where}\ \ X_j=\begin{pmatrix}X_{jb}&0\\0&X_{ja} \end{pmatrix}.\]
Note that the matrices $X_j$ are $n\times n$ matrices and the matrices $X_{jb}$ are $(n-1)\times (n-1)$ matrices. Moreover, $Z'_{11}=(D_b^{m-1}X_b|\cdots|D_bX_b|X_b)$ and $Z'_{22}=(D_a^{m-1}X_a|\cdots|D_aX_a|X_a)$ for 
\[X_b=\begin{pmatrix}X_{1b}\\ \cdots\\ X_{mb}\end{pmatrix}\ \ \text{and}\ \ X_a=\begin{pmatrix}X_{1a}\\ \cdots\\ X_{ma}\end{pmatrix}.\]

Now, for any $1\le i\le m$, by a direct computation we have

\[E_{D}(\mu_i)=\left\{\begin{pmatrix} q_1(\mu_i)v\\ q_2(\mu_i)v\\ q_3(\mu_i)v\\ \cdots\\ q_{m}(\mu_i)v \end{pmatrix}\ \middle|\ \begin{aligned} v=(z_1,z_2,\cdots,z_{n-1},0)^t\\q_1(z)=1,\  q_{j+1}(z)=z^{j}-\sum_{\ell=0}^{j-1}z^{\ell}b_{j-\ell}\end{aligned}\right\}\]

A similar result holds for $E_D(\lambda_i)$ via replacing $\mu_i$ with $\lambda_i$, replacing $b_j$ with $a_j$ for all $j$, and setting $v=(0,0,\cdots,0,z_n)^t$. Similarly, replacing $v$ with  $v_b=(z_1,z_2,\cdots,z_{n-1})$ in $E_{D}(\mu_i)$ we get $E_{D_b}(\mu_i)$, and replacing $v$ with $v_a=z_n$ in $E_{D}(\lambda_i)$ recovers $E_{D_a}(\lambda_i)$.

Let $v_{ij}\in E_{D}(\mu_i)$ be the vector obtained by setting $v=e_j$ in the formula above, and likewise $v_{ijb}\in E_{D_b}(\mu_i)$ be the vector obtained by setting $v_b=e_j$ in the analogous fomula. In these bases, the action of $Z_X$ on $E_{D}(\mu_i)$ and the action of $Z_{11}'$ on $E_{D_b}(\mu_i)$ are given by the same matrix.  Denote this matrix by $M_i$, and let $B_{ib}=(v_{i1b}|v_{i2b}|\cdots|v_{i(n-1)b})$. We then see that
$Z'_{11}B_{ib}=B_{ib}M_i$ and so 
\begin{equation} \label{eq:first-M}
(I_{n-1}|0|\cdots|0)Z'_{11}B_{ib}=M_i.\end{equation}

We would like to write $M_i$ explicitly in terms of $X_b$ and $\mu_i$. Toward this end, we observe that in block diagonal form we have the following expressions for $B_{ib}$ and $Z_{11}'B_{ib}$:
\[B_{ib}=\begin{pmatrix}q_1(\mu_i)I_{n-1}\\ q_2(\mu_i)I_{n-1}\\ q_3(\mu_i)I_{n-1}\\ \cdots\\ q_{m}(\mu_i)I_{n-1}\end{pmatrix}\quad\text{and}\quad Z_{11}'B_{ib}=\left(\sum_{j=1}^{m}q_j(\mu_i)D_b^{m-j}\right)X_b.\]
Let $R_j=\left(0\mid 0\mid \cdots\mid I_{n-1}\mid\cdots\mid 0\right)$ where $I_{n-1}$ is the $j$-th block of size $(n-1)\times (n-1)$. Then we see that $R_j D_b=b_jR_1+R_{j+1}$.  Since $D_b$ is a square matrix, the expression $q_{j}(D_b)$ computed as a formula for square matrices makes sense. Moreover by a direct computation, we have 
$R_{j}=R_1q_{j}(D_b)$, which implies 
\begin{equation}\label{eq:rearrange-qj}\sum_{j=1}^{m}\mu_i^{m-j}R_j=R_1\sum_{j=1}^m\mu_i^{m-j}q_{j}(D_b)=R_1\left(\sum_{j=1}^{m}q_j(\mu_i)D_b^{m-j}\right)\end{equation}
where the second equality follows from expanding the middle summation and checking that the coefficient of $D_b^{m-j}$ is equal to $q_j(\mu_i)$.
Applying the equality \eqref{eq:rearrange-qj} to \eqref{eq:first-M} and the block diagonal form of $Z'_{11}B_{ib}$ above now gives the following expression for $M_i$ in terms of $X_b$ and $\mu_i$:
\begin{equation}\label{eqn:gemini-helped}
M_i=R_1Z'_{11}B_{ib}=R_1\left(\sum_{j=1}^{m}q_j(\mu_i)D_b^{m-j}\right)X_b=\sum_{j=1}^{m}\mu_i^{m-j}R_jX_b=\sum_{j=1}^{m}\mu_i^{m-j}X_{jb}.\end{equation}

Now, recall that for $X\in K_{D,\RR}$,  $M_i$ is a Hermitian matrix for all $i$. The set of matrix equations $M_i=\overline{M}_i^T$ for $i=1,\dots,m$ may be broken down into the following linear systems of equations. For each pair $1\le k<\ell\le n-1$ by considering the $k\ell$ entry of the equation $M_i=\overline{M}_i^T$ for each $i$, we obtain $m$ equations in $2m$ variables which correspond to the $k\ell$ and $\ell k$ entries of the matrices $X_{jb}$ for $j=1,\cdots,m$. It is easy to check that the matrix of coefficients has rank $m$, implying that the kernel has complex dimension $m$. Similarly for each value of $k$ such that $1 \le k \le n-1$, by considering the imaginary part of the $k$th diagonal entry of each equation $M_i=\overline{M}_i^T$, we get a system of $m$ equations with real coefficients in $2m$ variables associated to the $k$th diagonal entries of the matrices $X_{jb}$. So the kernels of these systems have real dimension at least $m$. Consequently the total real dimension of the set of $X_b$ with the property that $M_i=\overline{M}_i^t$ for all $i$ is at least $2m\left(\frac{(n-2)(n-1)}{2}\right)+m(n-1)=m(n-1)^2$.

Similarly, given an eigenvector $w$ in $E_{D}(\lambda_i)$, the same line of argument as above shows that we have $Z_Xw=c_iw$ where 
\[c_i=\sum_{j=1}^m \lambda_i^{m-j}X_{ja}.\]
So, given that $c_i\in\RR$ for all $i$, we obtain a system of $m$ equations with real coefficients in $2m$ variables. The kernel of this system is of real dimension at least $m$. 

Summing the two computations above, we conclude that the dimension of $K_{D,\RR}$ is at least $m(n-1)^2+m$. Note that for any $Y \in \cY_{m,n,\tau}$, the bundle $K_Y$ has complex dimension $m(n-1)^2+m$, because $V$ has complex dimension $mn^2$ and $\mathcal Y_{m,n,\tau}$ has complex dimension $2m(n-1)$. It follows that $K_{D,\RR}\oplus iK_{D,\RR} = K_{D}$. Above we argued that checking this equality at a single $Y=D$ sufficed to check it over all $Y$, so we conclude that $K_{\RR}\oplus iK_{\RR} = K$. Thus we have constructed a real structure on $K$ as promised. \end{proof}


\begin{corollary} \label{cor:tangent-real}
Consider the real structure $V_{\RR}\subset V$ such that $V_{\RR}$ is trivial. The tangent bundle $T\cY_{m,n,\tau}$ stably admits a real structure via setting the virtual bundle $T_{\RR}\cY_{m,n,\tau}:=[V_{\RR}]-[K_{\RR}]$. 
	\end{corollary} 
\begin{proof} Recall that we have an exact sequence of complex vector bundles 
\[ 0\rightarrow K \rightarrow V \xrightarrow{f} T\cY_{m,n,\tau}\rightarrow 0.\]
Since short exact sequences of complex vector bundles over a paracompact base split, we have a noncanonical isomorphism of complex bundles $V\simeq T\cY_{m,n,\tau}\oplus K$. It is not necessarily the case that any given such isomorphism is compatible with the real structures on $V$ and $K$. However, we may let $T_{\RR}\cY_{m,n,\tau}:=[V_{\RR}]-[K_{\RR}] \in \mathit{KO}(\cY_{m,n,\tau})$. Then under the map  $\mathit{KO} \rightarrow K$, we see that $[T_{\RR}\cY_{m,n,\tau}]$ is sent to $[T\cY_{m,n,\tau}]$ as promised. Hence after stabilization we have achieved a real structure on the tangent bundle.
 \end{proof}

We are furthermore interested in the behavior of the tangent and normal bundles to the fixed set $\cY_{m,n,\tau}^{\fix}$ under $\iota_k$.

\begin{corollary} \label{cor:fixed-real}
	For any $Y\in\cY^{\fix}_{m,n,\tau}$ let 
	\[K_{Y}^{T}=\{X\in K_Y\mid T_kXt_k=X\}\ \ \text{and}\ \ K_{Y}^{N}=\{X\in K_Y \mid\ T_kXt_k=-X\}.\]
	The vector bundles $K^T$ and $K^N$ on $\cY^{\fix}_{m,n,\tau}$ admit real structures, and therefore the vector bundles $T\cY^{\fix}_{m,n,\tau}$ and $N\cY^{\fix}_{m,n,\tau}$ stably admit real structures.
\end{corollary}
\begin{proof} We have $K \simeq K^T \oplus K^N$. We observe that the real structure $K_{\RR}$ on $K$ constructed in Proposition~\ref{p:K_real} is compatible with this decomposition into the tangent and normal bundles; if $X$ lies in $K^T$ then $iX$ also lies in $K^T$, and likewise for $K^N$. Therefore we can express $K^T$ (resp. $K^N$) as a direct sum of $K_{\RR}^T \oplus iK_{\RR}^T$ (resp. $K_{\RR}^N \oplus iK_{\RR}^N$). The same logic decomposes  $V\simeq V^T\oplus V^N$, and $V^{T}$ (resp. $V^N$) as $V_{\RR}^T \oplus iV_{\RR}^T$ (resp. $V_{\RR}^N \oplus iV_{\RR}^N$). From there the result follows by applying the argument of Corollary \ref{cor:tangent-real} to the tangent and normal bundles of $\cY_{m,n,\tau}^{\fix}$. 
 \end{proof}

We now turn to upgrading the stable real structures of Corollaries~\ref{cor:tangent-real} and \ref{cor:fixed-real} to stable real polarizations by showing compatibility of the stable real structures with the tangent bundles to Lagrangians associated to matchings. As in the previous discussion, we first do the simplest case, which is $m=1$ with the standard Lagrangian, as a warm-up. In particular, it suffices to consider traceless matrices.

 \begin{example}[Compatibility between $T_{\mathbb{R}}\mathcal Y_z$ and the standard Lagrangian] \label{eg:Lag}
Assume $m=1$, $\tau_z=(\{z\},\{-(n-1)z\})$, $\cY_z:=\cY_{\tau_z}$ and let   $D_z=\mathit{diag}(z,z,\cdots,z,-(n-1)z)\in \cY_z$. Then, for $z\neq 0$, let 	\[\eL_z=\{UD_zU^{-1}\mid U\in U(n)\}\]
	be the standard Lagrangian in $\cY_z$ with the symplectic form constructed in Section \ref{sec:ambient}. (We denoted this Lagrangian $\eL_{\alpha}$ in Section~\ref{sec:construction}, but in the arguments that follow it will be helpful to index by the coordinate $z$.) For any $Y=UD_zU^{-1} \in \eL_z$, the first $n-1$ columns of $U$ give a orthonormal basis for $E_Y(z)$ and the last columns give a orthonormal basis for $E_{Y}(-(n-1)z)$. By Example~\ref{example:realstructurem1}, the fibre $K_{Y,\RR}$ is the subspace of all Hermitian matrices that preserve $E_{Y}(z)$, which is to say matrices of the form $U\begin{pmatrix}X_{11}'&0\\ 0&X_{22}'\end{pmatrix}U^{-1}$ with $X_{11}'=\overline{X_{11}'}^t$ and $X_{22}'\in\RR$. Consequently,
	\[[T_{Y,\RR}\cY_{z}]=-[K_{Y,\RR}]=\left[\left\{U\begin{pmatrix}0&X'_{12}\\ X'_{21}&0\end{pmatrix}U^{-1}\mid X'_{21}=\overline{X'_{12}}^t\right\}\right].\]
%
On the other hand, $T_Y\eL_z$ consists of all matrices $[Y,B]$ where $B$ is an element of $\mathfrak{u}(n)$, so that in particular $B+\overline{B}^t=0$. Note that 
\[[UD_zU^{-1},B]=U[D_z,B']U^{-1}=U\begin{pmatrix}0&nzB'_{12}\\ -nzB'_{21}&0\end{pmatrix} U^{-1}\]
where $B'=U^{-1}BU$ is also skew-Hermitian. Therefore, $f$ gives isomorphism between the subbundle 
\[i\left\{U\begin{pmatrix}0&X'_{12}\\ X'_{21}&0\end{pmatrix}U^{-1}\mid X'_{21}=\overline{X'_{12}}^t\right\}\subset V_{Y}\]
and $T_{Y}\eL_z$. So, $[T_Y\eL_z]=[T_{Y,\RR}\cY_z]$.

\end{example}

Continuing with the case $m=1$, we now examine how the real structure behaves when two eigenvalues in the configuration coincide.

\begin{example}[Extending $K_{\mathbb{R}}|_{\eL}$ over the partial Grothendieck resolution]\label{eg:resolution-real-bundle}
Continuing with Example \ref{eg:Lag}, let $\cZ = \frsl(n)^{((n-1)1)}$ and let $\widetilde{\cZ}=\widetilde{\mathfrak{sl}}(n)^{((n-1)1)}$ be its restricted partial Grothendieck resolution. Each element of $\widetilde{\cZ}$ is of the form $(Y,F)$ where $Y\in \cZ$ and $F$ is a flag with the following properties. First, $F$ is of the form $0=F_0\subset F_1\subset F_2=\CC^n$ and preserved by $Y$. Moreover, $F_1$ has complex dimension $\dim_{\CC}F_1=n-1$ and furthermore if we let $z = \chi^{((n-1)1)}(Y)$ be the eigenvalue of $Y$ with multiplicity at least $n-1$, we have $F_1\subseteq E_{Y}(z)$. Let $r\co \widetilde{\cZ}\to \cZ$ denote the resolution map, which is defined by forgetting the flag. 
	
	By \cite[Proposition 2.2]{Manolescu:sln} the map $\widetilde{\chi}^{((n-1)1)}\co\widetilde{\cZ}\to \CC$ which is the composition of $r$ and $\chi^{((n-1)1)}$ is a smooth fibre bundle. Furthermore,  the preimage of $\eL=\cup_{z} \eL_z$ in $\widetilde{\cZ}$, denoted by $\widetilde{\eL}$, is a smooth $\mathbb{CP}^{n-1}$-bundle over $\CC$. To see this, we observe that for $z=0$, 
	\[r^{-1}(\eL_0)=r^{-1}\{0_{n\times n}\}=\{(0,0\subset F_1\subset \mathbb{C}^{n})\mid \dim_{\mathbb{C}}F_1=n-1\}\] is the set of $(n-1)$-dimensional subspaces of $\CC^n$, which is diffeomorphic to $\mathbb{CP}^{n-1}$.  We define a smooth extension of the vector bundle $r^*K_{\RR}|_{\eL\setminus \eL_0}$ over $\widetilde{\eL}$ as follows.  Let $\widetilde{V}$ be the extension of $r^*V|_{\eL\setminus \eL_0}$ to a trivial bundle over $\widetilde{\eL}$. We then let $\widetilde{K}_{(Y,F)}\subset \widetilde{V}_{(Y,F)}$ be the subspace of matrices that are in $K_Y$, which is to say $[Y,X]=0$, and furthermore preserve the flag $F$. Note that this is the kernel of the map 
\[f\co\widetilde{V}_{(Y,F)}\to T_{(Y,F)}\widetilde{\cZ}_z  \qquad \qquad f(X) = ([Y,X],X).\] Here the first entry $[Y,X]$ lies in $\frg$ and the second entry $X$ lies in $\frg/\mathfrak{p}(F)$, where $\mathfrak{p}(F)$ is the set of matrices that preserve the flag.
	
Having defined $\widetilde{K}_{(Y,F)}$, we make an important caveat. Note	that if $Y\neq 0$ then $\widetilde{K}_{(Y,F)}\cong r^*K_{Y}$ and so $\dim \widetilde{K}_{(Y,F)}=(n-1)^2+1$. However, $\dim\widetilde{K}_{(0,F)}=n^2-(n-1)$, since it is the set of upper triangular matrices with respect to the flag. Therefore $\widetilde{K}$ constructed by taking the union of $\widetilde{K}_{(Y,F)}$ over $(Y,F) \in \widetilde{\eL}$ is not a vector bundle. 
	
However, we will see that the situation of extending the real structure over $\widetilde{\eL}$ is better. Let $\widetilde{V}_{\RR}\subset \widetilde{V}$ be trivial and define $\widetilde{K}_{(Y,F),\RR}$ to be the subspaces of matrices in $\widetilde{K}_{(Y,F)}$ that are Hermitian. Then, for any $(Y,F)$, the space $\widetilde{K}_{(Y,F),\RR}$ can be described as the subspace of all Hermitian metrices of the form \[U\begin{pmatrix}X_{11}'&0\\ 0&X_{22}'\end{pmatrix}U^{-1}\] where $U$ is a unitary matrix whose first $n-1$ columns gives a orthonormal basis for $F_1$. We see that unlike $\widetilde{K}$, the real bundle $\widetilde{K}_{\RR}$ forms a smooth vector bundle over $\widetilde{\eL}$. Consequently, we define an extension of the virtual bundle $r^*T_{\RR}\cZ|_{\eL\setminus\eL_0}$ over $\widetilde{\eL}$ as $-[\widetilde{K}_{\RR}]$. For every $z$ this is stably isomorphic to the subbundle
\[\left\{U\begin{pmatrix}0&X'_{12}\\ X'_{21}&0\end{pmatrix}U^{-1}\mid X'_{21}=\overline{X'_{12}}^t\right\}\subset \widetilde{V}_{(Y,F)}.\]
Moreover, $f$ composed with $i$ multiplication gives an isomorphism between this subbundle and $T_{(Y,F)}\eL_{z}$ for all $z$.
\end{example}

We now deal with the general case of colliding pairs of eigenvalues. Let $\tau=(\vec{\mu}, \vec{\lambda})$ be a point in $\mathfrak{h}^{\pi}/{W^{\pi}}$ such that $\mu_i=\lambda_i$ for all $1\le i\le m$  and $\lambda_i \neq \lambda_j$ if $i \neq j$. Recall that the smallest orbit in $(\chi^{\pi})^{-1}(\tau)$ where all $m$ eigenspaces are $n$-dimensional is denoted by $\mathcal{O}^{min}$. Recall further from the discussion in Section~\ref{sec:construct-lagrangians} that $\mathcal{O}^{min}\cap\mathcal{S}_{m,n}$ consists of only one point, which we here denote by $\mathbf{Y}$. As in the discussion surrounding the local isomorphism~\eqref{eq:local isom} from that section, consider the polydisk $P\subset \mathfrak{h}^{\pi}/W^{\pi}$ consisting of $(\vec{\mu},\vec{\lambda})$ such that
\begin{equation} \label{eq:tau-for-z} \vec{\mu}=(\lambda_1+z_1,\cdots,\lambda_m+z_m) \text { and }  \vec{\lambda}=(\lambda_1-(n-1)z_1,\cdots,\lambda_m-(n-1)z_m)\end{equation}
for $(z_1,\cdots,z_m)\in\CC^m$ in a small disk around $0$. Let $\Upsilon$ be the local isomorphism of Equation \eqref{eq:local isom}, which identifies a neighborhood  $\mathrm{nbd}(\mathbf{Y})\subset (\chi^{\pi})^{-1}(P)\cap\mathcal{S}_{m,n}$ with a neighborhood $\mathrm{nbd}(0)\subset\left(\mathfrak{sl}(n)^{((n-1)1)}\right)^{m}$.

Analogously to the conventions in Example~\ref{eg:Lag} for $m=1$, given $z=(z_1,\dots,z_m) \in \CC^m$, let $\cY_{z}$ indicate $\cY_{m,n,\tau}$ where $\tau$ is the configuration associated to $z$ under the coordinates~\eqref{eq:tau-for-z}. Likewise let $\eL_z$ be the product Lagrangian $\eL_{z_1} \times \dots \times \eL_{z_m}$ constructed via simultaneous collapse. Following the notation of the example let the slice of the $i$th factor $\mathfrak{sl}(n)^{((n-1)1)}$, which corresponds to the coordinate $z_i$, be $\cZ_i$ and its restriction to a particular value of $z_i$ be $\cY_{z_i}$. Above this factor we construct bundles $V_{i,\RR}$, $K_{i,\RR}$, and $T_{\RR}\cY_{z_i}$ as in the Example~\ref{example:realstructurem1}.

\begin{lemma}\label{l:isomorphic_real_bundle} For any $z=(z_1,\cdots,z_m)$ sufficiently close to $0$, the pullback bundle $\Upsilon^*T(\eL_{z})$ is stably isomorphic to the virtual bundle $T_{\RR}\cY_{z}|_{\Upsilon^{-1}(\eL_{z})}$.
	\end{lemma}

\begin{proof}
Recall that $\eL_z=\eL_{z_1}\times\cdots\times \eL_{z_m}$ has the structure of a product. By applying the conclusion of Example \ref{eg:Lag} to each factor, we see that \[T(\eL_{z_1}\times\cdots \times \eL_{z_m})\simeq T_{\RR}\cY_{z_1}|_{\eL_{z_1}}\times \dots \times T_{\RR}\cY_{z_m}|_{\eL_{z_m}}.\] What remains to be checked is that the local product of real structures on the righthand side agrees with the restriction of the global real structure on the tangent bundle. In particular, we will show
\begin{equation}\label{eq:real-structure-is-product}\Upsilon^*([V_{1,\RR}\oplus \cdots\oplus V_{m,\RR}]-[K_{1,\RR}\oplus\cdots\oplus K_{m,\RR}])|_{\Upsilon^{-1}(\eL_z)}=([V_{\RR}]-[K_{\RR}])|_{\Upsilon^{-1}(\eL_z)}.\end{equation} 
Toward this end, we first lift the local isomorphism $\Upsilon$ to the partial Grothendieck resolutions. In particular we construct a diagram of maps as follows.

\[
\begin{tikzcd}[row sep=large, column sep=1 cm]
(\widetilde{\chi}^{\pi})^{-1}(\widetilde{P})\cap \cS_{m,n}\arrow[r,"r_1"] \arrow[d, "\widetilde{\chi}^{\pi}"] \arrow[rrr, bend left=15, "\widetilde{\Upsilon}"]& 
(\chi^{\pi})^{-1}(P)\cap \cS_{m,n} \arrow[r,"\Upsilon"] \arrow[d, "\chi^{\pi}"] & 
\left(\mathfrak{sl}(n)^{((n-1)1)}\right)^{m}\arrow[d, "\chi"] & 
\left(\widetilde{\mathfrak{sl}}(n,\CC)^{((n-1)1)}\right)^{m}  \arrow[d, "\widetilde{\chi}"]  \arrow[l,"r_2"]  \\
\widetilde{P} \arrow[r] & 
P \arrow[r, "z"'] & 
\CC^m & 
\CC^m \arrow[l, "="']  
\end{tikzcd}
\]
In the above, $\widetilde{P}\subset \mathfrak{h}^{\pi}$ is the component of the preimage of $P$ under the map $\mathfrak{h}^{\pi} \to \mathfrak{h}^{\pi}/W^{\pi}$ which is specified by the order of the indices. Moreover the maps $\chi$ and $\widetilde{\chi}$ are $\chi = (\chi^{((n-1)1)},\cdots,\chi^{((n-1)1)})$ and $\widetilde{\chi}=(\widetilde{\chi}^{((n-1)1)},\cdots,\widetilde{\chi}^{((n-1)1)})$ respectively. 

Now suppose $Y\in\mathrm{nbd}(\mathbf{Y})$, and let $\chi^{\pi}(Y)$ be taken to $(z_1,z_2,\cdots,z_{m})$ under the coordinate map $z$ on the middle of the bottom row of the figure. Let $\Upsilon(Y)=(A_1,A_2,\cdots, A_m)$. We will see that it follows from the definition of $\Upsilon$ given in Section~\ref{sec:construct-lagrangians} that there are isomorphisms
\[ E_Y(\lambda_i+z_i) \simeq E_{A_i}(z_i) \qquad \qquad E_Y(\lambda_i - (n-1)z_i) \simeq E_{A_i}(-(n-1)z_i).\]
To see this, recall that the first map in the definition of $\Upsilon$ is a local isomorphism between $(\mathbf{Y}+\mathcal{R})\cap\frg^{\pi}$ and $\cS^{ss}\cap\frg^{\pi}$ that identifies $Y_A=\mathbf{Y}+A\in \cS^{ss}\cap\frg^{\pi}$ for small $A$ with $Y=g^{-1}Y_Ag$, where $g \in G$.
Thus, $g$ identifies the eigenspaces $E_{Y_A}(\lambda_i+z_i)\subseteq E_{\mathbf{Y}}(\lambda_i)$ with $E_{Y}(\lambda_i+z_i)$, and likewise identifies the eigenspaces $E_{Y_A}(\lambda_i-(n-1)z_i)\subseteq E_{\mathbf{Y}}(\lambda_i)$ with $E_{Y}(\lambda_i-(n-1)z_i)$. Having made this identification, we next see that projection $p:\CC^{mn}\to \CC^n$ on the first $n$ coordinates identifies $E_{Y_A}(\lambda_i+z_i)$  with $E_{A_i}(z_i)$ and  $E_{Y_A}(\lambda_i-(n-1)z_i)$ with $E_{A_i}(-(n-1)z_i)$, as promised.  

Note that if $z_i=0$, the logic above gives an identification of the generalized $\lambda_i$-eigenspace of $Y$ with the generalized $0$-eigenspace of $A_i$. Any choice of flag $F$ such that $(Y,F)\in (\widetilde{\chi}^{\pi})^{-1}(\widetilde{P})\cap \cS_{m,n}$ specifies an $(n-1)$-dimensional subspace $F_i$ of $E_{Y}(\lambda_i+z_i)$ for all $i$ and so correspondingly specifies an $(n-1)$-dimensional subspace $F_{i}'$ of $E_{A_i}(z_i)$. We define the map $\widetilde{\Upsilon}$ to be
\[
\widetilde{\Upsilon}(Y,F)=((A_1,F^1),\cdots, (A_m,F^m))\]
where $F^i=\{0\subset F_{i}'\subset  \CC^n\}$. Note that if all of the coordinate values  $z_1,\dots,z_m$ associated to $\chi^{\pi}(Y)$ are non-zero then $r_1^{-1}(Y)$ and $r_2^{-1}(\Upsilon(Y))$ both consist of a single point, and our definitions immediately give that $\widetilde{\Upsilon}(r_1^{-1}(Y))=r_2^{-1}(\Upsilon(Y))$.

Let $\widetilde{\eL}=r_2^{-1}(\cup_z\eL_z)$ in $\left(\widetilde{\mathfrak{sl}}(n,\CC)^{((n-1)1)}\right)^m$, which is a $(\mathbb{CP}^{n-1})^m$-bundle over $\CC^m$.  As discussed in Example \ref{eg:resolution-real-bundle} we have extensions $\widetilde{V}_{1,\RR}\oplus\cdots\oplus\widetilde{V}_{m,\RR}$ and $\widetilde{K}_{1,\RR}\oplus\cdots\oplus\widetilde{K}_{m,\RR}$ of the vector bundles $V_{1,\RR}\oplus \cdots\oplus V_{m,\RR}$ and $K_{1,\RR}\oplus\cdots\oplus K_{m,\RR}$, respectively, over $\widetilde{\eL}$. Let $\widetilde{\eL}'=\widetilde{\Upsilon}^{-1}(\widetilde{\eL}\cap r_{2}^{-1}(\mathrm{nbd}(0)))$. We define extensions $\widetilde{V}_{\RR}$ and $\widetilde{K}_{\RR}$ over $\widetilde{\eL}'$ as follows. First, let $\widetilde{V}$ and $\widetilde{V}_{\RR}\subset \widetilde{V}$ be trivial vector bundles over $\widetilde{\eL}'$ extending $V$ and $V_{\RR}$, respectively.

Next, note that at each point $(Y,F)$ if the $i$th coordinate $z_i\circ \chi^{\pi}(Y)$ specifying the configuration associated to $Y$ is nonzero, then $F_i=E_{Y}(\lambda_i+z_i)$ and $F_i'=E_{A_i}(z_i)$. Moreover since $(Y,F)\in\widetilde{\eL}'$, the matrices $A_i$ are normal, and so $E_{A_i}(-(n-1)z_i)$ is the orthogonal complement of $F_i'$.  Thus, we define $\overline{F}_i\subset \CC^{mn}$ to be equal $E_{Y}(\lambda_i-(n-1)z_i)$ if $z_i\neq 0$. Otherwise, consider the identification of $E_{Y}(\lambda_i)$ and $E_{A_i}(0)$ and define $\overline{F}_i$ to be the subspace corresponding to the orthogonal complement of $F_i'\subset E_{A_i}(0)$. 

Consider a basis $\{B_{ij}\mid 1\le i\le m,\  1\le j\le n\}$ for $\CC^{mn}$ such that $\{B_{i1},\cdots,B_{i(n-1)}\}$ and $\{B_{in}\}$ are bases for $F_i$ and $\overline{F}_i$, respectively, and the projection $p$ maps them to orthonormal bases for $p(F_i)$ and $p(\overline{F}_i)$.  Let $B=(B_{11}|B_{12}|\cdots|B_{mn})$. 
 We define $\widetilde{K}_{(Y,F),\RR}$ to consist of every $X\in \widetilde{V}_{(Y,F)}$ such that $B^{-1}Z_{X}B$ is block diagonal whose blocks alternate between $(n-1)\times (n-1)$ and $1\times 1$ Hermitian matrices. Since $B$ is chosen with respect to $F_i$ and $\overline{F}_i$ this is equivalent to asking that $[Y,Z_X]=0$, $Z_X$ preserves $F$, and $B^{-1}Z_{X}B$ is Hermitian. As illustrated by Example \ref{eg:resolution-real-bundle}, the spaces $\widetilde{K}_{(Y,F),\RR}$ form a vector bundle over $\widetilde{\eL}'$ which extends $K_{\RR}$.

At last we are ready to that prove~\eqref{eq:real-structure-is-product} holds. Note that $\widetilde{\eL}'$ is a smooth $(\mathbb{CP}^{n-1})^m$ bundle over a disk in $\widetilde{P}$, so in order to prove the claim for a small $z$ it is enough to check that over $\widetilde{\Upsilon}^{-1}(\widetilde{\eL}_0)$ there is an isomorphism between $\widetilde{V}$ and $\widetilde{\Upsilon}^*(\widetilde{V}_{1}\oplus\cdots\oplus\widetilde{V}_{m})$ that maps $\widetilde{K}_{\RR}$ to $\widetilde{\Upsilon}^{*}(\widetilde{K}_{1,\RR}\oplus\cdots\oplus\widetilde{K}_{m,\RR})$. 

Every point in $\widetilde{\Upsilon}^{-1}(\widetilde{\eL}_0)$ is of the form $(\mathbf{Y},F)$. To define such an isomorphism, we pick a basis for $E_{\lambda_i}(\mathbf{Y})$ that projects to the standard basis for $\CC^n$ by $p$. Then, we send every $X\in \widetilde{V}_{(\mathbf{Y},F)}$ to the action of $Z_X$ on $E_{\lambda_i}(\mathbf{Y})$ in this basis for $i=1,\cdots, n$. More precisely, let $M_{\lambda_i}$ be the $(mn)\times n$ matrix whose columns form a basis for $E_{\lambda_i}(\mathbf{Y})$ such that if $R_1 = (I_n \ 0 \dots 0)$ then we have $R_1M_{\lambda_i}=I_n$ is the identity matrix. We define 
\[\begin{split}\eta:\widetilde{V}|_{\widetilde{\Upsilon}^{-1}(\widetilde{\eL}_0)}\rightarrow \widetilde{\Upsilon}^{*}(\widetilde{V}_{1}\oplus\cdots\oplus\widetilde{V}_{m})|_{\widetilde{\Upsilon}^{-1}(\widetilde{\eL}_0)}&\\
	X\mapsto (R_1Z_XM_{\lambda_1},R_1Z_XM_{\lambda_2},\cdots,R_1Z_XM_{\lambda_m})&
\end{split}\]

Since, $[\mathbf{Y},Z_{X}]=0$ for any 
$X$, this is an isomorphism. Furthermore, we have that $X\in \widetilde{K}_{(\mathbf{Y},F),\RR}$ if and only if $Z_X$ preserves the flag $F$ and the action of $Z_{X}$ on every $E_{\lambda_i}(\mathbf{Y})$ in a basis $B_i=(B_{i1}|B_{i2}|\cdots|B_{in})$ constructed using the bases for $F_i$ and $\overline{F_i}$ as above is Hermitian. Note that $Z_X$ preserves the flag $F$ if and only if the $i$-th component of $\eta(X)$ preserves the flag $F^i$ for all $i$. Here, $\widetilde{\Upsilon}(\mathbf{Y},F)=((0,F^1),\cdots,(0,F^m))$. Moreover, $B_i=M_{\lambda_i}B_i'$ where $B_i'=R_1B_i$ and $B_i'$ is a unitary matrix.  Let $X'_i$ be the action of  $Z_X$ on $E_{\lambda_i}(\mathbf{Y})$ with respect to basis $B_i$, that is $Z_XB_i=B_iX_i'$. Then,
\[R_1Z_XB_i=R_1B_iX_i'=B'_iX'_i\ \ \Rightarrow (R_1Z_XM_{\lambda_i})B_i'=B_i'X_i'.\]
Since $B_i'$ is unitary, $R_1Z_XB_i$ is Hermitian if and only if $X_i'$ is Hermitian. Thus, $\eta$ gives an isomorphism between $\widetilde{K}_{(\mathbf{Y},F),\RR}$ and $\widetilde{K}_{(0,F^1),\RR}\oplus\cdots \widetilde{K}_{(0,F^m),\RR}$. This concludes the argument.

\end{proof}

\begin{lemma} \label{lem:isomorphism-is-isotopy}
For any crossingless matching $\mathfrak{m}$ with associated Lagrangian $\eL_m$, the tangent bundle $T\eL_{\mathfrak{m}}$ is stably isotopic to $(T_{\RR} \cY_{m,n,\tau})|_{\eL_{\mathfrak{m}}}$ as stable Lagrangian subbundles of $(T\cY_{m,n,\tau})|_{\eL_{\mathfrak{m}}}$. 
\end{lemma}

\begin{proof}
First we show that $T\eL_{\mathfrak{m}}$ is stably isomorphic to $(T_{\RR} \cY_{m,n,\tau})|_{\eL_{\mathfrak{m}}}$ as stable real vector bundles.
By Lemma \ref{l:isomorphic_real_bundle}, $T\eL_z$ is stably isomorphic to $(T_{\RR} \cY_{z})|_{\eL_z}$ where $\eL_z$ is the standard Lagrangian constructed using the isomorphism~\eqref{eq:local isom} for $z=(z_1,\dots,z_m)$ such that $z_i \neq 0$ for all $i$.
Let $\gamma \co [0,1] \to \frh^{\pi,reg}/W^{\pi}$ be a path such that $h_{\gamma}^{resc}(\eL_{\mathfrak{m}})$ is Hamiltonian isotopic to $\eL_z$.
Using this path, we have that $T(h_{\gamma}^{resc}(\eL_{\mathfrak{m}}))$ is stably isomorphic to $(T_{\RR} \cY_{z})|_{h_{\gamma}^{resc}(\eL_{\mathfrak{m}})}$. Since by Proposition \ref{p:K_real} the bundles $V$, $K$ and $K_\RR$ are well-defined on $
\cS_{m,n}^{\pi, reg} = (\chi^{\pi})^{-1}(\frh^{\pi,reg}/W^{\pi}) \cap \cS_{m,n}$, we see that $(h_{\gamma}^{resc})^*(T_{\RR} \cY_{z})|_{h_{\gamma}^{resc}(\eL_{\mathfrak{m}})}$ is 
isotopic to $T_{\RR} \cY_{m,n,\tau}|_{\eL_{\mathfrak{m}}}$ as stable real subbundles inside the full bundle $T(\cS_{m,n}^{\pi, reg})|_{\eL_{\mathfrak{m}}}$. Therefore, $T\eL_{\mathfrak{m}}$ is stably isomorphic to $(T_{\RR} \cY_{m,n,\tau})|_{\eL_{\mathfrak{m}}}$. 

We now argue that this is sufficient; that is, that two stable real subbundles of a complex bundle over a product of complex projective spaces are isotopic as Lagrangian subbundles after stabilization if and only if they are isomorphic as real vector bundles. Toward this end, we review a little bundle theory. For $N \geq 0$, let $U=\cup_N U(N)$, $O=\cup_N O(N)$ and $c \co O \to U$ be the complexification map. 
We have the induced long exact sequence 
\[
\dots \to U \to U/O \to BO \to BU \to B(U/O) \to \dots. 
\]
The bundle $(T\cY_{m,n,\tau})|_{\eL_{\mathfrak{m}}}$ defines a homotopy class in $[\eL_{\mathfrak{m}}, BU]$.
Let $r \co \eL_{\mathfrak{m}} \to BU$ be a representative of the homotopy class.
The stable Lagrangian subbundles $T\eL_{\mathfrak{m}}$ and $T_{\RR} \cY_{m,n,\tau}|_{\eL_{\mathfrak{m}}}$ define
maps $\eL_{\mathfrak{m}} \to BO$, denoted by $f_0$ and $f_1$ respectively, such that for $i=0,1$ there is a homotopy $H_i$ from $c \circ f_i$ to $r$.
Since $f_0$ and $f_1$ represent the same class in $[\eL_{\mathfrak{m}},BO]$, there is furthermore a homotopy $H_f$ from $f_0$ to $f_1$.
Therefore the concatenation of the homotopy $H_0$ with the paramatrization reversed followed by $c \circ H_f$ and then $H_1$ in total defines a based loop at $r$, and hence defines an element
$\gamma \in \pi_1(\operatorname{Map}(\eL_{\mathfrak{m}}, BU),r)$.
The two Lagrangian subbundles are isotopic through Lagrangian subbundles exactly when the class $\gamma$ is trivial.

Since $BU$ is an infinite loop space, it is group-like up to homotopy.
Therefore, we can consider the inverse $-r \co \eL_{\mathfrak{m}} \to BU$ and the homotopy equivalence that $-r$ induces between $\operatorname{Map}(\eL_{\mathfrak{m}}, BU)_r$, the component of $\operatorname{Map}(\eL_{\mathfrak{m}}, BU)$ containing $r$, to 
$\operatorname{Map}(\eL_{\mathfrak{m}}, BU)_e$, the component of $\operatorname{Map}(\eL_{\mathfrak{m}}, BU)$ containing a constant map $e$.
As a result, we have
\[\pi_1(\operatorname{Map}(\eL_{\mathfrak{m}}, BU),r)=\pi_1(\operatorname{Map}(\eL_{\mathfrak{m}}, BU),e)=[\eL_{\mathfrak{m}}, \Omega BU]=[\eL_{\mathfrak{m}}, U].\]
Since $\eL_{\mathfrak{m}}$ is diffeomorphic to the $m$-fold product of the complex projective space $\mathbb{CP}^{n-1}$ and \[[(\mathbb{CP}^{n-1})^m,U]=\widetilde{K}^{-1}((\mathbb{CP}^{n-1})^m)=0,\] we conclude that the class $\gamma$ is trivial. In other words, $T\eL_{\mathfrak{m}}$ is stably isotopic to $T_{\RR} \cY_{m,n,\tau}|_{\eL_{\mathfrak{m}}}$ as stable Lagrangian subbundles.

We also mention the following alternate argument using the same key facts. Fixing $r \co \eL_{\mathfrak{m}} \to BU$, the space of stable Lagrangian subbundles is the space of lifts $f \co \eL_{\mathfrak{m}} \to BO$ such that $c \circ f=r$ strictly. Since the complexification map $c:BO \to BU$ is compatible with the infinite loop space structures on $BO$ and $BU$, we can take the difference of two such lifts. This defines an element in $\operatorname{Map}(\eL_{\mathfrak{m}}, U/O)$.
 Since $[\eL_{\mathfrak{m}},U]=0$, the map $[\eL_{\mathfrak{m}}, U/O] \to [\eL_{\mathfrak{m}}, BO]$ is injective.
Therefore, two stable Lagrangian subbundles are isotopic if and only if they are isomorphic as stable real vector bundles. 
\end{proof}

\begin{corollary} \label{cor:stable-normal-real} Let $\iota_k$ be the intrinsic involution on $\cY_{m,n,\tau}$ defined in \eqref{eq:intrinsic}. Over every component of $\cY^{\fix}_{m,n,\tau}$, the tangent bundle $T(\eL^{\fix}_{\mathfrak{m}})$ is stably isotopic to $(T_{\RR} \cY_{m,n,\tau}^{\fix})|_{\eL^{\fix}_{\mathfrak{m}}}$ as stable Lagrangian subbundles of $(T\cY^{\fix}_{m,n,\tau})|_{\eL^{\fix}_{\mathfrak{m}}}$. Moreover, the normal bundle $N(\eL^{\fix}_{\mathfrak{m}})$ is stably isotopic to $(N_{\RR} \cY_{m,n,\tau}^{\fix})|_{\eL^{\fix}_{\mathfrak{m}}}$ as stable Lagrangian subbundles of $(N\cY^{\fix}_{m,n,\tau})|_{\eL^{\fix}_{\mathfrak{m}}}$.
\end{corollary}

\begin{proof} This follows immediately using the same logic as Lemma~\ref{lem:isomorphism-is-isotopy}, and the fact that rescaled parallel transport is $\iota_k$-equivariant from Corollary~\ref{cor:invariant-parallel}. Note in particular that the proof of Lemma~\ref{lem:isomorphism-is-isotopy} uses the fact that the Lagrangians are a product of complex projective spaces, and does not require that all factors of the product have the same complex dimension. It therefore applies to the connected components of the fixed Lagrangians, which are of the form $(\CP^{k-1})^s\times (\CP^{n-k-1})^{m-s}$ for some $0 \leq s \leq m$. \end{proof}

We finally conclude by collecting the statements above into the proof of Theorem~\ref{thm:stable_real}.

\begin{proof} [Proof of Theorem~\ref{thm:stable_real}] By Corollary~\ref{cor:tangent-real}, the tangent bundle $T\cY_{m,n,\tau}$ admits a stable real structure, and by Corollary~\ref{cor:fixed-real} the same is true of the tangent and normal bundles of the fixed set under the involution $\iota_k$. By Lemma~\ref{lem:isomorphism-is-isotopy} the real structure on the tangent bundle $T\cY_{m,n,\tau}$ is compatible with the Lagrangians and therefore in fact a stable real polarization, and similarly by Corollary~\ref{cor:stable-normal-real} the real structures on the tangent and normal bundles to the fixed sets are compatible with the Lagrangians and therefore stable real polarizations.
\end{proof}

\subsection{Proof of Theorem~\ref{thm:mainss}} \label{sec:assemble-proof-thmmain}

We are now ready to prove the existence of the localization spectral sequence from Theorem~\ref{thm:mainss}. 

\begin{proof} [Proof of Theorem \ref{thm:mainss}] Choose $0<k<n$, and consider the involution on $\mathcal Y_{m,n,\tau}$ given by $\iota_k$ as defined in \eqref{eq:intrinsic}. Let $K$ be expressed as the union of matchings $\mathfrak{m}_0$ and $\mathfrak{m}_{\beta}$, and consider the Lagrangians $\eL_0$ and $\eL_{\beta}$. By Proposition~\ref{p:fixedFloer}, the Floer cohomology of $\eL_0^{\fix}$ and $\eL_{\beta}^{\fix}$ in $\cY_{m,n,\tau}^{\fix}$ computes $H_{k,symp}(K) \oplus H_{n-k,symp}(K)$. By Corollary~\ref{cor:stable-normal-real}, the four-tuple $(\cY_{m,n,\tau}, \eL_0, \eL_{\beta}, \iota_k)$ admits a set of stable normal real polarizations; that is to say, there is a stable real polarization of $(N\cY_{m,n,\tau}^{\fix}, N\eL_0^{\fix}, N\eL_{\beta}^{\fix})$ where as usual we have a triple of bundles over each connected component of the fixed set for which the statement holds. Applying Theorem~\ref{thm:localization-intro} we obtain the claimed spectral sequence of Theorem~\ref{thm:mainss}. \end{proof}

\subsection{Remarks on the case of links} \label{sec:link-discussion}

We now touch on the case that $L$ is a link of more than one component rather than a knot. On the combinatorial side, there are analogs of the spectral sequences \eqref{eq:krssbig} and \eqref{eq:krsssplit} for links. These converge to a direct sum indexed over the sublinks $L'$ of $L$. For example, the analog of \eqref{eq:krsssplit} is of the form 
\begin{equation} \label{eq:combinatorial-links}
H_n(L; \mathbb C) \rightrightarrows \bigoplus_{L' \subset L} \left(H_{k}(L'; \mathbb C) \otimes H_{n-k}(L \setminus L'; \CC)\right)\end{equation}
where the homology of the empty sublink is understood to be the ground field $\CC$; one may find treatments of special cases in \cite{Gornik,Wu:generic} and the general case in \cite[Section 1.1]{RW:foams}. This generalizes the well-known behavior of the Lee spectral sequence for a link \cite{Leedeformation, Rasmussen:Milnor}.

In the symplectic case, we can apply the arguments of Section~\ref{sec:real} in the case of a link $L$ to obtain a spectral sequence which appears to parallel \eqref{eq:combinatorial-links}. However, the nature of the final page is somewhat complicated. Consider the generalization of Proposition~\ref{p:fixedFloer} to links. In this case, we may suppose that we are interested in $\eL_0$ constructed from the identity rainbow matching, with associated permutation the identity, and $\eL_{\beta}$ constructed from a matching $\mathfrak{m}_{\beta}$ whose associated permutation has $\ell$ cycles, where $\ell$ is the number of components of $L$. Rerunning the argument of the proposition now shows that $\eL_0^{\fix}$ and $\eL_{\beta}^{\fix}$ intersect in $2^{\ell}$ connected components of $\cY_{m,n,\tau}^{\fix}$. These $2^{\ell}$ are associated to the $2^{\ell}$ different sublinks of the link. The direct sum of the Lagrangian Floer groups in these  $2^{\ell}$ components is the final page of the localization spectral sequence. 

However, the extra components in which $\eL_0^{\fix}$ and $\eL_{\beta}^{\fix}$ intersect are not diffeomorphic to $\cY^{\pi_{std}}_{m,k,\tau}$ or $\cY^{\pi_{std}}_{m,n-k,\tau}$, but rather are products $\cY_{m,k,\tau_{\pi\setminus\pi'}}\times\cY_{m,n-k,\tau_{\pi'}}$ where neither $\pi'$ nor $\pi \setminus \pi'$ is standard. Therefore, Lagrangian Floer cohomology in these components is not necessarily equivalent to the usual symplectic Khovanov--Rozansky homology in the appropriate rank, but rather is another possibly different link invariant. We plan to investigate these link invariants in future work.

\begin{example} To emphasize the point above, we conclude by examining the example of the Hopf link $L$. Let $n=5, m=2, k=3$ as in Examples~\ref{example:n5m2} and \ref{example:n5m2unknot}. Let $L$ be constructed by letting $\mathfrak{m}_0$ be the identity rainbow matching shown in Figure~\ref{fig:matching_example} and $\mathfrak{m}_{\beta}$ be a small Hamiltonian translate of the same matching, so that $\mathfrak{m}_0$ and $\mathfrak{m}_{\beta}$ both have associated permutation the identity. Since the permutation determines which components of the fixed set of $\cY_{m,n,\tau}$ contain a component of the fixed Lagrangian, it follows that the components of $\eL_0^{\fix}$ and $\eL_{\beta}^{\fix}$ are located as follows.

\begin{enumerate}

\item The component of $\cY_{2,5,\tau}^{\fix}$ corresponding to $(\pi',f_1)$ where $\pi'=(2,2)$ is diffeormorphic to $\cY^{\pi_{std}}_{2,3,\tau}\times \{\mathrm{pt}\}$ and contains a component of $\eL_{0}^{\fix}$ and of $\eL_{\beta}^{\fix}$, each of which is diffeomorphic to $\CC \PP^2 \times \CC \PP^2$. 

\item The component of $\cY_{2,5,\tau}^{\fix}$ corresponding to $(\pi', f_2)$ where $\pi'=(1,1,1,1)$ is diffeomorphic to $\{\mathrm{pt}\} \times \cY^{\pi_{std}}_{2,2,\tau}$ and contains a component of $\eL_0^{\fix}$ and of $\eL_{\beta}^{\fix}$, each of which is diffeomorphic to $\CC \PP^1 \times \CC \PP^1$.

\item The component of the fixed set of $\cY_{2,5,\tau}^{\fix}$ corresponding to $(\pi',f_3)$ with $\pi'=(2,1,1)$ does not contain any components of the fixed Lagrangians.

\item The component of fixed set of $\cY_{2,5,\tau}^{\fix}$ corresponding to $(\pi',f_4)$ with $\pi'=(2,1,1)$ contains a component of $\eL_0^{\fix}$ and a component of $\eL_{\beta}^{\fix}$, each of which appears in $\cY_{2,3,\tau_{\pi\setminus\pi'}}\times\cY_{2,2,\tau_{\pi'}}$ as a product of a copy of $\CC \PP^{2}$ in the first factor with a copy of $\CC \PP^1$ in the second factor.

\item The component of fixed set of $\cY_{2,5,\tau}^{\fix}$ corresponding to $(\pi',f_5)$ with $\pi'=(2,1,1)$ contains a component of $\eL_0^{\fix}$ and a component of $\eL_{\beta}^{\fix}$, each of which appears in $\cY_{2,3,\tau_{\pi\setminus\pi'}}\times\cY_{2,2,\tau_{\pi'}}$ as a product of a copy of $\CC \PP^{2}$ in the first factor with a copy of $\CC \PP^1$ in the second factor.

\item The component of fixed set of $\cY_{2,5,\tau}^{\fix}$ corresponding to $(\pi',f_6)$ with $\pi'=(2,1,1)$ does not contain any components of the fixed Lagrangians.

\end{enumerate}

\end{example}

Finally, we mention one case in which the discussion above simplifies: if $n=2$ and $k=1$, each fixed Lagrangian is zero-dimensional with $2^m$ components. The set of intersections between $\eL_0^{\fix}$ and $\eL_{\beta}^{\fix}$ is exactly $2^{\ell}$ points, each of which is in a component of $\cY_{m,2,\tau}^{\fix}$ corresponding to a sublink of the link. The Lagrangian Floer cohomology in the fixed set is therefore exactly $\FF_2^{\oplus 2^{\ell}
}$, and we recover the spectral sequence~\eqref{eq:symplecticbn} of \cite[Theorem 7.16]{HLS:flexible} discussed in Remark~\ref{rem:Khovanov-case}.

\section{Extrinsic actions on symplectic Khovanov--Rozansky homology} \label{sec:extrinsic}

In this section we study an extrinsic action on the symplectic Khovanov--Rozansky homology of a two-periodic link. Recall that a link $L \subseteq S^3$ is said to be \emph{two-periodic} or \emph{doubly-periodic} if there is an orientation-preserving action of $\mathbb Z/2\mathbb Z$ on $S^3$ which preserves $L$ setwise and whose fixed set is an unknotted axis $A$ disjoint from the link. Taking the quotient under the action gives a quotient map $p \co (S^3, L) \rightarrow (S^3, \overline{L})$ which is a double branched covering map over $A$. The link $\overline{L} = p(L)$ is the \emph{quotient link}. Any two-periodic link admits a braid decomposition which is itself a second power; that is, it comes from stacking two copies of a braid $\beta$ vertically.

\subsection{The action for two-periodic links} \label{sec:periodic-action} Let $T$ be the $2mn\times 2mn$ block diagonal matrix $T=\mathrm{diag}(I_n,-I_n,I_n,\cdots,-I_n)$ and $\iota_p \co \mathfrak{g}\to \mathfrak{g}$ be the involution $\iota_p(Y)=-TYT$. This induces an involution on $\cS_{2m,n}$ given by  
\begin{align}\label{eq:action2}
\begin{pmatrix} 
	Y_{11} & I_n & &\\
	Y_{21} & 0 & I_n &&\\
	\dots & & &\dots   &\\
	Y_{(2m-1)1} & & & &I_n\\
	Y_{(2m)1}& 0 & 0 & 0&0
\end{pmatrix}  \mapsto \begin{pmatrix} 
	-Y_{11} & I_n & &\\
	Y_{21} & 0 & I_n &&\\
	\dots & & &\dots   &\\
	-Y_{(2m-1)1} & & & &I_n\\
	Y_{(2m)1}& 0 & 0 & 0&0
\end{pmatrix}.
\end{align}
Moreover, for any configuration $\tau=(\vec{\mu},\vec{\lambda})$, the involution $\iota_p$ maps $(\chi^{\pi})^{-1}(\tau)$ to $(\chi^{\pi})^{-1}(-\tau)$ and so $\cY_{2m,n,\tau}$ to $\cY_{2m,n,-\tau}$ where $-\tau=(-\vec{\mu},-\vec{\lambda})$.
	Thus, if $\tau$ is of the form $\tau=(\{\pm \mu_1, \dots, \pm \mu_m \},\{\pm \lambda_1, \dots \pm \lambda_m\})$ then there is an induced involution on $\cY_{2m,n,\tau}$. Moreover, the fixed set $\cY_{2m,n,\tau}^{\fix}$ consists of matrices $Y$ such that the block $Y_{2i-1,1}=0$ for all $1\le i\le m$.

\subsection{A stable tangent-normal isomorphism for the periodic case} \label{sec:periodic-tniso}

Let $f \co \cS_{m,n} \to \cS_{2m,n}$ be the map given by
\begin{equation}\label{eq:embed-fixed-set}
\begin{pmatrix} 
Y_{11} & I_n & &\\
Y_{21} & 0 & I_n &&\\
\dots & & &\dots   &\\
 Y_{(m-1)1} & & & &I_n\\
Y_{m1}& 0 & 0 & 0&0
\end{pmatrix}  \mapsto \begin{pmatrix} 
0 & I_n & &\\
Y_{11} & 0 & I_n &&\\
\dots & & &\dots   &\\
 0 & & & &I_n\\
Y_{m1}& 0 & 0 & 0&0
\end{pmatrix} 
\end{equation}
Furthermore let $\tau=(\{\pm \mu_1, \dots, \pm \mu_m \},\{\pm \lambda_1, \dots \pm \lambda_m\})$ and $\tau^2:=(\{ \mu_1^2, \dots,\mu_m^2\},\{ \lambda_1^2 \dots, \lambda_m^2\})$.

\begin{lemma}\label{l:fixed_locus}
The map $f$ restricts to an isomorphism $f \co \cY_{m,n,\tau^2} \to \cY_{2m,n,\tau}^{\fix}$.
\end{lemma}

\begin{proof}
Let $W \in \cY_{m,n,\tau^2}$. Suppose that $x \in \CC^{mn}$ is an eigenvector of $W$ with eigenvalue $\alpha=\beta^2$. Write $x$ in the form $x:=(x_1,\dots, x_m)^T$ with each $x_i \in \mathbb{C}^n$. Consider the vector $y \in \CC^{2mn}$ of the form
\[
y:=(y_1,y_2, \dots, y_{2m-1}, y_{2m})^T=(x_1, \beta x_1, x_2, \beta x_2, \dots, x_m, \beta x_m)^T.
\] 
One can check that $Wx=\alpha x$ implies that $f(W)y=\beta y$.

Since the $ \lambda_i^2$-eigenspace of $W$ is $1$-dimensional for all $i$ and the $  \mu_i^2$-eigenspace of $W$  is $(n-1)$-dimensional for all $i$, we conclude that the $ \pm \lambda_i$-eigenspace of $f(W)$ is $1$-dimensional for all $i$ and the $\pm  \mu_i$-eigenspace of  $f(W)$ is $(n-1)$-dimensional for all $i$.
Therefore, $f(W) \in \cY_{2m,n,\tau}^{\fix}$.

Conversely, if $Y \in \cY_{2m,n,\tau}^{\fix}$, let $W \in \cS_{m,n} $ be the matrix such that $f(W)=Y$. By comparing the eigenspaces as above, we conclude that $W \in \cY_{m,n,\tau^2}$.
\end{proof}

Applying the first claim of Corollary~\ref{cor:bundle-descriptions} to $\cY_{2m,n,\tau}$ we see that the tangent space of $Y \in \cY_{2m,n,\tau}$ is
\[
T_Y \cY_{2m,n,\tau}
= \{[Y,Z_X]\ |\ Z_X:=(Y^{2m-1}X|\cdots|YX|X),\ X\in M_{(2mn)\times n}\}.
\]

Now suppose that $Y \in \cY_{2m,n, \tau}^{\fix}$, and let $R=[Y,Z_X]$. We consider the behavior of $Z_X$. Observe that
\begin{equation}\label{eq:sign-change}-TRT = -T[Y,Z_X]T = [-TYT, TZ_XT] = [Y, TZ_XT].\end{equation} This means that the map that sends $R$ to $-TRT$ corresponds to the map that sends $Z_X$ to $TZ_XT$, with no additional change of sign. We may split the tangent space at $Y$ as the direct sum $T_Y\cY_{2m,n,\tau}=V\oplus W$ where
\[V=T_{Y}\cY_{2m,n,\tau}\cap\{[Y,Z_X]:TZ_XT=Z_X\} \qquad \qquad W=T_{Y}\cY_{2m,n,\tau}\cap\{[Y,Z_X]:TZ_XT=-Z_X\}.\]

\noindent This brings us to the following analysis of the tangent and normal bundles to the fixed set.

\begin{lemma} \label{lem:tangent-normal}
	For any $Y\in \cY_{2m,n,\tau}^{\fix}$, the tangent and normal spaces of $\cY_{2m,n,\tau}^{\fix}$ at $Y$ are the subspaces 
	\[\begin{split}
	&T_Y\cY^{\fix}_{2m,n,\tau}=V=\{[Y,Z_X]\ |\ Z_X:=(Y^{2m-1}X|\cdots|YX|X),\ X\in M_{(2mn)\times n}\ \text{and}\ \ TZ_XT=Z_X\}\\
	&N_Y\cY^{\fix}_{2m,n,\tau}=W=\{[Y,Z_X]\ |\ Z_X:=(Y^{2m-1}X|\cdots|YX|X),\ X\in M_{(2mn)\times n}\ \text{and}\ \ TZ_XT=-Z_X\}
\end{split}\]
	 and for $R = [Y,Z_X]$ in $T_Y\cY^{\fix}_{2m,n,\tau}$ the map $i(R)=YR$ induces an isomorphism of the form 
	\[
	i:
	\begin{pmatrix} 
		0 & 0& &\\
		R_{21} & 0 & 0 &&\\
		\dots & & &\dots   &\\
		0 & & & &0\\
		R_{2m,1}& 0 & 0 & 0&0
	\end{pmatrix}  \mapsto \begin{pmatrix} 
		R_{21} & 0 & &\\
		0 & 0 & 0 &&\\
		\dots & & &\dots   &\\
		R_{2m,1} & & & &0\\
		0& 0 & 0 & 0&0
	\end{pmatrix}
	\]
	 from the tangent bundle to the normal bundle.
\end{lemma}

\begin{proof}
	Let  $I_Y=d_Y\iota_p$ on $T_Y\cY_{2m,n,\tau}$, so that for $R=[Y,Z_X]\in T_Y\cY_{2m,n,\tau}$ we have $I_Y(R)=-TRT$. The tangent space and normal space of  $\cY_{2m,n,\tau}^{\fix}$ at $Y$ are the $+1$ and $-1$ eigenspaces of $I_Y$ respectively.  Recalling from \eqref{eq:sign-change} that $-T[Y,Z_X]T=[Y,TZ_XT]$, we see that we have $V\subset E_{+1}(I_Y)$. Similarly, $W\subset E_{-1}(I_Y)$.  Therefore, the tangent and normal spaces of $\cY_{2m,n,\tau}^{\fix}$ at $Y$ are equal to $V$ and $W$ respectively. 
	
Now consider $R=[Y,Z_X]$ in $T_Y\cY_{2m,n,\tau}$. We recall that $R$ is of the form $(R_{ij})_{1\leq i,j\leq 2m}$ where the $R_{ij}$ are $n\times n$ matrices with the property that $R_{ij}=0$ if $j\neq 1$. Indeed, we see that $TRT=-R$ if and only if the blocks $R_{(2i-1)1}$ are also zero for all $i$. So if $R \in V$, we see that $R$ must have the form 
	\[R=\begin{pmatrix} 
		0 & 0& &\\
		R_{21} & 0 & 0 &&\\
		\dots & & &\dots   &\\
		0 & & & &0\\
		R_{(2m)1}& 0 & 0 & 0&0
	\end{pmatrix}. \]
Similarly if $R \in W$, then $R$ has the form	
		\[R=\begin{pmatrix} 
		R_{11} & 0& &\\
		0 & 0 & 0 &&\\
		\dots & & &\dots   &\\
		R_{(2m-1)1} & & & &0\\
		0& 0 & 0 & 0&0
	\end{pmatrix}. \]	
For the isomorphism, suppose $R=[Y,Z_X]$ is in $V$. It follows that $YR=Y[Y,Z_X]=[Y,YZ_X]$ is in $W$, because $T(YZ_X)T=TYTTZ_XT=-YZ_X$. Indeed, we note that $R$ is determined by the last column of $Z_X$, that is $X$,  and the map $Z_X\mapsto YZ_X=Z_{YX}$ is an isomorphism on the last block column. 	
	\end{proof}


 

\begin{lemma}\label{lem:localstr-2} 
	Consider the polydisk $P\subset \mathfrak{h}^{\pi}/W^{\pi}$ consisting of points 
	\begin{align*}
\begin{pmatrix}\{\lambda_1+z_1,-\lambda_1+w_1,\cdots, \lambda_{\ell}+z_{\ell},-\lambda_{\ell}+w_{\ell},\pm\mu_{\ell+1},\cdots,\pm\mu_{m}\}, \\
\{\lambda_1-(n-1)z_1,-\lambda_1-(n-1)w_1,\cdots,\lambda_{\ell}-(n-1)z_{\ell},-\lambda_{\ell}-(n-1)w_{\ell},\pm\lambda_{\ell+1},\cdots,\pm\lambda_{m}\}
\end{pmatrix}
\end{align*}
	 where $\lambda_i,\mu_j$ pairwise distinct and $|z_i|,|w_i|\le \epsilon$ for some small $\epsilon$. Let $\mathcal{O}^{min}\subset (\chi^{\pi})^{-1}(P)$ be the orbit consisting of the matrices whose $\pm\lambda_i$-eigenspaces are $n$-dimensional for all $i=1, \cdots,\ell$. Then the isomorphism \eqref{eq:local isom}
	\[(\chi^{\pi})^{-1}(P)\cap \mathcal S_{2m,n} \to (\mathcal{O}^{min}\cap \mathcal S_{2m,n})\times (\frsl(n)^{((n-1)1)})^{2\ell}\]
	can be made equivariant with respect to $\iota_p$ on the left and the involution on the right given by 
\begin{equation}\label{eq:local-involution}(Y,A_1,A_1',\cdots,A_{\ell},A_{\ell}')\to (\iota_p(Y),-A_1',-A_1, \cdots,-A_{\ell}',-A_{\ell}).\end{equation}
	\end{lemma}	

\begin{proof}
The involution $\iota_p$ preserves $\mathcal{O}^{min}\cap \mathcal S_{2m,n}$ and $\mathcal S_{2m,n}$ setwise. Therefore the bundle $T\mathcal S_{2m,n}|_{\mathcal{O}^{min}\cap \mathcal S_{2m,n}}$ is an $\iota_p$-equivariant holomorphic bundle over $\mathcal{O}^{min}\cap \mathcal S_{2m,n}$ and $T(\mathcal{O}^{min}\cap \mathcal S_{2m,n})$ is an $\iota_p$-equivariant subbundle. Thus by Lemma \ref{lem:equiv-normal} there is holomorphic $\iota_p$-equivariant subbundle $\cR$ of  $T\mathcal S_{2m,n}|_{\mathcal{O}^{min}\cap \mathcal S_{2m,n}}$  complementary to $T(\mathcal{O}^{min}\cap \mathcal S_{2m,n})$. 
Notice that $\iota_p$ acts on $\mathcal{S}^{ss}\cap\mathfrak{g}^{\pi} $ and the canonical isomorphism $\cR\cap \mathfrak{g}^{\pi} \simeq \mathcal{S}^{ss}\cap\mathfrak{g}^{\pi} $
is $\iota_p$-equivariant.

We now recall the second step of the isomorphism, which relates $\mathcal{S}^{ss}\cap\mathfrak{g}^{\pi} $ to  $(\mathcal{O}^{min}\cap \mathcal S_{2m,n})\times (\frsl(n)^{((n-1)1)})^{2\ell}$. Given $Y+A \in \mathcal{S}^{ss}\cap\mathfrak{g}^{\pi}$, recall that we consider the action of $A$ on the eigenspaces $E_Y( \pm \lambda_i)$ for $i=1,\dots,\ell$. The eigenspaces $E_{Y}(\pm\lambda_i)$ can each be canonically identified with $\mathbb{C}^n$ by projection $p$ on the first $n$ coordinates using Lemma~\ref{lem:injectivity} as usual. With this in mind, we obtain the canonical isomorphisms
\[\frsl(E_Y(\pm\lambda_i))^{((n-1)1)} \cong \frsl(n)^{((n-1)1)}. \]
The matrix $Y+A$ is sent to $(Y, A|_{E_Y(\lambda_1)},A|_{E_Y(-\lambda_1)}, \dots, A|_{E_Y(\lambda_{\ell})}, A|_{E_Y(-\lambda_{\ell})})$ under these canonical isomorphisms.

Note that if $v\in E_Y(\pm\lambda_i)$ then $Tv\in E_{Y}(\mp\lambda_i)$. Consider a basis $B=\{v_1,v_2,\cdots, v_n\}$ for $E_Y(\lambda_i)$ such that the projection of $v_j$ on the first $n$ coordinates is equal to $e_j$. We see that $TB=\{Tv_1,\cdots, Tv_n\}$ will be a basis for $E_Y(-\lambda_i)$ with the same property. Once we have chosen the bases in this way, an element $A$ in  $\frsl(E_Y(\lambda_i))^{((n-1)1)} \times \frsl(E_Y(-\lambda_i))^{((n-1)1)}$ in the basis $\{B,TB\}$ for $E_Y(\lambda_i)\oplus E_Y(-\lambda_i)$ appears in a block diagonal form $\diag(A_{i,+},A_{i,-})$ where $A_{i,\pm}\in \frsl(n)^{((n-1)1)}$ and the above identification maps $A$ to $(A_{1,+},A_{1,-}, \dots, A_{\ell,+},A_{\ell,-})$. Note that $\iota_p(Y+A)=\iota_p(Y)+\iota_p(A)=\iota_p(Y)-TAT$. Observe that under the canonical isomorphisms, $TAT$ is mapped to the presentation of $A$ in the basis $\{TB,B\}$ which is the block diagonal matrix $\mathrm{diag}(A_{i,-},A_{i,+})$. It follows that $-TAT$ is identified with $(-A_{i,-},-A_{i,+})$. \end{proof}

\begin{corollary}\label{c:preserving Lag} 
Let $\mathfrak{m}=\{\pm\delta_1,\cdots,\pm\delta_{m}\}$ be a matching with endpoints on $\tau$ such that endpoints of $\delta_i$ are $\lambda_i$ and $\mu_i$.
Then $\iota_p$ preserves $\eL_{\mathfrak{m}}$ and there is a $\ZZ/2\ZZ$-equivariant diffeomorphism $\phi$ between $\eL_{\mathfrak{m}}$ and $(\mathbb{CP}^{n-1}\times \mathbb{CP}^{n-1})^m$ equipped with the $\ZZ/2\ZZ$ action 
\begin{align}\label{eq:i_pLag}(z_1,z_1',\dots,z_m,z_m') \mapsto (z_1',z_1, \dots, z_m',z_m).\end{align}

\end{corollary}

\begin{proof}
This is analogous to Lemma \ref{lem:lagrangian-involution}. We can show that the parallel transport is equivariant so it reduces to the case that $\mathfrak{m}$ is as in Lemma \ref{l:productLag} for which by Lemma \ref{lem:localstr-2}, $\eL_{\mathfrak{m}}$  can be equivariantly identified up to equivariant Hamiltonian isotopy with the product of the vanishing projective spaces in $(\frsl(n)^{((n-1)1)})^{2m}$ with the action \eqref{eq:i_pLag}.

To identify $\eL_{\mathfrak{m}}^{\fix}$ with $f(\eL_{\mathfrak{m}^2})$, it suffices to work near $(\mathcal{O}^{min}\cap \mathcal S_{2m,n})$ so that the local isomorphism \eqref{eq:local isom} with $\ell=m$ applies.
We consider the composition
\begin{align}\label{eq:embed}
(\frsl(n)^{((n-1)1)})^{m} \to (\chi^{\bar{\pi}})^{-1}(\bar{P}) \cap \mathcal S_{m,n} \to (\chi^{\pi})^{-1}(P) \cap  \mathcal S_{2m,n} \to (\frsl(n)^{((n-1)1)})^{2m} 
\end{align}
where the first and last arrows are the respective local isomorphisms, the middle arrow is $f$, $\bar{P}$ consists of points $(\{\lambda_1^2+z_1,\dots, \lambda_m^2+z_m\} ,\{\lambda_1^2-(n-1)z_1,\dots, \lambda_m^2-(n-1)z_m\})$ and $\bar{\pi}$ is the corresponding partition.
Note that $f$ covers the map 
\[
\bar{P} \to P, \quad
(\{\lambda_i^2+z_i\}, \{\lambda_i^2-(n-1)z_i\}) \to \left(\left\{\pm \sqrt{\lambda_i^2+z_i}\right\}, \left\{\pm \sqrt{\lambda_i^2-(n-1)z_i}\right\}\right).
\]
Since the last local isomorphism is equivariant and the image of $f$ is the fixed locus of $\iota_p$, we know that the image of the composition \eqref{eq:embed} is the fixed antidiagonal
\[
\Delta:=\{(A_1,-A_1,\dots,A_m,-A_m) \}
\]

The Lagrangians $\eL_{\mathfrak{m}^2}$ and $\eL_{\mathfrak{m}}$ are identified up to Hamiltonian isotopy with the products of vanishing projective spaces in 
$(\frsl(n)^{((n-1)1)})^{m}$ and $(\frsl(n)^{((n-1)1)})^{2m}$ respectively.
It is easy to check that $\eL_{\mathfrak{m}}^{\fix}$ is Hamiltonian isotopic to the product of vanishing projective spaces in $\Delta$ with respect to the product K\"ahler form induced from $(\frsl(n)^{((n-1)1)})^{2m}$.
Since the construction of vanishing projective space is well-defined up to Hamiltonian isotopy,  $f(\eL_{\mathfrak{m}^2})$ is Hamiltonian isotopic to $\eL_{\mathfrak{m}}^{\fix}$ in $\Delta$ and hence in $\cY_{2m,n,\tau}^{\fix}$.

\end{proof}


We are now ready to show that the isomorphism $i$ between the tangent and normal bundles of $\cY^{\fix}_{2m,n,\tau}$ from Lemma \ref{lem:tangent-normal} is compatible with the tangent and normal bundles of the Lagrangians; that is, $i$ determines an isomorphism of polarizations.

\begin{lemma} \label{lem:lagrangian-tangent-normal}
For any equivariant crossingless matching $\mathfrak{m}$, that is invariant under $z\mapsto -z$ as in Corollary \ref{c:preserving Lag}, with associated Lagrangian $\eL_{\mathfrak{m}}$, the normal bundle $N\eL_{\mathfrak{m}}^{\fix}$ is stably isotopic to $i(T\eL_{\mathfrak{m}}^{\fix})$ as stable Lagrangian subbundles of $(N\cY_{2m,n,\tau}^{\fix})|_{\eL^{\fix}_{\mathfrak{m}}}$. 
\end{lemma}

\begin{proof}
Given Corollary \ref{c:preserving Lag} and the fact that the tangent and normal bundles of the diagonal submanifold of a product of a manifold with itself are isomorphic, it is clear that $N\eL_{\mathfrak{m}}^{\fix}$ and $T\eL_{\mathfrak{m}}^{\fix}$ are isomorphic as real vector bundles. 
Therefore, we can apply the same argument in Lemma \ref{lem:isomorphism-is-isotopy} to show that $N\eL_{\mathfrak{m}}^{\fix}$ and $i(T\eL_{\mathfrak{m}}^{\fix})$ are isotopic as  stable Lagrangian subbundles.
\end{proof}

\subsection{Proof of Theorem \ref{thm:periodicss}} \label{sec:proof-periodic}

\begin{figure}
\scalebox{.8}{
\small
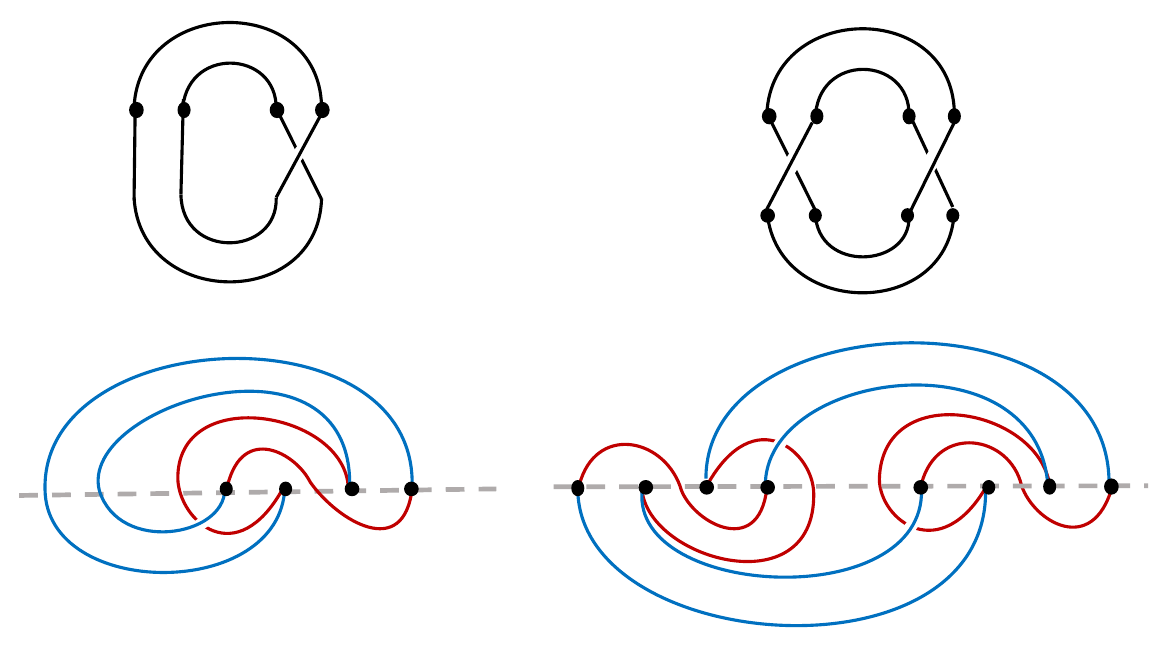}
\caption{An example of matchings for a periodic link $L$ and its quotient $\overline{L}$. On the left we see the quotient link $\overline{L}$ with blue identity matching $\mathfrak{m}_0$ and red braid matching $\mathfrak{m}_{\beta}$; these are the same as the matchings in Figure~\ref{fig:matching_example} after a rotation. On the right we see lifts of these two matchings. The blue overcrossing matching is $\mathfrak{m}_0'$ and the red undercrossing matching is $\mathfrak{m}_{\beta^2}$.}
\label{fig:periodic_matchings_example}
\end{figure}

We conclude with a proof of Theorem~\ref{thm:periodicss}. Let $L$ be a two-periodic link, and choose a braid representative $\beta^2$ on $m$ strands for $L$ which is the result of stacking two copies of a braid $\beta$. Let $\overline{L}$ be the quotient link of $L$ under the action, such that $\beta$ is a braid representative for $\overline{L}$. Consider a pair of configurations
\begin{align*} 
\tau &= ( \{\pm \mu_1, \dots, \pm \mu_m \}, \{\pm \lambda_1, \dots, \pm \lambda_m\}) \\
\tau^2 &= (\{\mu_1^2,\dots, \mu_m^2\}, \{\lambda_1^2, \dots, \lambda_m^2\})\end{align*}
with the property that all of the $\mu_i$ and $\lambda_j$ are real numbers and furthermore $0<\lambda_1 < \lambda_2 < \dots <\lambda_m < \mu_1 < \dots < \mu_m$. These configurations, and the matchings we shortly describe, are illustrated in Figure~\ref{fig:periodic_matchings_example}.

Because of the ordering we have chosen on the variables $\lambda_i, \mu_j$ and hence on the variables $\lambda_i^2, \mu_j^2$, we observe that the identity rainbow matching $\mathfrak{m}_0$ appears as the blue set of paths curved around the origin shown on the lefthand side of Figure~\ref{fig:periodic_matchings_example}. Choose a matching $\mathfrak{m}_{\beta} = (\epsilon_1, \dots, \epsilon_m)$ on $\tau^2$ such that pairing $\mathfrak{m}_{\beta}$ and $\mathfrak{m}_0$ on $\tau^2$ gives $\overline{L}$ with the property that the paths $\epsilon_i$ lie entirely in the right-hand side of the plane; that is, such that the real coordinate of a point $\epsilon_i(t)$ is always strictly positive. We illustrate this choice in the red matching on the lefthand side of the figure.

Now we consider matchings on $\tau$. Let $\mathfrak{m}_0'$ be the blue matching shown on the righthand side of Figure \ref{fig:periodic_matchings_example} which connects $\mu_i$ with $-\lambda_i$ by a rainbow arc above the real axis and connects $-\mu_i$ with $\lambda_i$ by a rainbow arc below the real axis. Observe that $(\mathfrak{m}_0')^2=\mathfrak{m}_0$. Since every path $\epsilon_i$ in $\mathfrak{m}_{\beta}$ lies entirely in the right-hand side of the plane, there are two distinct connected paths squaring to $\epsilon_i$ one of which connects the positive square roots of $\epsilon_i(0)$ and $\epsilon_i(1)$ and the other of which connects the negative square roots. Let these paths be $\delta_i$ and $-\delta_i$ respectively. Set $\mathfrak{m}_{\beta^2} = (\pm \delta_1, \dots, \pm \delta_m)$; this is illustrated in the red matching on the righthand side of the figure. We see that $\mathfrak{m}_{\beta^2}$ and $\mathfrak{m}_0'$ pair to give $L$. Moreover by construction $\mathfrak{m}_{\beta^2}^2=\mathfrak{m}_{\beta}$.

\begin{proof}[Proof of Theorem \ref{thm:periodicss}] Let $L$ be a two-periodic link and $\overline{L}$ its quotient. Choose matchings $\mathfrak{m}_0$ and $\mathfrak{m}_{\beta}$ for $\overline{L}$ and $\mathfrak{m}_0'$ and $\mathfrak{m}_{\beta^2}$ with $2m$ strands for $L$ as described above. Consider the triple $(\cY_{2m,n,\tau}, \eL_{\mathfrak{m}_0'}, \eL_{\mathfrak{m}_{\beta^2}})$, whose Floer cohomology computes $H_{n, symp}(L)$. 

By Lemma~\ref{l:fixed_locus} and Corollary~\ref{c:preserving Lag}, the fixed sets $(\cY_{2m,n,\tau}^{\fix}, \eL_{\mathfrak{m}_0'}^{\fix}, \eL_{\mathfrak{m}_{\beta^2}}^{\fix})$ of the involution $\iota_p$ are the image under the isomorphism $f$ of \eqref{eq:embed-fixed-set} of the triple $(\cY_{m,n,\tau^2}, \eL_{\mathfrak{m}_0}, \eL_{\mathfrak{m}_\beta})$, the Floer cohomology of which computes $H_{n, symp}(\overline{L})$. By Lemma~\ref{lem:tangent-normal} there is an isomorphism between the tangent bundle $T\cY_{2m,n,\tau}^{\fix}$ and the normal bundle $N\cY_{2m,n,\tau}^{\fix}$. By Lemma~\ref{lem:lagrangian-tangent-normal} this isomorphism is compatible with the Lagrangians and gives a stable isomorphism between the tangent and normal polarization. We may now apply Large's localization theorem, which is to say Theorem~\ref{thm:large-localization}, to obtain the claimed spectral sequence.
\end{proof}

\begin{remark} The reader may be aware that previous studies of two-periodic links in symplectic Khovanov homology have often used a picture which starts from a bridge diagram for the two-periodic link $L$ as in for example \cite[Figure 5.1]{HMR}; one uses the endpoints of the bridges as the points of the configuration $\tau$, the set of bridges as the identity matching, and the set of undercrossing arcs as the symmetric matching determining $L$. We have avoided this picture in favor of a more braid-centric discussion, since the $\frsl(n)$ symplectic link homologies currently take as input a braid representative of the link (or flattened braid diagram), not an arbitrary bridge diagram. \end{remark}

\bibliographystyle{amsalpha}
\bibliography{biblio}

\end{document}